\documentclass[11pt]{article}

\usepackage{amsmath, amsthm, amssymb, amscd, paralist}
\usepackage{enumerate}
\usepackage{enumitem}
\usepackage{graphicx}
\usepackage[export]{adjustbox}
\usepackage[dvipsnames]{xcolor}
\usepackage{bm}
\usepackage[title]{appendix}
\usepackage{hyperref}
\usepackage{tikz}
\usepackage{wrapfig}
\usetikzlibrary{matrix,arrows,positioning}
\definecolor{Green}{rgb}{0.0, 0.5, 0.0}

\newcommand{\eps}{\varepsilon}

\newcommand{\diam}{\operatorname{diam}}

\newcommand{\wt}{\widetilde}

\newcommand{\Var}{\operatorname{Var}}

\newcommand{\E}{\mathbb{E}}

\renewcommand{\P}{\mathbb{P}}

\newcommand{\Unif}{\operatorname{Unif}}

\newcommand{\KL}{\operatorname{KL}}

\newcommand{\B}{\mathbb{B}}

\newcommand{\R}{\mathbb{R}}
\newcommand{\N}{\mathbb{N}}

\newcommand{\bt}{\mathbf{t}}

\newcommand{\bd}{\mathbf{d}}
\newcommand{\bk}{\mathbf{k}}

\newcommand{\mA}{\mathcal{A}}
\newcommand{\mB}{\mathcal{B}}
\newcommand{\mC}{\mathcal{C}}
\newcommand{\mD}{\mathcal{D}}

\newcommand{\mF}{\mathcal{F}}

\newcommand{\mH}{\mathcal{H}}

\newcommand{\mL}{\mathcal{L}}
\newcommand{\mM}{\mathcal{M}}
\newcommand{\mN}{\mathcal{N}}

\newcommand{\mS}{\mathcal{S}}

\newcommand{\mW}{\mathcal{W}}

\newcommand{\Ind}{\mathbf{1}}

\newcommand{\bbB}{\mathbb{B}}
\newcommand{\bbH}{\mathbb{H}}

\newcommand{\bu}{\mathbf{u}}
\newcommand{\bv}{\mathbf{v}}

\newcommand{\bx}{\mathbf{x}}
\newcommand{\bX}{\mathbf{X}}
\newcommand{\by}{\mathbf{y}}
\newcommand{\bY}{\mathbf{Y}}

\newcommand{\balpha}{\bm{\alpha}}
\newcommand{\bbeta}{\bm{\beta}}

\newcommand{\bxi}{\bm{\xi}}

\newcommand{\frakr}{{\mathfrak r}}

\numberwithin{equation}{section}
\newtheorem{theorem}{Theorem}
\newtheorem{lemma}{Lemma}
\newtheorem{remark}{Remark}

\newtheorem{assump}{Assumption}

\newcommand{\footremember}[2]{%
   \footnote{#2}
    \newcounter{#1}
    \setcounter{#1}{\value{footnote}}%
}

\begin{document}

% "Title of the paper"

\title{Posterior contraction for deep Gaussian process priors}

\author{%
    Gianluca Finocchio\footremember{gf}{University of Vienna}%
    \ and Johannes Schmidt-Hieber\footremember{jsh}{University of Twente}%
}

\date{\today}
\maketitle

\begin{abstract}
We study posterior contraction rates for a class of deep Gaussian process priors applied to the nonparametric regression problem under a general composition assumption on the regression function. It is shown that the contraction rates can achieve the minimax convergence rate (up to $\log n$ factors), while being adaptive to the underlying structure and smoothness of the target function. The proposed framework extends the Bayesian nonparametric theory for Gaussian process priors.
\end{abstract}

\paragraph{Keywords:} Bayesian inference; nonparametric regression; contraction rates; deep Gaussian processes; uncertainty quantification; neural networks.

\textbf{MSC 2020:} Primary: 62G08, 62G20; Secondary: 62C20, 62R07.
% MSC 2020: 62-XX Statistics
%
%           62Cxx Statistical decision theory
%           62C10 Bayesian problems; characterization of Bayes procedures
% (Second.) 62C20 Minimax procedures in statistical decision theory
%
%           62Gxx Nonparametric inference
%           62G05 Nonparametric estimation
% (Primary) 62G08 Nonparametric regression and quantile regression
% (Primary) 62G20 Asymptotic properties of nonparametric inference
%           
%           62Rxx Statistics on algebraic and topological structures
% (Second.) 62R07 Statistical aspects of big data and data science
%

\section{Introduction}

In the multivariate nonparametric regression model with random design distribution $\mu$ supported on $[-1,1]^d,$ we observe $n$ i.i.d. pairs $(\bX_i,Y_i) \in[-1,1]^d\times \R,$ $i=1,\dots,n,$ with $\bX_i\sim\mu$, 
\begin{align}\label{eq.reg_model}
    Y_i = f^*(\bX_i) + \eps_i, \quad i=1,\dots,n
\end{align}
and $\eps_i$ independent standard normal random variables that are independent of the design vectors $(\bX_1,\dots,\bX_n)$. We aim to recover the true regression function $f^*:[-1,1]^d\to\R$ from the sample. Here it is assumed that the regression function $f^*$ itself is a composition of a number of unknown simpler functions. This comprises several important cases including (generalized) additive models. As proved in the recent work \cite{2022arXiv220507764G}, Gaussian process priors are suboptimal to learn the regression function. Meanwhile, \cite{SH19} proved that sparsely connected deep neural networks are able to pick up the underlying composition structure and achieve near-minimax $L^2(\mu)$-estimation rates over a compositional function class.

Deep Gaussian process priors (DGPs), cf. \cite{Neal96, damianou2013deep, CutajarDGP16}, can be viewed as a Bayesian analogue of deep networks. While deep nets are built on a hierarchy of individual network layers, DGPs are based on iterations of Gaussian processes. Compared to neural networks, DGPs have moreover the advantage that the posterior can be used for uncertainty quantification. This makes them potentially attractive for AI applications with a strong safety aspect, such as automated driving and health.

In Bayesian nonparametric regression without compositional constraints, Gaussian process priors are a natural choice and a comprehensive literature is available, see for instance Chapter~11 in~\cite{GVDV17}. The seminal contribution of~\cite{vvvz} fully characterizes the contraction rate of the posterior distribution for a Gaussian process prior by the small-ball probabilities of the prior and a property of the associated reproducing kernel Hilbert space (RKHS). From this, one can deduce that the best possible rates are achieved when the smoothness of the prior matches that of the unknown regression function. Adaptation to the smoothness can be achieved via rescaled smooth Gaussian processes and an inverse gamma hyperprior on the rescaling parameter as in~\cite{vvv2009adaptive}.

In this work we extend the theory of Gaussian process priors to derive posterior contraction rates for DGPs to learn compositional functions. Inspired by model selection priors, we propose a hierarchical prior construction. First, we assign prior weights to possible composition structures. Then, given a composition structure, the final DGP prior puts suitable Gaussian processes on all functions in this model.

The main contribution of this work is to provide an explicit construction for a large class of deep Gaussian process priors and to characterize the corresponding posterior contraction rates. By controlling the composition of Gaussian processes, we obtain tight bounds for the entropy and the decentered small-ball probabilities of DGPs. We also provide examples of DGP priors inducing nearly minimax posterior contraction rates. In particular, if there is some low-dimensional structure in the composition, the posterior will not suffer from the curse of dimensionality.

\begin{figure}[ht]
    \centering
    \includegraphics[height=0.35\textwidth]{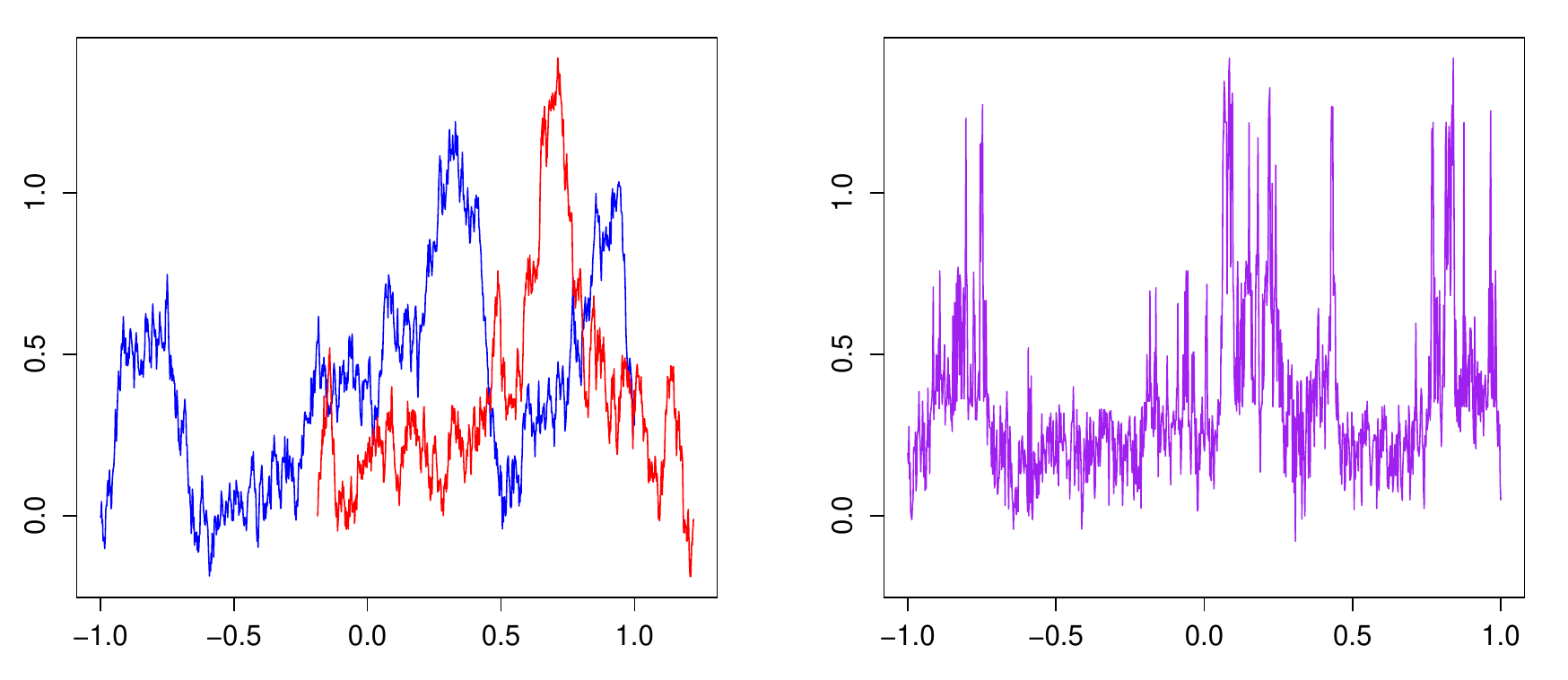}
    \vspace{-0.5cm}
    \caption{Composition of Gaussian processes results in rougher and more versatile sample paths. On the left: trajectories of two independent copies of a standard Brownian motion. On the right: the composition (red $\circ$ blue) of the trajectories.}
    \label{fig:fbm_dgp}
\end{figure}

Stabilization enhancing methods such as dropout and batch normalization are crucial for the performance of deep learning. In particular, batch normalization guarantees that the signal sent through a trained network cannot explode. We argue that for deep Gaussian processes similar effects play a role. In Figure~\ref{fig:fbm_dgp} below we visualize the effect of composing independent copies of a Gaussian process, the resulting trajectories are rougher and more versatile than those generated by the original process alone. This may, however, lead to wild behavior of the sample paths. As we aim for a fully Bayesian approach, the only possibility is to induce stability through the selection of the prior. We enforce stability by conditioning each individual Gaussian process to lie in a set of 'stable' paths. To achieve near optimal contraction rates, these sets have to be carefully selected and depend on the optimal contraction rate itself. 

The mathematical analysis requires to answer various questions involving the composition of GPs. To our knowledge, the closest results in the literature are bounds on the centred small-ball probabilities of iterated processes. They have been obtained for self-similar processes in \cite{MR2550290} and for time-changed self-similar processes in \cite{MR3474750}. A good reference on the literature of iterated processes is given by \cite{MR3638373}. In a different line of research, iterated Brownian motions (IBMs) occur in \cite{MR533542} as solutions of high-order parabolic stochastic differential equations (SDEs) and the path properties are studied in \cite{MR1278077}. The composition of general processes in relation with high-order parabolic and hyperbolic SDEs has been studied in \cite{MR1385409}. More recently, the \textit{ad libitum} (infinite) iteration of Brownian motions has been studied in \cite{MR3195821, MR3549711}.

The article is structured as follows. In Section~\ref{sec.model} we formalize the model and give an explicit parametrization of the underlying graph and the smoothness index. Section~\ref{sec.DGP} provides a detailed construction of the deep Gaussian process prior. In Section~\ref{sec.main} we state the main posterior contraction results. In Section~\ref{sec.dgp_minimax} we present a construction achieving optimal contraction rates and provide explicit examples in Section~\ref{sec.examples}. Section~\ref{sec.disc_dgp_prior} compares Bayes with DGPs and deep learning. It also contains a subsection on computational challenges. All proofs are deferred to the appendix.

\textit{Notation:} Vectors are denoted by bold letters, e.g. $\bx:=(x_1,\ldots,x_d)^\top.$ For $S\subseteq \{1,\ldots,d\},$ we write $\bx_S=(x_i)_{i \in S}$ and $|S|$ for the cardinality of $S.$ As usual, we define $|\bx|_p:= (\sum_{i=1}^d |\bx_i|^p)^{1/p},$ $|\bx|_\infty:= \max_i|\bx_i|,$ $|\bx|_0:= \sum_{i=1}^d \Ind(\bx_i \neq 0),$ and write $\|f\|_{L^p(D)}$ for the $L^p$ norm of $f$ on $D$. If there is no ambiguity concerning the domain $D,$ we also write $\|\cdot\|_p.$ For two sequences $(a_n)_n$ and $(b_n)_n$ we write $a_n \lesssim b_n$ if there exists a constant $C$ such that $a_n \leq C b_n$ for all $n.$ Moreover, $a_n\asymp b_n$ means that $(a_n)_n\lesssim(b_n)_n$ and $(b_n)_n\lesssim(a_n)_n.$ For positive sequences $(a_n)_n$ and $(b_n)_n$ we write $a_n \ll b_n$ if $a_n/b_n$ tends to zero when $n$ tends to infinity.

\section{Composition structure on the regression function}\label{sec.model}

In this section, we introduce the compositional class for the regression function $f^*$ in the nonparametric model~\eqref{eq.reg_model}. As we are interested in prediction, the aim is to learn $f^*$ but not its underlying compositional structure.

Consider the class of functions $f$ which can be written as the composition of $q +1$ functions, that is, $f  = g_{q }\circ g_{q -1}\circ \ldots\circ g_1\circ g_0,$ for functions $g_i:[a_i,b_i]^{d_i }\to [a_{i+1},b_{i+1}]^{d_{i+1} }$ with $d_0=d$ and $d_{q+1}=1.$ If $f $ takes values in the interval $[-1,1],$ rescaling $h_i=g_i(\|g_{i-1}\|_\infty \cdot )/\|g_i\|_\infty$ with $\|g_{-1}\|_\infty:=1$ leads to the alternative representation
\begin{align}\label{eq.f0}
	f  = h_{q }\circ h_{q -1}\circ \ldots\circ h_1\circ h_0
\end{align}
for functions $h_i:[-1,1]^{d_i }\to [-1,1]^{d_{i+1} }.$ We also write $h_i=(h_{ij})_{j=1,\ldots,d_{i+1} }^\top,$ with $h_{ij}:[-1,1]^{d_i} \to [-1,1].$ The representation can be modified if $f$ takes values outside $[-1,1],$ but to avoid unnecessary technical complications, we do not consider this case here. Although the function $h_i$ in the representation of $f$ is defined on $[-1,1]^{d_i},$ we allow each component function $h_{ij}$ to possibly only depend on a subset of $t_i$ 'active' variables $\mS_{ij} \subseteq \{1,\ldots,d_i\}$ for some $t_i \leq d_i$. Writing for a subset of indices $S,$ $(\cdot)_S: \bx\mapsto\bx_S=(x_i)_{i \in S},$ define
\begin{align}
        \overline h_{ij}:[-1,1]^{t_i}\to[-1,1], \ \  \bx_{\mS_{ij}}\mapsto h_{ij}(\bx_{\mS_{ij}},\bx_{{\mS_{ij}^c}}).
        \label{eq.overline_h_def}
\end{align}
As $h_{ij}$ does not depend on $\bx_{{\mS_{ij}^c}},$ this function is well-defined.

To define suitable function classes and priors on composition functions, it is natural to first associate to each composition structure a directed graph. The nodes in the graph are arranged in $q+2$ layers with $q+1$ the number of components in \eqref{eq.f0}. The number of nodes in each layer is given by the integer vector $\bd := (d,d_1 ,\ldots,d_q ,1) \in \N^{q +2}$ storing the dimensions of the components $h_i$ appearing in~\eqref{eq.f0}. In the graph, we draw an edge between the $j$-th node in the $i+1$-st layer and the $k$-th node in the $i$-th layer iff $k \in \mS_{ij}.$ The number $t_i$ is the in-degree of the nodes in layer $i+1.$

For any  $i,$ the subsets corresponding to different nodes $j=1,\ldots,d_{i+1},$ are combined into $\mS_i := (\mS_{i1},\ldots,\mS_{id_{i+1}})$ and $\mS := (\mS_0 ,\ldots,\mS_q).$ Setting $\bt:=(t_0,\dots, t_q),$ we summarize the previous quantities into the hyper-parameter
\begin{align}\label{def.lambda}
    \lambda := (q ,\bd ,\bt ,\mS ),
\end{align}
which we refer to as the graph of the function $f$ in~\eqref{eq.f0}. The set of all possible graphs is denoted by $\Lambda.$ 

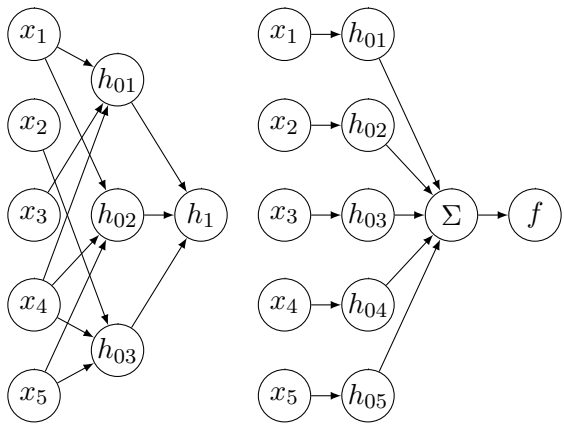
\begin{wrapfigure}{r}{0.5\textwidth}
\hspace{0.8cm}
\vspace{-0.8cm}
\centering

\begin{tikzpicture}[scale=0.5,
plain/.style={
  draw=none,
  fill=none,
  },
net/.style={
  matrix of nodes,
  nodes={
    draw,
    circle,
    inner sep=7pt
    },
  nodes in empty cells,
  column sep=0.4cm,
  row sep=-3pt
  },
>=latex,
]

\matrix[net] (mat)
{
|[plain]| & |[plain]|                 \\
          & |[plain]| & |[plain]| &   &         \\
|[plain]| &                              \\
          & |[plain]| & |[plain]| &   &      \\
|[plain]| & |[plain]|                       \\
          &           &           &   &   &   &      \\
|[plain]| & |[plain]|                       \\
          & |[plain]| & |[plain]| &   &      \\
|[plain]| &                                     \\
          & |[plain]| & |[plain]| &   &      \\
};

\draw (mat-2-1) -- (mat-2-1) node {$x_1$}; 
\draw (mat-4-1) -- (mat-4-1) node {$x_2$}; 
\draw (mat-6-1) -- (mat-6-1) node {$x_3$}; 
\draw (mat-8-1) -- (mat-8-1) node {$x_4$}; 
\draw (mat-10-1) -- (mat-10-1) node {$x_5$}; 

\draw (mat-3-2) -- (mat-3-2) node {$h_{01}$}; 

\draw (mat-6-2) -- (mat-6-2) node {$h_{02}$}; 

\draw (mat-9-2) -- (mat-9-2) node {$h_{03}$}; 

\draw (mat-6-3) -- (mat-6-3) node {$h_{1}$}; 

\draw[->] (mat-2-1) -- (mat-3-2);
\draw[->] (mat-6-1) -- (mat-3-2);
\draw[->] (mat-8-1) -- (mat-3-2);

\draw[->] (mat-2-1) -- (mat-6-2);
\draw[->] (mat-8-1) -- (mat-6-2);
\draw[->] (mat-10-1) -- (mat-6-2);

\draw[->] (mat-4-1) -- (mat-9-2);
\draw[->] (mat-8-1) -- (mat-9-2);
\draw[->] (mat-10-1) -- (mat-9-2);

\draw[->] (mat-3-2) -- (mat-6-3);
\draw[->] (mat-6-2) -- (mat-6-3);
\draw[->] (mat-9-2) -- (mat-6-3);

\draw (mat-2-4) -- (mat-2-4) node {$x_1$}; 
\draw (mat-4-4) -- (mat-4-4) node {$x_2$}; 
\draw (mat-6-4) -- (mat-6-4) node {$x_3$}; 
\draw (mat-8-4) -- (mat-8-4) node {$x_4$}; 
\draw (mat-10-4) -- (mat-10-4) node {$x_5$}; 

\draw (mat-2-5) -- (mat-2-5) node {$h_{01}$}; 
\draw (mat-4-5) -- (mat-4-5) node {$h_{02}$}; 
\draw (mat-6-5) -- (mat-6-5) node {$h_{03}$}; 
\draw (mat-8-5) -- (mat-8-5) node {$h_{04}$}; 
\draw (mat-10-5) -- (mat-10-5) node {$h_{05}$}; 

\draw (mat-6-6) -- (mat-6-6) node {$\Sigma$}; 

\draw (mat-6-7) -- (mat-6-7) node {$f$}; 

\draw[->] (mat-2-4) -- (mat-2-5);
\draw[->] (mat-4-4) -- (mat-4-5);
\draw[->] (mat-6-4) -- (mat-6-5);
\draw[->] (mat-8-4) -- (mat-8-5);
\draw[->] (mat-10-4) -- (mat-10-5);

\draw[->] (mat-2-5) -- (mat-6-6);
\draw[->] (mat-4-5) -- (mat-6-6);
\draw[->] (mat-6-5) -- (mat-6-6);
\draw[->] (mat-8-5) -- (mat-6-6);
\draw[->] (mat-10-5) -- (mat-6-6);

\draw[->] (mat-6-6) -- (mat-6-7);
\end{tikzpicture}

\caption{Graph representation of the example function (left) and generalized additive model (right).}
\label{fig:f_network}
\end{wrapfigure}

As an example consider the function $f(x_1,\ldots,x_5) = h_1(h_{01}(x_1,x_3,x_4),\newline h_{02}(x_1,x_4,x_5), h_{03}(x_2,x_4,x_5) )$ with corresponding graph representation displayed in Figure~\ref{fig:f_network}. This function is the composition of two layers and takes a five-dimensional input, while each of its components never involves more than three variables at once (the in-degree of all nodes is $3$). In this case, we have $q=1,$ $d_0=5, d_1=3, d_2=1$ and $t_0=t_1=3.$ The active sets are $\mS_{01}=\{1,3,4\}, \mS_{02}=\{1,4,5\}, \mS_{03}=\{2,4,5\},$ and $\mS_{11}=\{1,2,3\}$. Another example shown in Figure~\ref{fig:f_network} are generalized additive models (GAMs), where the functions are of the form $f(x_1,\ldots,x_d) = h_{21}(\sum_{i=1}^d h_{0i}(x_i))$. The first layer of the graph computes the function $x_i\mapsto h_{0i}(x_i)$, the second layer sums these functions, and the last layer finally computes the GAM $h_{21}(\sum_{i=1}^d h_{0i}(x_i))$. In this case we have $q=2,$ $\bd=(d,d,1,1),$ $\bt=(1,d,1),$ $\mS_{0i}=\{i\},$ for $i=1,\ldots,d,$ $\mS_{11}=\{1,\ldots,d\},$ and $\mS_{21}=\{1\}.$

We impose smoothness conditions on all functions in the composition. A function has H\"older smoothness index $\beta>0$ if all partial derivatives up to order $\lfloor\beta\rfloor$ exist and are bounded, and the partial derivatives of order $\lfloor\beta\rfloor$ are $(\beta-\lfloor\beta\rfloor)$-H\"older. Here, the ball of $\beta$-smooth H\"older functions of radius $K$ is defined as
\begin{align}\label{def.holder_ball}
	\begin{split}
		\mC^{\beta}_r(K) &= \Big\{f:[-1,1]^r\to [-1,1]: \\
		& 2r\sum_{\bm{\alpha}:|\bm{\alpha}|<\lfloor\beta\rfloor} \|\partial^{\bm{\alpha}} f\|_{\infty} + 2^{\beta-\lfloor \beta \rfloor} \sum_{\bm{\alpha}:|\bm{\alpha}|=\lfloor\beta\rfloor} \sup_{\substack{\bx,\by\in [-1,1]^r \\ \bx\neq \by}} \frac{|\partial^{\bm{\alpha}} f(\bx) - \partial^{\bm{\alpha}} f(\by)|}{|\bx-\by|_{\infty}^{\beta-\lfloor\beta\rfloor}} \leq K \Big\},
	\end{split}
\end{align}
where $\bm{\alpha}=(\alpha_1,\ldots,\alpha_r)\in\N^r$ is a multi-index, $|\bm{\alpha}|:=|\bm{\alpha}|_1$ and $\partial^{\bm{\alpha}} = \partial^{\alpha_1}\cdots\partial^{\alpha_r}.$ The factors $2r$ and $2^{\beta-\lfloor \beta \rfloor}$ guarantee the embedding $\mC^{\beta}_r(K)\subseteq \mC^{\beta'}_r(K)$ whenever $\beta'\leq \beta,$ see Lemma~\ref{lem.Di_embedding}. 

For $\beta=1,$ we recover the Lipschitz functions $|f(\bx)-f(\by)|\leq \tfrac K2 |\bx-\by|_\infty.$ If for a positive integer $\beta$ all the partial derivatives $\partial^{\bm{\alpha}}f$ with $|\bm{\alpha}|\leq \beta$ exist and are continuous, then it follows from the definition of the H\"older functions and first order Taylor expansion that $f\in \mC_r^\beta(K)$ for all sufficiently large $K.$ Moreover, if $f(x_1,\ldots,x_r)=g_1(x_1)\cdot \ldots \cdot g_r(x_r)$ with $g_1,\ldots,g_r\in \mC_1^\beta(K),$ then there exists a finite $K'$ such that $f\in \mC_r^\beta(K'),$ see Lemma~\ref{lem.holder.example} for a detailed statement and a proof.

We impose now H\"older smoothness on the functions $\overline h_{ij}$ in~\eqref{eq.overline_h_def}, assuming that $\overline h_{ij} \in \mC_{t_i}^{\beta_i}(K),$ with $\beta_i \in [\beta_{-},\beta_{+}]$ for some known and fixed lower and upper bounds $\beta_-,\beta_+$ satisfying $0<\beta_{-}\leq\beta_{+}<+\infty.$ The smoothness indices of the $q+1$ components are collected into the vector
\begin{align}\label{def.smooth}
    \bbeta  := (\beta_0 ,\ldots,\beta_q) \in [\beta_{-},\beta_{+}]^{q+1} =: I(\lambda).
    \end{align}
Combined with the graph parameter~\eqref{def.lambda}, the compositional functions \eqref{eq.f0} are completely described by
\begin{align}\label{def.model}
    \eta := (\lambda ,\bbeta ) = (q ,\bd ,\bt ,\mS , \bbeta ).
\end{align}
We refer to $\eta$ as the composition structure of the regression function $f$ in~\eqref{eq.f0}. The set of all possible choices of $\eta = (\lambda,\bbeta)$ with $\lambda\in\Lambda$ and $\bbeta \in I(\lambda)$ is denoted by $\Omega.$

Throughout the following, we assume that the true regression function $f^*$ belongs to the function space $\mF(\eta^*,K)$, for some unknown $\eta^*\in\Omega$, where
\begin{align}
\begin{split}\label{eq.support}
	\mF\big(\eta,K\big) := \Big\{f=h_{q }\circ h_{q -1}\circ \ldots\circ h_1\circ h_0:\ h_i=(h_{ij})_j:[-1,1]^{d_i }\to [-1,1]^{d_{i+1} },& \\ 
	\overline h_{ij}\in\mC_{t_i}^{\beta_i}\big(K\big) \Big\},
\end{split}
\end{align}
for some known $K>0$. Given a true regression function $f^*$, there might be several compositional function decompositions $\eta$ for which $f^*\in\mF(\eta,K).$ Our main result on posterior contraction rates states that if $f^*\in\mF(\eta,K)$ for a given $\eta$, then the posterior contracts with a rate $\varepsilon_n(\eta).$ This means that if $f^*\in\mF(\eta_j,K)$ for different choices $\eta_1,\ldots,\eta_m$, the smallest rate $\wedge_{j=1}^m\varepsilon_n(\eta_j)$ can be achieved.

%In principle, there might be an infinite number of models $\eta$ such that $f^*\in\mF(\eta,K)$ but, as we will see in Theorem~\ref{thm.post_contr} below, different choices of $\eta$ lead to different posterior contraction rates and the regression function $f^*$ will always be associated with the composition structure $\eta^*$ leading to the fastest contraction rate (in this sense, the most parsimonious model).

\section{Deep Gaussian process priors}\label{sec.DGP}

In this section, we construct deep Gaussian process priors as prior on composition functions. The hierarchical prior construction assigns first specific prior weights to all composition structures and, for any fixed structure, Gaussian process priors on all the individual functions that occur in the representation. To achieve fast contraction rates, the prior weights on the composition structures need to be selected carefully. Due to the complexity, the construction is split into several steps. 

%To achieve fast contraction rates, the prior weight assigned to a specific composition structure depends on the smoothness properties and the sample size. Therefore, the composition structure can not be decoupled from the estimation problem. We provide a discussion at the end of the section.

{\bf Step 0. Choice of Gaussian processes.} For a centered Gaussian process $X=(X_t)_{t\in T},$ the covariance operator viewed as a function on $T\times T,$ that is, $(s,t) \mapsto k(s,t)=\E[X_sX_t]$ is a positive semidefinite function. The reproducing kernel Hilbert space (RKHS) generated by $k$ is called the RKHS corresponding to the Gaussian process $X,$ see~\cite{vvvzRKHS} for more details.

For any dimension $r=1,2,\dots,$ and any $\beta>0,$ pick a a centered Gaussian process $\widetilde G^{(\beta,r)}=(\widetilde G^{(\beta,r)}(\bu))_{\bu\in[-1,1]^r}$ on the Banach space of continuous functions from $[-1,1]^r$ to $\R$ equipped with the supremum norm. We will use this process later as a prior for $\beta$-H\"older smooth functions defined on $[-1,1]^r$ in the compositional structure. Write $\|\cdot \|_{\bbH^{(\beta,r)}}$ for the RKHS-norm of the reproducing kernel Hilbert space $\bbH^{(\beta,r)}$ corresponding to $\widetilde G^{(\beta,r)}.$ For positive H\"older radius $K,$ we call
\begin{align}\label{eq.def_uniform_GP_conc}
    \varphi^{(\beta,r,K)}(u)
    := \sup_{f\in \mC_r^\beta(K)} \, 
    \inf_{g: \|g-f\|_\infty \leq u} \, \|g\|_{\bbH^{(\beta,r)}}^2 - \log \P\big(\big\|\widetilde G^{(\beta,r)}\big\|_\infty\leq u\big),
\end{align}
the concentration function over $\mC_r^\beta(K).$ This is the global version of the local concentration function appearing in the posterior contraction theory for Gaussian process priors, see \cite{vvvz}. For any $0<\alpha\leq 1,$ let $\eps_n(\alpha,\beta,r)$ be such that 
\begin{align}
\varphi^{(\beta,r,K)}\big(\eps_n(\alpha,\beta,r)^{1/\alpha}\big) \leq n \eps_n(\alpha,\beta,r)^2.
    \label{eq.def_eps_alpha}
\end{align}

{\bf Step 1. Deep Gaussian processes.} We now define a corresponding DGP $G^{(\eta)}$ on a given composition structure  $\eta=(q,\bd,\bt,\mS,\bbeta).$ Let $\B_\infty(R):=\{f:\sup_{\bx \in [-1,1]^r} |f(\bx)|\leq R\}$ be the supremum unitary ball with radius $R.$ For simplicity, we suppress the dependence on $r.$ Recall that the H\"older ball radius $K$ in \eqref{eq.support} is assumed to be known. With $\alpha_i :=\prod_{\ell=i+1}^q (\beta_{\ell} \wedge 1),$ the subset of paths 
\begin{align}
	\mD_i(\eta,K) := \B_\infty(1) \cap \left( \mC_{t_i}^{\beta_i}(K) + \B_\infty\big( 2\eps_n(\alpha_i,\beta_i,t_i)^{1/\alpha_i} \big) \right), 
	\label{eq.def.Di}
\end{align}
contains all functions that belong to the supremum unitary ball $\B_\infty(1)$ and are at most $2\eps_n(\alpha_i,\beta_i,t_i)^{1/\alpha_i}$-away in supremum norm from the H\"older-ball $\mC_{t_i}^{\beta_i}(K)$. With $\widetilde G^{(\beta,r)}$ the centred Gaussian process in Step~0, write $\overline{G}_i^{(\beta_i,t_i)}$ for the process $\widetilde G^{(\beta_i,t_i)}$ conditioned on the event $\{\widetilde G^{(\beta_i,t_i)} \in \mD_i(\eta,K)\}.$ Recall that for an index set $S,$ the function $(\cdot)_{S}$ maps a vector to the components in $S.$ For each $i=0,\ldots,q,$ $j=1,\ldots,d_{i+1},$ define the component functions $G_{ij}^{(\eta)}$ to be independent copies of the processes $\overline{G}_i^{(\beta_i,t_i)}\circ (\cdot )_{\mS_{ij}}:[-1,1]^{d_i}\to [-1,1].$ Finally, set $G_i^{(\eta)} := (G_{ij}^{(\eta)})_{j=1}^{d_{i+1}}$ and define the deep Gaussian process $G^{(\eta)} := G_q^{(\eta)}\circ\ldots\circ G_0^{(\eta)}:[-1,1]^d\to [-1,1].$ We denote by $\Pi(\cdot|\eta)$ the distribution of $G^{(\eta)}.$

{\bf Step 2. Structure prior.} We now construct a hyperprior on the underlying composition structure.
For $I(\lambda)$ as in \eqref{def.smooth} and for any function $a(\eta)=a(\lambda,\bbeta),$ it is convenient to define
\begin{align*}
    \int a(\eta) \, d\eta := \sum_{\lambda\in\Lambda} \int_{I(\lambda)} a(\lambda,\bbeta) \, d\bbeta. 
\end{align*}
Let $\gamma$ be a probability density on the possible composition structures, that is, $\int \gamma(\eta)\, d\eta = 1.$ We can construct such a measure $\gamma$ by first choosing a distribution on the number of compositions $q$. Given $q$ one can then select distributions on the ambient dimensions $\bd,$ the effective dimensions $\bt,$ the active sets $\mS$ and finally the smoothness indices $\bbeta \in [\beta_{-}, \beta_{+}]^{q+1}$ via the conditional density formula $\gamma(\eta) = \gamma(\lambda)\gamma(\bbeta|\lambda) = \gamma(q)\gamma(\bd|q)\gamma(\bt|\bd,q)\gamma(\mS|\bt,\bd,q)\gamma(\bbeta|\lambda).$
For a sequence $\eps_n(\eta)$ satisfying
\begin{align}
    \eps_n(\eta) \geq 
    \max_{i=0, \dots,q} \eps_n(\alpha_i,\beta_i,t_i), \ \ \text{with}\ \ \alpha_i :=\prod_{\ell=i+1}^q \big(\beta_{\ell} \wedge 1\big),
    \label{eq.def_rate}
\end{align}
and $|\bd|_1=1+\sum_{i=0}^q d_i,$ consider the hyperprior
\begin{align}
\begin{split}\label{def.prior}
    \pi(\eta) & := \frac{e^{-\Psi_n(\eta)} \gamma(\eta)}{
	\int e^{-\Psi_n(\eta)} \gamma(\eta) \, d\eta}, \ \ \text{with}\ \ \Psi_n(\eta) := n\eps_n(\eta)^2 + e^{e^{|\bd|_1}}. 
\end{split}
\end{align}
The denominator is positive and finite, since $0 < e^{-\Psi_n(\eta)} \leq 1$ and $\int \gamma(\eta) \, d\eta = 1.$ 

\textbf{Step 3. DGP prior.} We consider deep Gaussian process priors of the form
\begin{align}\label{def.dgp_prior}
    \Pi(df) := \int_{\Omega} \Pi(df|\eta) \pi(\eta)\, d\eta,
\end{align}
where $\Omega$ is the set of all valid composition structures, $\Pi(\cdot|\eta)$ is the distribution of the DGP $G^{(\eta)}$ and $\pi(\eta)$ is the structure prior on $\eta.$
\newline

\noindent
\textbf{Some remarks on the DGP prior.} Step~0 allows for a large class of Gaussian processes, since the only requirement is that the concentration function inequality \eqref{eq.def_eps_alpha} admits solutions for any $0<\alpha\leq1$. Lemma~\ref{lem.smallest_solution} shows that it is often enough to check~\eqref{eq.def_eps_alpha} for $\alpha=1$ only. 

Composing Gaussian processes leads to wild sample paths with non-negligible probability. We overcome the issue in Step~1, by conditioning the Gaussian processes onto the sets $\mD_i(\eta,K)$. This is well-defined since $\mD_i(\eta,K)$ contains the sup-norm ball $\B_\infty\big( 2\eps_n(\alpha_i,\beta_i,t_i)^{1/\alpha_i} \big)$ and Gaussian processes with continuous sample paths give positive mass to $\B_\infty(R)$ for any positive radius $R>0.$ By choosing suitable Gaussian processes, the conditioning is not very restrictive. Brownian motion, for instance, is known to have sample paths that are almost surely $\beta$-H\"older continuous for any $\beta<1/2.$ The H\"older norm is, however, random. The conditioning step requires then to constraint the prior to Brownian motion sample paths with H\"older norm bounded by $K.$

Step~2 of the DGP prior construction introduces the rate $\eps_n(\eta).$ This sequence will later be shown to be the posterior contraction rate if the true regression function $f^*$ is an element of the compositional space $\mF(\eta,K).$ To achieve fast posterior contraction rates, it is desirable to pick a small $\eps_n(\eta)$ in~\eqref{eq.def_rate}. On the contrary, Assumption \ref{ass.rate} below imposes additional restrictions on the rate.

The prior $\pi(\eta)$ on the composition structure $\eta$ should be viewed as a model selection prior, see also Section~10 in~\cite{GVDV17}. As always, some care is needed to avoid posterior contraction on models that are too large and consequently lead to overfitting and suboptimal posterior contraction rates. This is achieved by the carefully chosen exponent $\Psi_n$ in \eqref{def.prior}, which depends on the sample size and penalizes large composition structures.

\section{Main results} \label{sec.main}

Denote by $\Pi\big(\cdot|\bX,\bY\big)$ the posterior distribution corresponding to a DGP prior $\Pi$ constructed as above and $(\bX,\bY) = (\bX_i,Y_i)_i$ a sample from the nonparametric regression model \eqref{eq.reg_model}. For normalizing factor $Z_n: = \int p_f/p_{f^*}(\bX,\bY) \, \Pi(df)$ and any Borel measurable $\mA$ in the Banach space of continuous functions on $[-1,1]^d,$ 
\begin{align}
\begin{split}\label{def.post}
    \Pi\big(\mA|\bX,\bY\big) &= Z_n^{-1}  \int_\mA \frac{p_f}{p_{f^*}}(\bX,\bY) \, \Pi(df),\quad \Pi(df) := \int_{\Omega} \Pi(df|\eta)\pi(\eta)\, d\eta,
\end{split}
\end{align}
where $(p_f/p_{f^*})(\bX,\bY)$ denotes the likelihood ratio. With a slight abuse of notation, for any subset of composition structures $\mM\subseteq\Omega,$ we set
\begin{align}\label{def.post_model}
    \Pi\big(\eta\in\mM|\bX,\bY\big) &:= Z_n^{-1} \int \frac{p_f}{p_{f^*}}(\bX,\bY) \int_\mM \Pi(df|\eta) \pi(\eta) \, d\eta,
\end{align}
which is the contribution of the composition structures $\mM\subseteq\Omega$ to the posterior mass. 

Before we can state the results, we first need to impose some conditions. The first condition is on the  distribution $\gamma$ appearing in the graph prior in Step~2 of the DGP prior construction. While uniformity of $\gamma(\cdot|\lambda)$ is the most natural choice and simplifies the mathematical analysis considerably, the next condition also states that all graphs have to be charged with non-negative mass and requires $\int \sqrt{\gamma(\eta)} \, d\eta$ to be finite. The latter restricts the amount of prior mass that can be assigned to complex composition structures.
\begin{assump}\label{ass.prior}
    We assume that, for any graph $\lambda,$ the measure $\gamma(\cdot|\lambda)$ is the uniform distribution on the hypercube of possible smoothness indices $I(\lambda) = [\beta_{-},\beta_{+}]^{q+1}.$ Furthermore, we assume that the distribution $\gamma$ is independent of $n,$ that it assigns positive mass $\gamma(\eta)>0$ to all composition structures $\eta,$ and that it satisfies $\int \sqrt{\gamma(\eta)} \, d\eta < +\infty.$ 
\end{assump}

The prior on $\gamma$ can be constructed hierarchically by using the decomposition $\gamma(\eta) = \gamma(q)\gamma(\bd|q)\gamma(\bt|\bd,q)\gamma(\mS|\bt,\bd,q)(\beta_+-\beta_-)^{-q-1}.$ It is natural to draw $q$ and all components of $\bd|q$ independently from distributions on the positive integers and generate $\bt|(\bd,q)$ and $\mS|(\bt,\bd,q)$ independently from uniform distributions on the respective sample spaces. The distribution $\gamma$ assigns positive mass to all composition structures. The next result shows that the previous assumption is satisfied under suitable moment conditions. 

\begin{lemma}
\label{lem.A1_checked}
If the prior $\gamma$ is generated hierarchically as above, where $q$ is drawn from a distribution with finite moment $\E_q[A^q]$ for all $A>0,$ and all components $d_1|q, \ldots,d_q|q$ are i.i.d. with $\E_{d_1|q}[d_1^32^{d_1^2}]<\infty,$ then, $\gamma$ satisfies Assumption \ref{ass.prior}.
\end{lemma}

%Our proof extends to the case when $\gamma(\cdot|\lambda)$ is not the uniform distribution.

The second assumption deals with the rates in Step~2 of the DGP prior construction. We impose a lower bound on $\eps_n(\alpha,\beta,r)$ and also require that, uniformly over all models $\eta=(\lambda,\bbeta),\eta'=(\lambda,\bbeta')$ sharing the same graph $\lambda$ and having similar smoothness indices $\bbeta,\bbeta'$, the corresponding rates $\eps_n(\eta),\eps_n(\eta')$ only differ by multiplicative factors. 

\begin{assump}\label{ass.rate}
    We assume the following on the rates appearing in the construction of the prior.
    \begin{enumerate}[noitemsep,topsep=0pt,label=(\roman*),ref=(\roman*)]
        \item \label{ass.rate_low_bound} For any positive integer $r,$ any $\beta>0,$ let $Q_1(\beta,r,K)$ the constant from Lemma~\ref{lem.Holder_entropy_sharp}. Then, the sequences $\eps_n(\alpha,\beta,r)$ solving the concentration function inequality~\eqref{eq.def_eps_alpha} are chosen in such a way that
        \begin{align}\label{eq.rate_low_bound}
            \eps_n(\alpha,\beta,r) \geq Q_1(\beta,r,K)^{\frac{\beta}{2\beta + r}} n^{-\frac{\beta\alpha}{2\beta\alpha+r}}.
        \end{align}
        \item \label{ass.rate_comparison} There exists a constant $Q\geq 1$ such that the following holds. For any $n>1,$ any graph $\lambda = (q,\bd,\bt,\mS)$ with $|\bd|_1 = 1 + \sum_{i=0}^q d_i \leq \log(2\log n),$ and any $\bbeta'=(\beta_0',\dots,\beta_q'),\bbeta=(\beta_0,\dots,\beta_q)\in I(\lambda)$ satisfying $\beta_i'\leq \beta_i \leq \beta_i'+1/\log^2 n$ for all $i=0,\dots,q,$ the rates relative to the composition structures $\eta=(\lambda,\bbeta)$ and $\eta'=(\lambda,\bbeta')$ satisfy 
        \begin{align}\label{eq.rate_comparison}
            \eps_n(\eta) \leq \eps_n(\eta') \leq Q \eps_n(\eta).
        \end{align}
    \end{enumerate}
\end{assump}

The previous assumptions are checked for a number of examples in Section \ref{sec.examples}. The rate $\eps_n(\eta)$ associated to a composition structure $\eta$ can be viewed as measure of the complexity of this structure, where larger rates $\eps_n(\eta)$ correspond to more complex models. Our first result states that the posterior concentrates on small models in the sense that all posterior mass is asymptotically allocated on a set
\begin{align}\label{def.M_n_C}
    \mM_n(C):=\big\{\eta: \eps_n(\eta) \leq C \eps_n(\eta^*) \big\} \cap \big\{\eta: |\bd|_1 \leq \log(2\log n) \big\}
\end{align}
with sufficiently large constant $C.$ This shows that the posterior not only concentrates on models with fast rates $\eps_n(\eta)$ but also on graph structures with number of nodes in each layer bounded by $\log(2\log n)$. The proof is given in Appendix~\ref{sec.proofs_main}. 

\begin{theorem}[Model selection]
	\label{thm.mod_selection}
	Let $\Pi\big(\cdot|\bX,\bY\big)$ be the posterior distribution corresponding to a DGP prior $\Pi$ constructed as in Section~\ref{sec.DGP} and satisfying Assumptions~\ref{ass.prior} and~\ref{ass.rate}. Let $\eta^*=(\lambda^*,\bbeta^*)$ for some $\bbeta^* \in (\beta_-,\beta_+)^{q^*+1}$ and suppose $\eps_n(\eta^*)\leq 1/(4Q)$. Then, for a positive constant $C=C(\eta^*),$ 
	\begin{align*}
	    \sup_{f^* \in \mF(\eta^*,K)} \ \E_{f^*}\big[\Pi\big( \eta \notin \mM_n(C)  \big|\bX,\bY\big)\big] \xrightarrow{n\to\infty} 0,
	\end{align*}
	where $\E_{f^*}$ denotes the expectation with respect to $\P_{f^*},$ the true distribution of the sample $(\bX,\bY).$ 
\end{theorem}
Denote by $\mu$ the distribution of the covariate vector $\bX_1$ and write $L^2(\mu)$ for the weighted $L^2$-space with respect to the measure $\mu.$ The next result shows that the posterior distribution achieves contraction rate $\eps_n(\eta^*)$ up to a $\log n$ factor. The proof is based on the testing approach developed in \cite{MR1790007} and given in Appendix~\ref{sec.proofs_main}.

\begin{theorem}[Posterior contraction]
	\label{thm.post_contr}
	Let $\Pi\big(\cdot|\bX,\bY\big)$ be the posterior distribution corresponding to a DGP prior $\Pi$ constructed as in Section~\ref{sec.DGP} and satisfying Assumptions~\ref{ass.prior} and~\ref{ass.rate}. Let $\eta^*=(\lambda^*,\bbeta^*)$ for some $\bbeta^* \in (\beta_-,\beta_+)^{q^*+1}$ and suppose $\eps_n(\eta^*)\leq 1/(4Q)$. Then, for a positive constant $L=L(\eta^*)$,
	\begin{align*}
	    \sup_{f^* \in \mF(\eta^*,K)} \ \E_{f^*}\left[\Pi\left( \|f-f^*\|_{L^2(\mu)} \geq L (\log n)^{1+\log K} \eps_n(\eta^*) \big|\bX,\bY\right)\right] \xrightarrow{n\to\infty} 0,
	\end{align*}
	where $\E_{f^*}$ denotes the expectation with respect to $\P_{f^*},$ the true distribution of the sample $(\bX,\bY).$
\end{theorem}

\begin{remark}\label{sec.weak_constr}
Our proving strategy allows for the following modification to the construction of the DGP prior. The concentration functions $\varphi^{(\beta,r,K)}$ in~\eqref{eq.def_uniform_GP_conc} are defined globally over the H\"older-ball $\mC_r^\beta(K).$ The concentration function inequality in~\eqref{eq.def_eps_alpha} essentially requires that the closure $\overline\bbH^{(\beta,r)}$ of the RKHS of the underlying Gaussian process $\widetilde G^{(\beta,r)}$ contains the whole H\"older ball. There are classical examples for which this is too restrictive, and one might want to weaken the construction by considering a subset $\mH_r^\beta(K)\subseteq\mC_r^\beta(K).$ This can be done by replacing $\mC_{t_i}^{\beta_i}(K)$ with the corresponding subset $\mH_{t_i}^{\beta_i}(K)$ in the definition of the conditioning sets $\mD_i(\eta,K)$ in~\eqref{eq.def.Di} and the
function class $\mF(\eta,K)$ in~\eqref{eq.support}. As a consequence, this also reduces the class of functions for which the posterior contraction rates derived in Theorem \ref{thm.post_contr} hold.
\end{remark}

We do not impose an a-priori known upper bound on the complexity of the underlying composition structure \eqref{eq.f0}. While we think that this is natural in practice, it causes extra technical complications resulting for instance in the appearance of the bound $|\bd|_1\leq\log(2\log n)$ in Assumption~\ref{ass.rate} and in the definition of the model class $\mM_n(C)$ in~\eqref{def.M_n_C}. If we additionally assume that the true composition structure satisfies $|\bd|_1\leq D$ for a known upper bound $D,$ then the factor $e^{e^{|\bd|_1}}$ in \eqref{def.prior} can be avoided. Moreover, the $(\log n)^{1+\log K}$-factor occurring in the posterior contraction rate is somehow an artifact of the proof, and could be replaced by $K^D,$ see the proof of Lemma \ref{lem.local_entropy} for more details. A trade-off regarding the choice of $K$ appears. To allow for larger classes of functions and a weaker constraint induced by the conditioning on \eqref{eq.def.Di}, we want to select a large $K.$ On the contrary, large $K$ results in slower posterior contraction guarantees.

The main result requires that all smoothness indices lie in the interval $[\beta_-,\beta_+]$ with $0<\beta_-<\beta_+<\infty.$  As commonly observed in nonparametric Bayes, extension to $(0,\infty)$ is highly non-trivial as many constants depend in an intricate way on $\beta_-, \beta_+$ and quickly explode in the limits $\beta_-\to 0$ and $\beta_+\to \infty.$ A prototypical example for the latter is the constant $e^{\beta_+}$ in Lemma \ref{lem.eps_ratio} below. One could also wonder what would happen in the misspecified case where the true $\bbeta^*$ lies outside of the compact set $[\beta_{-},\beta_{+}]^{q^*+1}$. For nonparametric regression using GPs, it is known that extending the range of $\bbeta^*$ while keeping the prior on $\bbeta$ as before, the rates become suboptimal by a polynomial factor in the sample size $n$ due to the mismatch between the smoothness of the prior and that of the true regression function, see~\cite{castillo2008lower} for the Riemann-Liouville prior.

%The nearly optimal construction involving the Gaussian processes in Section~\ref{sec.dgp_minimax} requires that some quantities are bounded uniformly over $\beta\in[\beta_{-},\beta_{+}]$, such as the constants $Q_1(\beta,r,K)$ and $Q$ in Assumption~\ref{ass.rate}. One can see that $Q_1(\beta_{-},r,K)$ is independent of the Gaussian priors and is proportional to $K^{2r/\beta{-}}$ and, under the nearly optimal construction, $Q=e^{\beta_{+}}$. Both these quantities are unbounded when $\beta_{-}^{-1},\beta_{+}$ are arbitrarily large. One could argue heuristically that, if these quantities grow to infinity slower than $\log\log n$, then the resulting posterior should still achieve nearly optimal rates since also the constant $L$ in Theorem~\ref{thm.post_contr} is proportional to $Q=e^{\beta_{+}}$.  

%As an example, consider the extreme case where all entries of $\bbeta^*$ are larger than $\beta_{+}$. The embeddings of the H\"older spaces lead to posterior contraction rates $\eps_n(\eta_{+})$ where $\eta_{+} = (\lambda^*, \bbeta_{+})$ and $\bbeta_{+}$ is the vector whose entries are all equal to $\beta_{+}$. This rate is slower than $\eps_n(\eta^*)$ due to the true regression function being smoother than the DGP prior.

For mathematical convenience, all functions in the function class~\eqref{eq.support} have range $[-1,1]$. In general, any function $f$ mapping to $[-a,a]$ can be represented as in~\eqref{eq.f0} where all intermediate components $h_{ij}$ in the compositional class have range $[-1,1]$ and only the very last $h_q$ has range $[-a,a]$. In the construction of the DGP prior, we then need to restrict for $i=q$ in \eqref{eq.def.Di} to the set $\B_\infty(a)$ instead. The mathematical analysis can be extended. In particular, in the Hellinger-$L^2(\mu)$ equivalence in~\eqref{eq.Hell_lb}, the constant factor becomes $e^{-a^2/2}$.

By making appropriate changes in the proofs, one can also include in the nonparametric model~\eqref{eq.reg_model} a known standard deviation $\sigma^*>0$ instead of taking standard normal noise. As the analysis is already quite technical, we have chosen to avoid additional parameters that are not of primary focus.

We view the proposed Bayesian analysis rather as a proof of concept than something that is straightforward implementable or computationally efficient. The main obstacles towards a scalable Bayesian method are the combinatorial nature of the set of graphs as well as conditioning the sample paths to neighborhoods of H\"older functions, see also Section~\ref{sec.disc_dgp_prior} for a more in-depth discussion.

\section{On nearly optimal contraction rates} \label{sec.dgp_minimax}

Theorem 2.5 in \cite{MR1790007} ensures the existence of a frequentist estimator converging to the true parameter with the posterior contraction rate. This implies that the fastest possible posterior contraction rate is the minimax estimation rate. For the prediction loss $\|f-g\|_{L^2(\mu)}$, the minimax estimation rate over the class $\mF(\eta,K)$ is, up to some logarithmic factors,
\begin{align}\label{eq.rate_minimax}
	\frakr_n(\eta) = \max_{ i=0,\ldots, q} n^{-\frac{\beta_i\alpha_i}{2\beta_i\alpha_i+t_i}}, \ \ \text{with}\ \ \alpha_i :=\prod_{\ell=i+1}^q \big(\beta_{\ell} \wedge 1\big),
\end{align}
see \cite{SH19}. This rate is attained by suitable estimators based on sparsely connected deep neural networks. It is also shown in the recent work \cite{2022arXiv220507764G} that there are compositional structures for which GP priors lead to a posterior contraction rate that is suboptimal by a polynomial factor compared with $\frakr_n(\eta).$

Below we derive sufficient conditions that are simpler than Assumption~\ref{ass.rate}, apply to standard examples of Gaussian processes in the Bayes literature and imply an optimal posterior contraction rate $\frakr_n(\eta)$ up to $\log n$-factors. 

The first result shows that the solution to the concentration function inequality for arbitrary $0< \alpha\leq 1$ can be deduced from the solution for $\alpha=1.$ The proof is in Appendix~\ref{sec.proofs_dgp_optimal}.
\begin{lemma}\label{lem.smallest_solution}
    Let $\eps_n(1,\beta,r)$ be a solution to the concentration function inequality~\eqref{eq.def_eps_alpha} for $\alpha=1.$ Then, any sequence $\eps_n(\alpha,\beta,r) \geq \eps_{m_n}(1,\beta,r)^\alpha$ where $m_n$ is chosen such that $m_n \eps_{m_n}(1,\beta,r)^{2-2\alpha} \leq n,$ solves the concentration function inequality for arbitrary $\alpha\in(0,1].$
\end{lemma}

The next result identifies sequences $\eps_n(\alpha,\beta,r)$ satisfying inequality~\eqref{eq.def_eps_alpha} provided the solution $\eps_n(1,\beta,r)$ for $\alpha=1$ is, up to $\log n$ factors, $n^{-\beta/(2\beta+r)}$. The proof is given in Appendix~\ref{sec.proofs_dgp_optimal}.
\begin{lemma}\label{lem.rates} Let $n\geq 3.$
\begin{compactitem}
\item[(i)] If the sequence $\eps_n(1,\beta,r) = C_1 (\log n)^{C_2} n^{-\beta/(2\beta+r)}$ solves the concentration function inequality~\eqref{eq.def_eps_alpha} for $\alpha = 1,$ with constants $C_1\geq1$ and $C_2\geq0,$ then, any sequence
    \begin{align*}
        \eps_n(\alpha,\beta,r) \geq C_1^2 (2\beta+1)^{ 2C_2} (\log n)^{C_2 (2\beta+2)} \ n^{-\frac{\beta\alpha}{2\beta\alpha+r}}
    \end{align*}
    solves the concentration function inequality for arbitrary $\alpha\in(0,1].$
\item[(ii)] If there are constant $C_1'\geq 1$ and $C_2'\geq 0$ such that the concentration function satisfies
    \begin{align} \label{eq.GP_conc_strong}
        \varphi^{(\beta,r,K)}(\delta) \leq C_1' (\log \delta^{-1})^{C_2'} \delta^{-\frac{r}{\beta}},   \ \ \text{for all} \ 0<\delta\leq 1,
    \end{align}
then, any sequence
    \begin{align*}
        \eps_n(\alpha,\beta,r) \geq C_1' (\beta\log n)^{C_2'} \, n^{-\frac{\beta\alpha}{2\beta\alpha+r}}
    \end{align*}
    solves the concentration function inequality~\eqref{eq.def_eps_alpha} for arbitrary $\alpha\in(0,1].$    
\end{compactitem}    
\end{lemma}
To verify Assumption~\ref{ass.rate}, we need to pick suitable sequences $\eps_n(\eta)\geq \max_{i=0, \dots,q} \eps_n(\alpha_i,\beta_i,t_i).$ In the previous lemma, $\eps_n(\alpha,\beta,r) = C_1(\beta,r)(\log n)^{C_2(\beta,r)}n^{-\beta\alpha/(2\beta\alpha+r)}$ for some constants $C_1(\beta,r)\geq1$ and $C_2(\beta,r)\geq0.$ A suitable choice in this case is
\begin{align}\label{eq.def_rate_enlarge}
    \eps_n(\eta)
    = \widetilde C_1(\eta)(\log n)^{\widetilde C_2(\eta)} \frakr_n(\eta),
\end{align}
provided the constants $\widetilde C_j(\eta):=\max_{i=0,\dots,q} \sup_{\beta \in [\beta_-,\beta_+]} C_j(\beta,t_i),$ $j\in \{1,2\}$ are finite. This can be checked by verifying that $\beta \mapsto C_j(\beta,r),$ $j\in \{1,2\}$ are bounded functions on $[\beta_-,\beta_+].$

\begin{lemma} \label{lem.eps_ratio}
    The rates $\eps_n(\eta)$ in~\eqref{eq.def_rate_enlarge} satisfy condition~\ref{ass.rate_comparison} in Assumption~\ref{ass.rate} with $Q = e^{\beta_{+}}.$
\end{lemma}    

Let $\Pi$ be a DGP prior constructed with the Gaussian processes and rates given in this section. Then, the corresponding posterior satisfies Theorem~\ref{thm.post_contr} and contracts with rate $\frakr_n(\eta^*)$ up to the multiplicative factor $L(\eta^*) \widetilde C_1(\lambda^*,K) (\log n)^{1+\log K + \widetilde C_2(\lambda^*,K)}$. An examination of the proof gives $L(\eta^*) = M \cdot 10 C e^{\beta_{+}}$ for $M>0$ a sufficiently large universal constant.

\section{Examples of DGP priors} \label{sec.examples}

In this section we verify the abstract conditions for standard choices of Gaussian processes $\{\widetilde G^{(\beta,r)} : \beta\in [\beta_{-},\beta_{+}]\}$ and show that they achieve near optimal posterior contraction rates.

\subsection{L\' evy's fractional Brownian motion} 
Assume that the upper bound $\beta_+$ on the possible range of smoothness indices is bounded by one. A zero-mean Gaussian process $X^\beta$ is called a L\'evy fractional Brownian motion of order $\beta\in(0,1)$ if 
\begin{align*}
	X^\beta(0)=0,\quad \E\left[|X^\beta(\bu) - X^\beta(\bu')|^2\right] = |\bu-\bu'|_2^{2\beta}, \quad  \forall \bu,\bu'\in[-1,1]^r.
\end{align*}
The covariance function of the process is $\E[X^\beta(\bu) X^\beta(\bu')] = \tfrac{1}{2}(|\bu|_2^{2\beta} + |\bu'|_2^{2\beta} - |\bu-\bu'|_2^{2\beta} ).$ Chapter~3 in~\cite{cohen2013fractional} provides the following representation for a $r$-dimensional $\beta$-fractional Brownian motion. Denote by $\widehat{f}(\bxi) := (2\pi)^{-r/2} \int_{\R^r} e^{i\bu^\top\bxi}f(\bu)d\bu$ the Fourier transform of the function $f$. For $W = (W(\bu))_{\bu\in[-1,1]^r}$ a multidimensional Brownian motion, and $C_\beta$ a positive constant depending only on $\beta,\, r$,
\begin{align*}
	X^\beta(\bu) = \int_{\R^r}  \frac{e^{-i\bu^\top\bxi}-1}{C_\beta^{1/2}|\bxi|_2^{\beta+r/2}} \widehat{W}(d\bxi),
\end{align*}
in distribution, where $\widehat W(d\bxi)$ is the Fourier transform of the Brownian random measure $W(d\bu)$, see Section~2.1.6 in~\cite{cohen2013fractional} for definitions and properties. The same reference defines, for all $\varphi\in L^2([-1,1]^r),$ the integral operator
\begin{align*}
	(I^\beta \varphi)(\bu) &:= \int_{\R^r} \overline{\widehat{\varphi}(\bxi)} \frac{e^{-i\bu^\top\bxi}-1}{C_\beta^{1/2}|\bxi|_2^{\beta+r/2}} \frac{d\bxi}{(2\pi)^{r/2}}.
\end{align*}
The RKHS $\bbH^\beta$ of $X^\beta$ is given in Section~3.3 in~\cite{cohen2013fractional} by
\begin{align*}
	\bbH^\beta = \Big\{ I^\beta\varphi: \varphi\in L^2([-1,1]^r) \Big\},\quad\langle I^\beta\varphi, I^\beta\varphi' \rangle_{\bbH^\beta} = \langle \varphi,\varphi' \rangle_{L^2([-1,1]^r)}.
\end{align*}
The process $X^\beta$ is always zero at $\bu=0.$ To release it at zero, let $Z\sim\mN(0,1)$ be independent of $X^\beta$ and consider the process $\bu\mapsto Z+X^\beta(\bu).$ The RKHS of the constant process $\bu\mapsto Z$ is the set $\bbH^Z$ of all constant functions and, by Lemma~I.18 in~\cite{GVDV17}, the RKHS of $Z+X^\beta$ is the direct sum $\bbH^Z\oplus\bbH^\beta.$

The next result is proved in Appendix~\ref{sec.proofs_ex} and can be viewed as the multidimensional extension of the RKHS bounds in Theorem~4 in~\cite{castillo2008lower}. Whereas the original proof relies on kernel smoothing and Taylor approximations, we use a spectral approach. Write
\begin{align}\label{eq.Sobolev}
    \mW_r^{\beta}(K) := \bigg\{h:[-1,1]^r\to[-1,1] : \int_{\R^r} |\widehat h(\bxi)|^2 (1+|\bxi|_2)^{2\beta} \frac{d\bxi}{(2\pi)^{r/2}} \leq K \bigg\},
\end{align}
for the $\beta$-Sobolev ball of radius $K.$

\begin{lemma}\label{lem.fbm_rkhs}
	Let $\beta\in[\beta_{-},1]$ and $Z+X^\beta = (Z+X^\beta(\bu))_{\bu\in[-1,1]^r}$ the fractional Brownian motion of order $\beta$ released at zero. Fix $h\in C_r^{\beta}(K)\cap\mW_r^{\beta}(K).$ Set $\phi_\sigma=\sigma^{-r}\phi(\cdot/\sigma)$ with $\phi$ a suitable regular kernel and $\sigma<1.$ Then, $\|h*\phi_\sigma-h\|_\infty \leq K R_\beta \sigma^{\beta}$ and $\|h*\phi_\sigma\|_{\bbH^Z\oplus\bbH^\beta}^2 \leq K^2 L_\beta^2 \sigma^{-r}$ for some constants $R_\beta,L_\beta$ that depend only on $\beta,\ r.$
\end{lemma}

The next lemma shows that for L\'evy's fractional Brownian motion released at zero near optimal posterior contraction rates can be obtained. For that we need to restrict the definition of the global concentration function to the smaller class $\mC_r^\beta(K) \cap \mW_r^\beta(K).$ The proof of the lemma is in Appendix~\ref{sec.proofs_ex}.

\begin{lemma}\label{lem.fbm_ass2}
Let $\beta_+\leq 1$ and work on the reduced function spaces $\mH_r^\beta(K)=\mC_r^\beta(K) \cap \mW_r^\beta(K)$ as outlined in Remark~\ref{sec.weak_constr}. For $\{\widetilde G^{(\beta,r)}:\beta\in[\beta_{-},\beta_{+}]\}$ the family of Levy's fractional Brownian motions $Z+X^\beta$ released at zero, there exist sequences $\eps_n(\eta)=C_1(\eta)(\log n)^{C_2(\eta)}\mathfrak{r}_n(\eta)$ such that Assumption~\ref{ass.rate} holds.   
\end{lemma}

\subsection{Truncated wavelet series}\label{sec.ex_truncated}

Let $\{\psi_{j,k}:j\in\N_+,\ k=1,\ldots 2^{jr}\}$ be an orthonormal wavelet basis of $L^2([-1,1]^r).$ For any $\varphi \in L^2([-1,1]^r),$ we denote by $\varphi = \sum_{j=1}^{\infty} \sum_{k=1}^{2^{jr}} \lambda_{j,k}(\varphi) \psi_{j,k}$ its wavelet expansion. The quantities $\lambda_{j,k}(\varphi)$ are the corresponding real coefficients. For any $\beta>0,$ we denote by $\mB_{\infty,\infty,\beta}$ the Besov space of functions $\varphi$ with finite
\begin{align*}
    \|\varphi\|_{\infty,\infty,\beta} &:= \sup_{j\in\N} 2^{j(\beta+\frac{r}{2})}\max_{k=1,\ldots 2^{jr}} |\lambda_{j,k}(\varphi)|.
\end{align*}
We assume that the wavelet basis is $s$-regular with $s>\beta_{+}.$ For i.i.d. random variables $Z_{j,k}\sim\mN(0,1),$ consider the Gaussian process induced by the truncated series expansion
\begin{align*}
    X^\beta(\bu) := \sum_{j=1}^{J_\beta} \sum_{k=1}^{2^{jr}} \frac{2^{-j(\beta + \frac{r}{2})}}{\sqrt{jr}} Z_{j,k} \psi_{j,k}(\bu),
\end{align*}
where the maximal resolution $J_\beta$ is chosen as the integer closest to the solution $J$ of the equation $2^{J} = n^{1/(2\beta+r)},$ see Section~4.5 in~\cite{vvvz}. The RKHS of the process $X^\beta$ is given in the proof of Theorem~4.5 in~\cite{vvvz} as the set $\bbH^\beta$ of functions $\varphi = \sum_{j=1}^{J_\beta}\sum_{k=1}^{2^{jr}} \lambda_{j,k}(\varphi) \psi_{j,k}$ with coefficients $\lambda_{j,k}(\varphi)$ satisfying $\|\varphi\|_{\bbH^\beta}^2 := \sum_{j=1}^{J_\beta} \sum_{k=1}^{2^{jr}} jr 2^{2j(\beta+r/2)} \lambda_{j,k}(\varphi)^2< \infty.$

For this family of Gaussian processes, it is rather straightforward to verify that conditioning on a neighbourhood of $\beta$-smooth functions as in Step~1 of the deep Gaussian process prior construction is not a restrictive constraint. The next result shows that, with high probability, the process $X^\beta$ belong to the $\mB_{\infty,\infty,\beta}$-ball of radius $(1+K')\sqrt{2\log2},$ with $K'>\sqrt{3}.$ In view of Section~4.3.6 in~\cite{gine2016mathematical}, the Besov space $\mB_{\infty,\infty,\beta}$ contains the H\"older space $\mC_r^\beta$ for any $\beta>0,$ and they coincide whenever $\beta$ is not an integer. Moreover, $\mB_{\infty,\infty,\beta'}\subset \mC_r^\beta$ for all $\beta'>\beta.$

\begin{lemma}\label{lem.gp_conditioned}
    Let $X^\beta$ be the truncated wavelet process. Then, for any $K'>\sqrt{3},$ 
    \begin{align*}
        \P\left(\|X^\beta\|_{\infty,\infty,\beta} \leq (1+K')\sqrt{2\log 2} \right) \geq 1 - \frac{4}{2^{r{K'}^2} - 4}.
    \end{align*}
\end{lemma}
The proof is postponed to Appendix~\ref{sec.proofs_ex}. The probability in the latter display converges quickly to one. As an example consider $K'=2.$ Since $r\geq 1,$ the bound implies that more than $2/3$ of the simulated sample paths $\bu\mapsto X^\beta(\bu)$ lie in the Besov ball $\mB_{\infty,\infty,\beta}(3\sqrt{2\log 2}).$

The next lemma shows that for the truncated series expansion, near optimal posterior contraction rates can be achieved. The proof of the lemma is deferred to Appendix~\ref{sec.proofs_ex}. 

\begin{lemma}\label{lem.trunc_conditioned_ass2}
For $\{\widetilde G^{(\beta,r)}:\beta\in[\beta_{-},\beta_{+}]\}$ the family of truncated Gaussian processes $X^\beta$ there exist sequences $\eps_n(\eta) = C_1(\eta)(\log n)^{C_2(\eta)}\mathfrak{r}_n(\eta)$ such that
Assumption~\ref{ass.rate} holds.   
\end{lemma}

\subsection{Stationary process}
A zero-mean Gaussian process $X^\nu = (X^\nu(\bu))_{\bu\in[-1,1]^r}$ is called stationary if its covariance function can be represented by a spectral density measure $\nu$ on $\R^r$ as
\begin{align*}
	\E[X^\nu(\bu)X^\nu(\bu')] = \int_{\R^r} e^{-i(\bu-\bu')^\top\bxi} \nu(\bxi) d\bxi,
\end{align*} 
see Example~11.8 in~\cite{GVDV17}. We consider stationary Gaussian processes with radially decreasing spectral measures that have exponential moments, that is, $\int e^{c|\bxi|_2} \nu(\bxi) d\bxi <+\infty$ for some $c>0.$ Such processes have smooth sample paths thanks to Proposition~I.4 in~\cite{GVDV17}. An example is the square-exponential process with spectral measure $\nu(\bxi) = 2^{-r}\pi^{-r/2} e^{-|\bxi|_2^2/4}.$ For any $\varphi\in L^2(\nu),$ set $(H^\nu \varphi)(\bu) := \int_{\R^r} e^{i\bxi^\top\bu} \varphi(\bxi) \nu(\bxi) d\bxi.$ The RKHS of $X^\nu$ is given in Lemma~11.35 in~\cite{GVDV17} as $\bbH^\nu = \{H^\nu \varphi: \varphi\in L^2(\nu) \}$ with inner product $\langle H^\nu\varphi, H^\nu\varphi'\rangle_{\bbH^\nu} = \langle \varphi, \varphi' \rangle_{L^2(\nu)}.$

For every $\beta \in [\beta_{-},\beta_{+}],$ take $\widetilde G^{(\beta,r)}$ to be the rescaled process $X^{\nu}(a\cdot) = (X^{\nu}(a\bu))_{\bu\in[-1,1]^r}$ with scaling 
\begin{align}\label{eq.scaling}
    a = a(\beta,r) = n^{\frac{1}{2\beta+r}} (\log n)^{-\frac{1+r}{2\beta+r}}.
\end{align} 
The process $\widetilde G^{(\beta,r)}$ thus depends on the sample size $n.$ We prove the next result in Appendix~\ref{sec.proofs_ex}.

\begin{lemma}\label{lem.stat_ass2}
    For $\{\widetilde G^{(\beta,r)}:\beta\in[\beta_{-},\beta_{+}]\}$ the family of rescaled stationary processes $X^{\nu}(a\cdot)$, there exist sequences $\eps_n(\eta)=C_1(\eta)(\log n)^{C_2(\eta)}\mathfrak{r}_n(\eta)$ such that Assumption \ref{ass.rate} holds.   
\end{lemma}

\section{DGP priors, wide neural networks and regularization} \label{sec.disc_dgp_prior}

In this section, we explore similarities and differences between deep learning and the Bayesian analysis based on (deep) Gaussian process priors. Both methods are based on the likelihood. It is moreover known that standard random initialization schemes in deep learning converge to Gaussian processes in the wide limit. Since the initialization is crucial for the success of deep learning, this suggests that the initialization could act in a similar way as a Gaussian prior in the Bayesian world. Next to a proper initialization scheme, stability enhancing regularization techniques such as batch normalization are widely studied in deep learning and a comparison might help us to identify conditions that constraint the potentially wild behavior of deep Gaussian process priors. Below we investigate these aspects in more detail. 

It has been argued in the literature that Bayesian neural networks and regression with Gaussian process priors are intimately connected. In Bayesian neural networks, we generate a function valued prior distribution by using a neural network and drawing the network weights randomly. Recall that a neural network with a single hidden layer is called shallow, and a neural network with a large number of units in all hidden layers is called wide. If the network weights in a shallow and wide neural network are drawn i.i.d., and the scaling of the variances is such that the prior does not become degenerate, then, it has been argued in~\cite{Neal96} that the prior will converge in the wide limit to a Gaussian process prior and expressions for the covariance structure of the limiting process are known. One might be tempted to believe that for a deep neural network one should obtain a deep Gaussian process as a limit distribution. If the width of all hidden layers tends simultaneously to infinity, \cite{alex2018gaussian} proves that this is false and that one still obtains a Gaussian limit. The covariance of the limiting process is, however, more complicated and can be given via a recursion formula, where each step in the recursion describes the change of the covariance by a hidden layer. \cite{alex2018gaussian} shows moreover in a simulation study that Bayesian neural networks and Gaussian process priors with appropriate choice of the covariance structure behave indeed similarly.

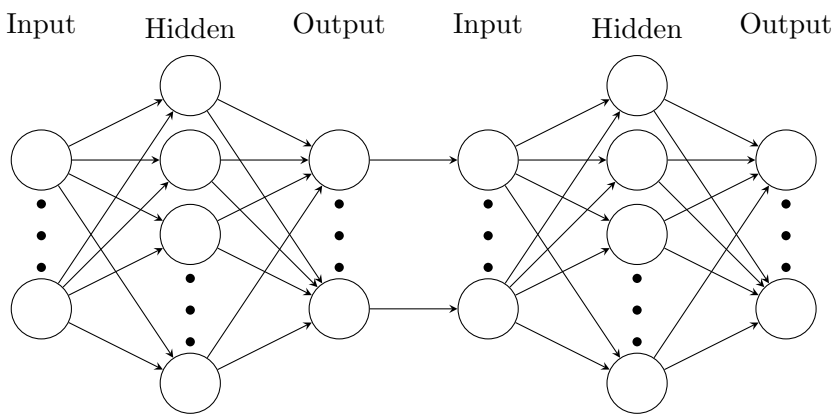
\begin{figure}[ht]
\hspace{-0.8cm}
    \centering
    \tikzset{%
  every neuron/.style={
    circle,
    draw,
    minimum size=0.8cm
  },
  neuron missing/.style={
    draw=none, 
    scale=3,
    text height=0.333cm,
    execute at begin node=\color{black}$\vdots$
  },
}

\begin{tikzpicture}[scale=0.9, x=1.1cm, y=1.1cm, >=stealth]

\foreach \m/\l [count=\y] in {1,missing,2}
  \node [every neuron/.try, neuron \m/.try] (input1-\m) at (0,1.5-\y) {};

\foreach \m/\l [count=\y] in {1,2,3,missing,4}
  \node [every neuron/.try, neuron \m/.try ] (hidden1-\m) at (2,2.5-\y) {};

\foreach \m [count=\y] in {1,missing,2}
  \node [every neuron/.try, neuron \m/.try ] (output1-\m) at (4,1.5-\y) {};
  
\foreach \m/\l [count=\y] in {1,missing,2}
  \node [every neuron/.try, neuron \m/.try] (input2-\m) at (6,1.5-\y) {};

\foreach \m/\l [count=\y] in {1,2,3,missing,4}
  \node [every neuron/.try, neuron \m/.try ] (hidden2-\m) at (8,2.5-\y) {};

\foreach \m [count=\y] in {1,missing,2}
  \node [every neuron/.try, neuron \m/.try ] (output2-\m) at (10,1.5-\y) {};

\node at (input1-1) {};
\node at (input1-2) {};

\foreach \l [count=\i] in {1,2,3}
  \node at (hidden1-\i) {};
\node at (hidden1-4) {};

\node at (output1-1) {};
\node at (output1-2) {};

\node at (input2-1) {};
\node at (input2-2) {};

\foreach \l [count=\i] in {1,2,3}
  \node at (hidden2-\i) {};
\node at (hidden2-4) {};
  
\node at (output2-1) {};
\node at (output2-2) {};

\foreach \i in {1,...,2}
   \foreach \j in {1,...,4}
     \draw [->] (input1-\i) -- (hidden1-\j);

\foreach \i in {1,...,4}
   \foreach \j in {1,...,2}
     \draw [->] (hidden1-\i) -- (output1-\j);
     
\draw [->] (output1-1) -- (input2-1);
\draw [->] (output1-2) -- (input2-2);
     
\foreach \i in {1,...,2}
   \foreach \j in {1,...,4}
     \draw [->] (input2-\i) -- (hidden2-\j);

 \foreach \i in {1,...,4}
   \foreach \j in {1,...,2}
     \draw [->] (hidden2-\i) -- (output2-\j);

\node [align=center, above] at (0,2) {Input};
\node [align=center, above] at (6,2) {Input};
\node [align=center, above] at (2,2) {Hidden};
\node [align=center, above] at (8,2) {Hidden};
\node [align=center, above] at (4,2) {Output};
\node [align=center, above] at (10,2) {Output};

\end{tikzpicture}
    \caption{Schematic stacking of two shallow neural networks.}
    \label{fig:nn_stack}
\end{figure}

It is conceivable that if one keeps the width of some hidden layers fixed and let the width of all other hidden layers tend to infinity, the Bayesian neural network prior will converge, if all variances are properly scaled, to a deep Gaussian process. By stacking for instance two shallow networks as indicated in Figure~\ref{fig:nn_stack} and making the first and last hidden layer wide, the limit is the composition of two Gaussian process and thus, a deep Gaussian process. In the hierarchical deep Gaussian prior construction in Section~\ref{sec.DGP}, we pick in a first step a prior on composition structures. For Bayesian neural networks this is comparable with selecting first a hyperprior on neural network architectures. 

Even more recently, \cite{peluchetti2019infinitely} as well as \cite{hayou2019impact} studied the behaviour of neural networks with random weights when both depth and width tend to infinity.

While the discussion so far indicates that Bayesian neural networks and Bayes with (deep) Gaussian process priors are similar methods, we now show that deep learning with randomly initialized network weights behaves similarly as a Bayes estimator with respect to a (deep) Gaussian process prior. The random network initialization means that the deep learning algorithm is initialized approximately by a (deep) Gaussian process. Since it is well-known that the initialization is crucial for the success of deep learning, this suggests that the initialization indeed acts as a prior. Denote by $-\ell$ the negative log-likelihood/cross entropy. Whereas in deep learning we fit a function by iteratively decreasing the cross-entropy using gradient descent method, the posterior is proportional to $\exp(\ell) \times$prior and concentrates on elements in the support of the prior with small cross entropy. Bayesian sampling methods such as stochastic gradient Langevin dynamics are closely related to noisy stochastic gradient descent, see \cite{WellingTeh} and Gibbs sampling of the posterior has a similar flavor as coordinate-wise descent methods for the cross-entropy. 

The only theoretical result that we are aware of examining the relationship between deep learning and Bayesian neural networks is \cite{NIPS2018_7372}. It proves that for a neural network prior with network weights drawn i.i.d. from a suitable spike-and-slab prior, the full posterior behaves similarly as the global minimum of the cross-entropy (that is, the empirical risk minimizer) based on sparsely connected deep neural networks.

As a last point we now compare stabilization techniques for deep Gaussian process priors and deep learning. In Step~1 of the deep Gaussian process prior construction, we have conditioned the individual Gaussian processes to map to $[-1,1]$ and to generate sample paths in small neighborhoods of a suitable H\"older ball. This induces regularity in the prior and avoids the wild behaviour of the composed sample paths due to bad realizations of individual components. We argue that this form of regularization has a similar flavor as batch normalization, cf. Section~8.7 in~\cite{Goodfellow2016}. The purpose of batch normalization is to avoid vanishing or exploding gradients due to the composition of several functions in deep neural networks. The main idea underlying batch normalization is to normalize the outputs from a fixed hidden layer in the neural network before they become the input of the next hidden layer. The normalization step is different from the conditioning proposed for the compositions of Gaussian processes. In fact, for batch normalization, mean and variance of the outputs are estimated based on a subsample and an affine transformation is applied such that the outputs are approximately centered and have variance one. One of the key differences is that this normalization invokes the distribution of the underlying design, while the conditioning proposed for deep Gaussian processes is independent of the distribution of the covariates. One suggestion that we can draw from this comparison is that instead of conditioning the processes to have sample paths in $[-1,1]$, it might also be interesting to apply the normalization $f\mapsto f(t)/\sup_{t\in [-1,1]^r} |f(t)|$ between any two compositions. This also ensures that the output maps to $[-1,1]$ and is closer to batch normalization. A data-dependent normalization of the prior cannot be incorporated in the fully Bayesian framework considered here and would result in an empirical Bayes method.

\subsection{Computational challenges of DGP priors} 

For complex tasks, sampling from the posterior is a notoriously hard computational problem. Several algorithmic breakthroughs resulted in computationally scalable procedures even in the particularly problematic case of high-dimensional model selection. While this is not the main focus of this work, we briefly discuss in this section aspects regarding the implementation of DGP priors and the  underlying computational challenges.\\

\noindent
\textbf{Model selection prior.} A major computational challenge is the large number of potential composition structures. A naive sampling method needs to visit all these composition structures. Below we describe several methods to reduce the number of candidate structures. 

Firstly, we argue that many models can in principle be omitted as they do not lead to faster contraction rates. As an example consider functions of the form $f=h_1\circ h_0$ for univariate functions $h_i$ having smoothness $\beta_i\leq 1,$ that is, $h_i\in \mC_1^{\beta_i}(K)$ for $i=0,1.$ From the definition of H\"older spaces, we can immediately infer that $f$ has H\"older smoothness $\beta_0\beta_1.$ This means that we can alternatively represent $f$ as a degenerate composition model $f=h_0$ with $q=0$ and $h_0$ a $\beta_0\beta_1$-smooth function. The optimal contraction rate $\frakr_n(\eta)$ defined in \eqref{eq.rate_minimax} is in both cases $n^{-\beta_0\beta_1/(2\beta_0\beta_1+1)}.$ Writing the function as a composition model, does therefore not lead to a faster contraction rate. The prior does consequently also not need to put any mass on this composition. 

For a general composition structure $f=h_q\circ h_{q-1}\circ \ldots \circ h_0,$ let $i^*$ be an index dominating the rate in the sense that $\frakr_n(\eta) = \max_{ i=0,\ldots, q} n^{-\beta_i\alpha_i/(2\beta_i\alpha_i+t_i)}=n^{-\beta_{i^*}\alpha_{i^*}/(2\beta_{i^*}\alpha_{i^*}+t_{i^*})}.$ In many cases, it can be checked that the rates for the two composition structures 
\begin{align*}
    f=\underbrace{h_q\circ \ldots \circ h_{i^*}}_{=:H_1}\circ \underbrace{h_{i^*-1}\circ \ldots \circ h_{0}}_{=:H_0} = H_1 \circ H_0
\end{align*}
are the same. In all these cases, it is therefore sufficient to consider the simpler composition structure $H_1\circ H_0$ only.

The next lemma provides conditions such that the number of compositions $q$ can be reduced by one.

\begin{lemma}\label{lem.redundant}
Suppose that $f$ is a function with composition structure $\eta= (q ,\bd ,\bt ,\mS, \bbeta)$ and assume that $\beta_+=K=1.$ If there exists an index $j\in \{1,\dots,q\}$ with $t_j=t_{j-1}=1,$ then, $f$ can also be written as a function with composition structure $\eta':= (q-1 ,\bd_{-j} ,\bt_{-j} ,\mS_{-j},\bbeta'),$ where $\bd_{-j}, \bt_{-j}, \mS_{-j}$ denote $\bd ,\bt ,\mS$ with entries $d_j, t_j, \mS_j$ removed, respectively, and $\bbeta':=(\beta_0,\dots,\beta_{j-2},\beta_{j-1}\beta_j, \beta_{j+1},\dots,\beta_q).$ Moreover the induced posterior contraction rates agree, that is, $\mathfrak{r}_n(\eta)=\mathfrak{r}_n(\eta').$
\end{lemma}

The discussion above calls for some notion of equivalence classes on composition structures. Two composition structures are then said to be equivalent, if one of them can be reduced to the other and both lead to the same rate $\mathfrak{r}_n(\eta).$ It is enough to assign prior mass to only one composition structure in each of the equivalence classes. The equivalence classes can be computed a priori.

Even for large sample sizes, it makes little difference whether the contraction rate is $n^{-1/10}$ or $n^{-1/11},$ say. While compositional structure can lead to much faster contraction rates, the biggest gains occur for structures with small $q$ and small $t_i.$ This indicates that if we are willing to loose a bit in the contraction rates, the number of candidate composition structures can be significantly further reduced. 

An alternative is to combine posterior sampling with a greedy model selection method. By first sampling from small composition structures, we can successively increase the complexity of the composition structures. In each iteration, we keep the composition structures with the most pronounced posterior concentration.

Additional information about the underlying composition structure can be incorporated as well. In nonparametric statistics, the (generalized) additive models and the single index model are two examples of widely studied structural assumptions imposed on the target function. Both constraints can also be rewritten as compositional constraint of the form $f=h_q\circ \ldots \circ h_0$ and correspond therefore to specific compositional structures, see \cite{SH19} for details. This means that if we suspect that those structural assumptions represent the true regression function, we can put higher (or all) prior weight on the corresponding composition structures.\\

\noindent
\textbf{Conditioning of Gaussian processes.} Step~1 of the deep Gaussian process prior construction requires to condition Gaussian processes to the unitary ball $\bbB_{\infty}(1)$ and some enlargement of a H\"older ball. Following Section~4 in~\cite{Swiler2020survey}, conditioning the sample paths to the unitary ball is a bound constraint where the range of the process is restricted to some compact interval. One proposed method is to add a warping function to the output of the Gaussian process in such a way to restrict the range to the desired interval, introducing nonlinear layers when dealing with compositions. Another solution involves discrete constraints using truncated Gaussian distributions. This requires an expensive rejection sampling procedure that has been recently made more efficient in~\cite{Ray2019efficient, lopezlopera2019gauss}. An alternative solution involves constrained maximum likelihood optimization. This approach has been proposed in~\cite{Pensoneault2020nonneg} to enforce non-negativity but can be easily extended to other inequality constraints. Another method, proposed by~\cite{maatouk2017gauss} is based on finite dimensional approximations of Gaussian processes.

Restricting the sample paths of Gaussian processes to neighborhoods of a H\"older ball is considerably more challenging from a computational point of view. Lemma~\ref{lem.gp_conditioned} shows that truncated wavelet processes, this constraint is not very restrictive.\\

%For these reasons, we believe it is more appropriate to solve this issue by finding processes that, at least with high probability, are already well-behaved.

%Throughout the following, we denote by $\{\widetilde G^{(\beta,r)}: \beta\in[\beta_{-},\beta_{+}], r\in\N\}$ such family of Gaussian processes. Here the quantities $0<\beta_{-}<1,$ $1<\beta_{+}<\infty$ can be chosen arbitrarily.

\noindent
\textbf{Variational Bayes.}  The variational Bayes approach introduced by~\cite{titsias09b} shows how to approximate the posterior distribution induced by a Gaussian process prior. \cite{nieman2021contr} gives contraction rates for such an approximation and shows that minimax rates are still achievable even when the KL divergence between the true posterior and the variational posterior does not vanish in the limit $n\to\infty$. This approach has been extended in~\cite{damianou2013deep}. Conditionally on a fixed composition structure $\eta,$ let $\Pi(\cdot|\bX,\bY,\eta)$ denote the posterior distribution. This is the same as putting all mass of the model selection prior $\pi$ on $\eta.$ Sampling from $\Pi(\cdot|\bX,\bY,\eta)$ does therefore not involve model selection. \cite{damianou2013deep} proposes a procedure to sample from an approximation of the distribution $\Pi(\cdot|\bX,\bY,\eta)$ arising from a DGP prior $\Pi(\cdot|\eta)$ obtained by composing square-exponential processes (radial basis function kernels). The same procedure has been reformulated in~\cite{CutajarDGP16} in terms of random features expansions.

To sample approximately from the full posterior $\Pi(\cdot|\bX,\bY)=\int \Pi(\cdot|\bX,\bY,\eta) \pi(\eta) \, d\eta$ for our hierarchical prior, one can first sample a compositional structure $\eta \sim \pi(\eta)$ and afterwards use the method from \cite{damianou2013deep} to approximate the distribution $\Pi(\cdot|\bX,\bY,\eta).$ Let $Q(\cdot|\bX,\bY,\eta)$ denote the distribution of the approximation and $Q(\cdot|\bX,\bY)$ be the induced approximation of the posterior $\Pi(\cdot|\bX,\bY)$. If for the variational Bayes step, we have some approximation guarantee of the form $\sup_{\eta} \|\Pi(\cdot|\bX,\bY,\eta)-Q(\cdot|\bX,\bY,\eta)\|_{TV}<\gamma$, then also $\|\Pi(\cdot|\bX,\bY)-Q(\cdot|\bX,\bY)\|_{TV}<\gamma$.

\section*{Acknowledgements}
We would like to thank Joris Bierkens and Isma\"el Castillo for helpful discussions. We are extremely grateful to the Associate Editor and three anonymous referees for their constructive comments and suggestions that resulted in a major improvement of the work. The research is part of the project \textit{Nonparametric Bayes for high-dimensional models: contraction, credible sets, computations} (with project number 613.001.604) financed by the Dutch Research Council (NWO).


\newpage

\appendix
\section{Proofs}

We provide here the proofs for the results presented in the main sections, together with auxiliary material.

\subsection{Proofs for Section~\ref{sec.model}}

\begin{lemma}\label{lem.holder.example}
If $f(x_1,\ldots,x_r)=g_1(x_1)\cdot \ldots \cdot g_r(x_r)$ with $g_1,\ldots,g_r\in \mC_1^\beta(K),$ then there exists a finite $K'$ such that $f\in \mC_r^\beta(K').$
\end{lemma}
\begin{proof}[Proof of Lemma~\ref{lem.holder.example}] Without loss of generality, assume $K\geq1$. For any $\balpha=(\alpha_1,\ldots,\alpha_r)\in\N^r$ and $\bx=(x_1,\ldots,x_r),$ we have $\partial^{\balpha}f(\bx)=\partial^{\alpha_1}g_1(x_1)\cdots\partial^{\alpha_r}g_r(x_r)$. By assumption, $\sum_{\alpha_j\leq\lfloor\beta\rfloor} \|\partial^{\alpha_j}g_j\|_{\infty}\leq K$ for all $j=1,\ldots,r$, thus
\begin{align}\label{eq.holder.ex.bound}
    \sum_{\bm{\alpha}:|\bm{\alpha}|\leq\lfloor\beta\rfloor} \|\partial^{\bm{\alpha}} f\|_{\infty} \leq \sum_{\bm{\alpha}:|\bm{\alpha}|\leq\lfloor\beta\rfloor} \left(\prod_{j=1}^r \|\partial^{\alpha_j}g_j\|_{\infty}\right) \leq \prod_{j=1}^r \left(\sum_{\alpha_j\leq\lfloor\beta\rfloor} \|\partial^{\alpha_j}g_j\|_{\infty}\right) \leq K^r.
\end{align}
For real numbers $a_1,\ldots,a_r,b_1,\ldots,b_r,$ we can expand the difference of the products by a telescoping sum, finding that $a_1\cdot \ldots \cdot a_r-b_1\cdot \ldots \cdot b_r=\sum_{j=1}^r (\prod_{s=1}^{j-1} a_s)(a_j-b_j)\prod_{t=j+1}^r b_t.$ Therefore, $|a_1\cdot \ldots \cdot a_r-b_1\cdot \ldots \cdot b_r|\leq (\prod_{s=1}^r |a_s|\vee |b_s|\vee 1)\sum_{j=1}^r|a_j-b_j|.$ Applying this to $\partial^{\balpha}f(\bx)=\partial^{\alpha_1}g_1(x_1)\cdots\partial^{\alpha_r}g_r(x_r)$ and $\partial^{\balpha}f(\by)=\partial^{\alpha_1}g_1(y_1)\cdots\partial^{\alpha_r}g_r(y_r)$ gives
\begin{align*}
    \sup_{\substack{\bx,\by\in [-1,1]^r \\ \bx\neq \by}} \frac{|\partial^{\bm{\alpha}} f(\bx) - \partial^{\bm{\alpha}} f(\by)|}{|\bx-\by|_{\infty}^{\beta-\lfloor\beta\rfloor}}\leq \left(\prod_{j=1}^r (\|\partial^{\alpha_j}g_j\|_{\infty}\vee1)\right) \sum_{j=1}^r \sup_{\substack{x_j,y_j\in [-1,1] \\ x_j\neq y_j}} \frac{|\partial^{\alpha_j}g_j(x_j) - \partial^{\alpha_j}g_j(x_j)|}{|x_j-y_j|^{\beta-\lfloor\beta\rfloor}}.
\end{align*}
For any $\balpha=(\alpha_1,\ldots,\alpha_r)$ such that $|\balpha|=\lfloor\beta\rfloor$, the function $\partial^{\alpha_j} g_j$ belongs to $\mC_r^{\beta-\alpha_j}(K)$ by assumption. Therefore, $\sup_{x\neq y\in [-1,1]} |\partial^{\alpha_j} g_j(x) - \partial^{\alpha_j} g_j(y)|/|x-y|^{\beta-\lfloor\beta\rfloor} \leq K$. With the bound in \eqref{eq.holder.ex.bound} and $K\geq 1,$ we obtain
\begin{align*}
    \sum_{\bm{\alpha}:|\bm{\alpha}|=\lfloor\beta\rfloor} &\sup_{\substack{\bx,\by\in [-1,1]^r \\ \bx\neq \by}} \frac{|\partial^{\bm{\alpha}} f(\bx) - \partial^{\bm{\alpha}} f(\by)|}{|\bx-\by|_{\infty}^{\beta-\lfloor\beta\rfloor}} \leq r K \cdot \sum_{\bm{\alpha}:|\bm{\alpha}|=\lfloor\beta\rfloor} \left(\prod_{j=1}^r (\|\partial^{\alpha_j}g_j\|_{\infty}\vee1)\right) \leq r K^{r+1}.
\end{align*}
The proof is complete by taking $K'\geq 2rK^r + 2^{\beta-\lfloor\beta\rfloor}rK^{r+1}$.

\end{proof}

\begin{proof}[Proof of Lemma~\ref{lem.A1_checked}]
Given $q,\bd$, the components of $\bt=(t_0,t_1,\ldots,t_q)$ are drawn independently from $\Unif\{1,\ldots,d_i\}$. Given $q,\bd, \bt$, $\mS=(\mS_0,\mS_1,\ldots,\mS_q)$ with $\mS_i=(\mS_{i1},\ldots,\mS_{id_{i+1}})$ are drawn independently from $\Unif\{1,\ldots,\tbinom{d_i}{t_i}\}$. By construction,
\begin{align*} 
    &\int \sqrt{\gamma(\eta)} \, d\eta \\
    &=\sum_{\lambda} \sqrt{(\beta_+-\beta_-)^{q+1}} \sqrt{\gamma(\lambda)} \\
    &= \sum_{q} \sqrt{(\beta_+-\beta_-)^{q+1}} \sqrt{\gamma(q)} \sum_{d_1,\ldots,d_q} \prod_{j=1}^q \sqrt{\gamma(d_j|q)}\sum_{t_0,\ldots,t_q} \sqrt{\prod_{i=0}^q \frac{1}{d_i} } \sum_{\mS_0,\ldots,\mS_q} \sqrt{\prod_{i=0}^q\prod_{j=1}^{d_{i+1}} \frac{1}{\binom{d_i}{t_i}}}.
\end{align*}
Because of the inequality $\tbinom{n}{k}\leq \sum_{q=0}^n \tbinom{n}{q}= 2^n,$ we have
\begin{align*} 
    \sum_{\mS_0,\ldots,\mS_q} \sqrt{\prod_{i=0}^q\prod_{j=1}^{d_{i+1}} \frac{1}{\binom{d_i}{t_i}}} = \prod_{i=0}^q\prod_{j=1}^{d_{i+1}} \binom{d_i}{t_i} \sqrt{\prod_{i=0}^q\prod_{j=1}^{d_{i+1}} \frac{1}{\binom{d_i}{t_i}}} = \sqrt{\prod_{i=0}^q \binom{d_i}{t_i}^{d_{i+1}}} \leq \sqrt{\prod_{i=0}^q 2^{d_id_{i+1}}}.
\end{align*}
The sum $\sum_{t_0,\ldots,t_q}$ is over $\prod_{i=0}^q d_i$ many terms. Using that $2^{d_id_{i+1}}\leq 2^{d_i^2/2}2^{d_{i+1}^2/2},$ the previous displays combined give
\begin{align} 
    &\int \sqrt{\gamma(\eta)} \, d\eta
    \leq \sqrt{d2^{d^2}}\sum_{q} \sqrt{(\beta_+-\beta_-)^{q+1}} \sqrt{\gamma(q)} \sum_{d_1,\ldots,d_q}  \prod_{j=1}^q \sqrt{\gamma(d_j|q) d_j 2^{d_j^2}}.
    \label{eq.99775}
\end{align}
Cauchy-Schwarz inequality gives
\begin{align*}
    \sum_{d_j} \sqrt{\gamma(d_j|q) d_j 2^{d_j^2}}
    \leq \sqrt{\sum_{d_j} \gamma(d_j|q) d_j^3 2^{d_j^2}}\sqrt{\sum_{d_j} \frac{1}{d_j^2}}
    =\sqrt{\frac{\pi^6}6\E_{d_1|q}[d_1^32^{d_1^2}]}=: \kappa
\end{align*}
Due to $\sum_{d_1,\ldots,d_q}  \prod_{j=1}^q=\prod_{j=1}^q\sum_{d_j},$ \eqref{eq.99775} combined with Cauchy-Schwarz inequality yields
\begin{align*}
    \int \sqrt{\gamma(\eta)} \, d\eta
    &\leq \sqrt{d2^{d^2}}\sum_{q} \sqrt{(\beta_+-\beta_-)^{q+1}} \sqrt{\gamma(q)} \kappa^q \\
    &\leq \sqrt{d2^{d^2}} \sqrt{\sum_{q} (\beta_+-\beta_-)^{2q+2}\gamma(q) e^{q}\kappa^{2q}}\sqrt{\sum_q e^{-q}} \\
    &= \sqrt{d2^{d^2}} \sqrt{\E_{q} \big[(\beta_+-\beta_-)^{2q+2}e^{q}\kappa^{2q}\big]}\frac{1}{\sqrt{e-1}} \\
    &<\infty,
\end{align*}
using for the last step that by assumption $\E_q[A^q]<\infty,$ for all $A>0.$ The proof is complete.
\end{proof}

\subsection{Proofs for Section~\ref{sec.main}}\label{sec.proofs_main}

{\bf Information geometry in the nonparametric regression model.} The following results are fairly standard in the nonparametric Bayes literature. As we are aiming for a self-contained presentation of the material, these facts are reproduced here. Let $P_f$ be the law of {\it one} observation $(\bX_i,Y_i).$ The Kullback-Leibler divergence in the nonparametric regression model is 
\begin{align*}
	\KL\big(P_f,P_g\big) = \frac{1}{2} \int (f(\bx)-g(\bx))^2 d\mu(\bx) \leq \frac 12 \|f-g\|_{L^\infty([-1,1]^d)}^2
\end{align*}
with $\mu$ the distribution of the covariates $\bX_1.$ Using that $\E_f[\log dP_f/dP_g]=\KL(P_f,P_g),$ $\Var(Z)\leq \E[Z^2],$ $\E_f[Y|\bX] = f(\bX)$ and $\E_f[(Y-f(\bX))^2|\bX] = 1,$ we also have that 
\begin{align*}
	V_2(P_f,P_g)&:=\E_f\Big[\Big|\log \frac{dP_f}{dP_g} -\KL(P_f,P_g)\Big|^2\Big] \\
	&\leq \E_f\Big[\Big|\log \frac{dP_f}{dP_g}\Big|^2\Big] \\
	&= \E_f\Big[ \Big(Y\big(f(\bX)-g(\bX)\big)-\frac{1}{2} f(\bX)^2 + \frac{1}{2} g(\bX)^2\Big)^2 \Big] \\
	&= \E_f\Big[ \Big(Y\big(f(\bX)-g(\bX)\big)-\frac{1}{2} \big(f(\bX)-g(\bX)\big) \big(f(\bX)+g(\bX)\big) \Big)^2 \Big] \\
	&= \E_f\Big[ \big(f(\bX)-g(\bX)\big)^2 \Big(Y-\frac{1}{2} \big(f(\bX)+g(\bX)\big) \Big)^2 \Big] \\
	&= \int \big(f(\bx)-g(\bx) \big)^2 +\frac{1}{4} \big(f(\bx)-g(\bx) \big)^4 \, d\mu(\bx). 
\end{align*}
In particular, for $\eps \leq 1,$ $\|f-g\|_\infty \leq \eps/2$ implies that $V_2(P_f,P_g)\leq \eps^2$ and therefore
\begin{align}\label{eq.B2}
    B_2(P_f,\eps) = \Big\{ g: \KL(P_f,P_g) < \eps^2, V_2(P_f,P_g) < \eps^2  \Big\}
    \supseteq \Big\{g:\|f-g\|_\infty\leq \frac{\eps}2 \Big\}. 
\end{align}

We derive posterior contraction rates for the Hellinger distance. This can then be related to the $\|\cdot\|_{L^2(\mu)}$-norm as explained below. Using the moment generating function of a standard normal distribution, the Hellinger distance for one observation $(\bX,Y)$ becomes
\begin{align*}
	d_H(P_f,P_g) = 1 - \int \sqrt{dP_f dP_g} 
	&= 1- \int \sqrt{dP_f/dP_g}\, dP_g \\
	&= 1- \E_g \Big[e^{\frac{1}{4} (Y -g(\bX))^2 - \frac{1}{4} (Y -f(\bX))^2}\Big] \\
	&= 1- \E\Big[\E_g\Big[e^{\frac{1}{2} (Y-g(\bX)) (f(\bX)-g(\bX))} \Big| \bX \Big] e^{-\frac{1}{4} (f(\bX)-g(\bX))^2} \Big] \\
	&= 1- \E\Big[e^{-\frac{1}{8} (f(\bX)-g(\bX))^2} \Big] \\
	&=1 - \int e^{-\frac{1}{8} (f(\bx)-g(\bx))^2} \, d\mu(\bx).
\end{align*}
Since $1-e^{-x}\leq x$ and $\mu$ is a probability measure, we have that $d_H(P_f,P_g) \leq \tfrac 18 \int (f(\bx)-g(\bx))^2 \, d\mu(\bx).$ Due to  $1-e^{-x} \geq e^{-x} x,$ we also find 
\begin{align}
	d_H(P_f,P_g)\geq \frac{e^{-Q^2/2}}8 \int \big(f(\bx)-g(\bx)\big)^2 \, d\mu(\bx), \ \ 
	\text{for all} \ f,g, \, \text{with} \  \|f\|_\infty,\|g\|_\infty \leq Q.
	\label{eq.Hell_lb}
\end{align}

By Proposition~D.8 in~\cite{GVDV17}, for any $f,g,$ there exists a test such that $\E_f\phi \leq \exp(-\tfrac n8 d_H(P_f,P_g)^2)$ and $\sup_{h:d_H(P_h,P_g)<d_H(P_f,P_g)/2} \E_h[1-\phi] \leq \exp(-\tfrac n8 d_H(P_f,P_g)^2).$ This means that for the Hellinger distance, the test condition in~(8.2) in~\cite{GVDV17} holds for $\xi=1/2.$

{\bf Function spaces.} 
The next result shows that the H\"older-balls defined in this paper are nested. 
\begin{lemma}~\label{lem.Di_embedding}
If $0<\beta'\leq\beta,$ then, for any positive integer $r$ and any $K>0,$ we have $\mC_r^{\beta}(K) \subseteq \mC_r^{\beta'}(K).$
\end{lemma}
\begin{proof}[Proof of Lemma~\ref{lem.Di_embedding}]
    If $\lfloor \beta'\rfloor = \lfloor \beta\rfloor$ the embedding follows from the definition of the H\"older-ball and the fact that $\sup_{\bx,\by \in [-1,1]^r} |\bx-\by|_\infty^{\beta-\beta'}=2^{\beta-\beta'}.$ If $\lfloor \beta'\rfloor < \lfloor \beta\rfloor,$ it remains to prove $\mC^{\beta}(K) \subseteq \mC^{\lfloor \beta' \rfloor+1}(K).$ This follows from first order Taylor expansion,
    \begin{align*}
        2\sum_{\bm{\alpha}:|\bm{\alpha}|=\lfloor\beta'\rfloor} \sup_{\substack{\bx,\by\in [-1,1]^r \\ \bx\neq \by}} \frac{|\partial^{\bm{\alpha}} f(\bx) - \partial^{\bm{\alpha}} f(\by)|}{|\bx-\by|_{\infty}}
        &\leq 2
        \sum_{\bm{\alpha}:|\bm{\alpha}|=\lfloor\beta'\rfloor}
        \big\| \big|\nabla (\partial^{\bm{\alpha}} f) \big|_1 \big\|_\infty \\
        &\leq 2r
        \sum_{\bm{\alpha}:|\bm{\alpha}|=\lfloor\beta'\rfloor+1}
        \big\| \partial^{\bm{\alpha}} f \big\|_\infty,
    \end{align*}
    and the definition of the H\"older-ball in~\eqref{def.holder_ball}.
\end{proof}

The following is a slight variation of Lemma~3 in \cite{SH19}.

\begin{lemma}
\label{lem.comp_approx}
Let $h_{ij}:[-1,1]^{t_i}\to [-1,1]$ be as in~\eqref{eq.f0}. Assume that, for some $K \geq 1$ and $\eta_i\geq0,$ $|h_{ij}(\bx)-h_{ij}(\by)|_\infty\leq \eta_i+K|\bx-\by|_\infty^{\beta_i\wedge 1}$ for all $\bx,\by\in [-1,1]^{t_i}.$ Then, for any functions $\widetilde h_i = (\widetilde h_{ij})_j^\top$ with $\widetilde h_{ij}:[-1,1]^{t_i}\rightarrow [-1,1],$
\begin{align*}
	&\big \| h_q \circ \ldots  \circ h_0 - \widetilde h_q \circ \ldots  \circ \widetilde h_0 \big\|_{L^\infty[-1,1]^d}
	\leq K^q
	\sum_{i=0}^q \eta_i^{\alpha_i}+ \big\| |h_i - \widetilde h_i |_\infty \big\|_{\infty}^{\alpha_i}.
\end{align*}
with $\alpha_i=\prod_{\ell = i+1}^q \beta_\ell \wedge 1.$
\end{lemma}
\begin{proof}[Proof of Lemma~\ref{lem.comp_approx}.]
We prove the assertion by induction over $q.$ For $q=0,$ the result is trivially true. Assume now that the statement is true for a positive integer $k.$ To show that the assertion also holds for $k+1,$ define $H_k = h_k \circ \ldots \circ h_0$ and $\widetilde H_k = \widetilde h_k \circ \ldots \circ \widetilde h_0.$ By triangle inequality, 
\begin{align*}
	&\big | h_{k+1}\circ H_k(\bx) - \widetilde h_{k+1}\circ\widetilde H_k(\bx)\big|_\infty \\
	&\leq | h_{k+1} \circ H_k(\bx) -  h_{k+1} \circ \widetilde H_k(\bx)\big|_\infty
	+ | h_{k+1} \circ \widetilde H_k(\bx) - \widetilde h_{k+1} \circ \widetilde H_k(\bx)\big|_\infty\\
	&\leq \eta_{k+1}+K \big| H_k(\bx) - \widetilde H_k(\bx) \big|_\infty^{\beta_{k+1}\wedge 1} + \| |h_{k+1}-\widetilde h_{k+1}|_\infty \|_{\infty}.
\end{align*}
Together with the induction hypothesis and the inequality $(y+z)^\alpha \leq y^\alpha + z^\alpha $ which holds for all $y, z \geq 0$ and all $\alpha \in [0,1],$ the induction step follows.
\end{proof}
 
The next result is a corollary of Theorem~8.9 in~\cite{GVDV17}.

\begin{lemma}
\label{lem.post_contr_outside_generic_set}
Denote the data by $\mD_n$ and the (generic) posterior by $\Pi(\cdot |\mD_n).$ Let $(A_n)_n$ be a sequence of events and $B_2(P_{f^*},\eps)$ as in~\eqref{eq.B2}. Assume that
\begin{align}
	e^{2n a_n^2}\frac{\Pi(A_n)}{\Pi(B_2(P_{f^*},a_n))} \xrightarrow{n\to\infty}0,
	\label{eq.post_contr_outside_generic_set_to_show}
\end{align}
for some positive sequence $(a_n)_n.$ Then,
\begin{align*}
	\E_{f^*}\big[\Pi(A_n|\mD_n)\big] \xrightarrow{n\to\infty}0,
\end{align*}
where $\E_{f^*}$ is the expectation with respect to $P_{f^*}.$
\end{lemma}

We now can prove Theorem~\ref{thm.mod_selection}.\\

\begin{proof}[Proof of Theorem~\ref{thm.mod_selection}] 
By definition~\eqref{def.post_model}, the quantity $\Pi\big(\eta\notin\mM_n(C)|\bX,\bY\big)$ denotes the posterior mass of the functions whose models are in the complement of $\mM_n(C).$ In view of Lemma~\ref{lem.post_contr_outside_generic_set}, it is sufficient to show condition~\eqref{eq.post_contr_outside_generic_set_to_show} for $A_n = \{\eta\notin\mM_n(C)\}=:\mM_n^c(C)$ and $a_n$ proportional to $\eps_n(\eta^*).$ We now prove that
\begin{align}
    e^{2n a_n^2} \frac{\int_{\mM_n^c(C)} \pi(\eta)\, d\eta}{\Pi(B_2(P_{f^*},a_n))} \rightarrow 0,
	\label{eq.mod_select_to_show}
\end{align}
for $a_n=4K^{q^*}(q^*+1)Q\eps_n(\eta^*)$ and $\Pi$ the deep Gaussian process prior. The next result deals with the lower bound on the denominator. For any hypercube $I,$ we introduce the notation $\diam(I):=\sup_{\bbeta,\bbeta'\in I}|\bbeta-\bbeta'|_\infty.$

\begin{lemma}\label{lem.small_ball_condition}
    Let $\Pi$ be a DGP prior satisfying the assumptions of Theorem~\ref{thm.mod_selection}. With $Q$ the universal constant from Assumption~\ref{ass.rate}, $R^* := 4K^{q^*}(q^*+1)Q$ and sufficiently large sample size $n,$ we have
    \begin{align*}
        \Pi\Big( B_2\big(P_{f^*}, 2R^*\eps_n(\eta^*)\big) \Big) \geq e^{- |\bd^*|_1 Q^2 n\eps_n(\eta^*)^2} \pi(\lambda^*,I_n^*),
    \end{align*}
    where $I_n^* := \{\bbeta=(\beta_0,\dots,\beta_{q^*}): \beta_i \in [\beta_i^*-b_n, \beta_i^*], \forall i \}$ and $b_n :=1/\log^2 n.$
\end{lemma}
\begin{proof}[Proof of Lemma~\ref{lem.small_ball_condition}]
    By construction~\eqref{eq.B2}, the set $B_2(P_{f^*}, 2R^*\eps_n(\eta^*))$ is a superset of $\{g:\|f^*-g\|_\infty\leq R^*\eps_n(\eta^*)\},$ so that
    \begin{align*}
        \Pi\Big( B_2\big(P_{f^*}, 2R^*\eps_n(\eta^*)\big) \Big) \geq \Pi\Big( f^* + \B_\infty\left(R^*\eps_n(\eta^*)\right) \Big).
    \end{align*}
    We then localize the probability in the latter display in the neighborhood $I_n^*$ around the true $\bbeta^*=(\beta_0^*,\dots,\beta_q^*).$ Since $\bbeta^* \in I(\lambda^*) = (\beta_-,\beta_+)^{q^*+1},$ we can always choose $n$ large enough such that $I_n^* \subseteq I(\lambda^*).$ With $f^* = h_{q^*}^* \circ \ldots \circ h_0^*$ and $R^* = 4K^{q^*}(q^*+1)Q,$
    \begin{align}
    \begin{split}\label{eq.sb_lb1}
    	\Pi&\Big(\big\{g:\|f^*-g\|_\infty \leq R^*\eps_n(\eta^*) \big\}\Big) \\
    	&\geq \int_{I_n^*} \P\Big(\big \|h_{q^*}^* \circ \ldots \circ h_0^*- G_{q^*}^{(\lambda^*,\bbeta)}\circ\ldots\circ G_0^{(\lambda^*,\bbeta)} \big\|_\infty \leq R^* \eps_n(\eta^*) \Big) \pi(\lambda^*,\bbeta) d\bbeta.
    \end{split}
    \end{align}
    Fix any $\bbeta\in I_n^*.$ Both $G_{ij}^{(\lambda^*,\bbeta)}=\overline G_{ij}^{(\lambda^*,\bbeta)}\circ (\cdot)_{\mathcal{S}_{ij}^*}$ and $h_{ij}^*=\overline h_{ij}^* \circ (\cdot)_{\mathcal{S}_{ij}^*}$ map $[-1,1]^{d_i^*}$ into $[-1,1].$ They also depend on the same subset of variables $\mS_{ij}^*.$ By construction, the process $\overline G_{ij}^{(\lambda^*,\bbeta)}$ is an independent copy of the conditioned Gaussian process $\widetilde G^{(\beta_i,t_i^*)}|\{\widetilde G^{(\beta_i,t_i^*)}\in \mD_i(\lambda^*,\bbeta,K)\}.$ The function $\overline h_{ij}^*$ belongs by definition to the space $\mC_{t_i^*}^{\beta_i^*}(K)$ and satisfies $|\overline h_{ij}^*(\bx)-\overline h_{ij}^*(\by)|\leq K|\bx-\by|_\infty^{\beta_i^*\wedge 1}$ for all $\bx,\by\in[-1,1]^{t_i^*}$ and all $i=0,\ldots,q^*;$ $j=1, \ldots,d_{i+1}^*.$ By Lemma~\ref{lem.comp_approx} with $\eta_i=0,$ we thus find
    \begin{align*}
    	\big\|f^*-G^{(\lambda^*,\bbeta)}\big\|_\infty \leq K^{q^*}\sum_{i=0}^{q^*}\ \max_{j=1,\ldots,d_{i+1}^*} \big\|\overline h_{ij}^*- \overline G_{ij}^{(\lambda^*,\bbeta)} \big\|_\infty^{\alpha_i^*},
    \end{align*}
    where $\alpha_i^*=\prod_{\ell = i+1}^{q^*} \beta_\ell^* \wedge 1.$ Since $\alpha_i^* \geq \alpha_i=\prod_{\ell = i+1}^{q^*} (\beta_\ell \wedge 1)$ for $\bbeta\in I_n^*,$ if $\|\overline h_{ij}^*- \overline G_{ij}^{(\lambda^*,\bbeta)}\|_\infty $ is smaller than one, the latter display is bounded above by
    \begin{align*}
    	\big\|f^*-G^{(\lambda^*,\bbeta)}\big\|_\infty \leq K^{q^*}\sum_{i=0}^{q^*}\ \max_{j=1,\ldots,d_{i+1}^*} \big\|\overline h_{ij}^*- \overline G_{ij}^{(\lambda^*,\bbeta)} \big\|_\infty^{\alpha_i}.
    \end{align*}
    
    Set $\delta_{in}:=\eps_n(\alpha_i,\beta_i,t_i^*)^{1/\alpha_i}.$ By Assumption~\ref{ass.rate} and the definition of $\eps_n(\eta)$ in~\eqref{eq.def_rate}, we have $\delta_{in} \leq (Q \eps_n(\eta^*))^{1/\alpha_i}$ and so $4\delta_{in} < 1$ since we are also assuming $\eps_n(\eta^*) < 1/(4Q)$. 
    If $\|\overline h_{ij}^*-\overline G_{ij}^{(\lambda^*,\bbeta)}\|_\infty \leq 4 \delta_{in}$ for all $i=0,\ldots,q^*$ and $j=1,\ldots,d_{i+1}^*$ then $\|f^*-G^{(\lambda^*,\bbeta)}\|_\infty \leq R^* \eps_n(\eta^*).$ Consequently, 
    \begin{align}
    \begin{split} \label{eq.sb_lb2}
        \P&\Big(\big \|h_{q^*}^* \circ \ldots \circ h_0^*- G_{q^*}^{(\lambda^*,\bbeta)}\circ\ldots\circ G_0^{(\lambda^*,\bbeta)} \big\|_\infty \leq R^* \eps_n(\eta^*) \Big) \\
        &\geq \prod_{i=0}^{q^*}\prod_{j=1}^{d_{i+1}^*} \P\Big( \big \| \overline h_{ij}^* - \overline G_{ij}^{(\lambda^*,\bbeta)} \big\|_\infty \leq 4 \delta_{in} \Big).
    \end{split}
    \end{align}
    We now lower bound the probabilities on the right hand side. If $\overline h_{ij}^*\in \mC_{t_i^*}^{\beta_i^*}(K),$ then, $(1-2\delta_{in})\overline h_{ij}^*\in \mC_{t_i^*}^{\beta_i^*}((1-2\delta_{in})K)$ and by the embedding property in Lemma~\ref{lem.Di_embedding}, we obtain $(1-2\delta_{in})\overline h_{ij}^*\in \mC_{t_i^*}^{{\beta_i}}(K).$
    
    When $\|(1-2\delta_{in}) \overline h_{ij}^*- \widetilde G^{(\beta_i,t_i^*)}\|_\infty\leq 2\delta_{in},$ the Gaussian process $\widetilde G^{(\beta_i,t_i^*)}$ is at most $2\delta_{in}$-away from $(1-2\delta_{in})\overline h_{ij}^*\in\mC_{t_i^*}^{\beta_i}(K).$ Consequently, $\{\|(1-2\delta_{in}) \overline h_{ij}^* - \widetilde G^{(\beta_i,t_i^*)}\|_{\infty}\leq 2\delta_{in}\}\subseteq \{\widetilde G^{(\beta_i,t_i^*)}\in \mD_i(\lambda^*,\bbeta,K) \},$ where the bound for the unitary ball follows from the triangle inequality. Moreover $\|(1-2\delta_{in}) \overline h_{ij}^*-\widetilde G^{(\beta_i,t_i^*)}\|_\infty\leq 2\delta_{in}$ implies $\| \overline h_{ij}^*- \widetilde G^{(\beta_i,t_i^*)}\|_\infty\leq 4\delta_{in}.$ The concentration function property in Lemma~I.28 in~\cite{GVDV17} gives
    \begin{align*}
        \P\Big(\big\| (1-\delta_{in}) \overline h_{ij}^*- \widetilde G^{(\beta_i,t_i^*)} \big\|_{\infty}\leq 2\delta_{in} \Big)
    	\geq \exp\Big(- \varphi^{(\beta_i,t_i^*,K)}\big(\delta_{in}\big) \Big).
    \end{align*}
    Together with the concentration function inequality in~\eqref{eq.def_eps_alpha}, we find
    \begin{align*}
    	\P\Big( \big \| \overline h_{ij}^* - \overline G_{ij}^{(\lambda^*,\bbeta)} \big\|_\infty \leq 4 \delta_{in} \Big) &\geq \frac{\P\big(\|(1-\delta_{in}\big) \overline h_{ij}^*- \widetilde G^{(\beta_i,t_i^*)}\|_{\infty}\leq 2\delta_{in})}{\P\big(\widetilde G^{(\beta_i,t_i^*)} \in \mD_i(\lambda^*,\bbeta,K)\big)} \\
    	&\geq \P\Big(\big\| (1-\delta_{in}) \overline h_{ij}^*- \widetilde G^{(\beta_i,t_i^*)} \big\|_{\infty}\leq 2\delta_{in} \Big) \\
    	&\geq \exp\Big(- \varphi^{(\beta_i,t_i^*,K)}\big(\delta_{in}\big) \Big)\\
    	&\geq \exp\Big( - n\eps_n(\alpha_i,\beta_i,t_i^*)^2 \Big) \\
    	&\geq \exp\Big( - Q^2 n\eps_n(\eta^*)^2 \Big),
    \end{align*}
    where $Q$ is the universal constant from Assumption~\ref{ass.rate}. With $|\bd^*|_1=1+\sum_{j=0}^{q^*} d_j^*,$ \eqref{eq.sb_lb1}, \eqref{eq.sb_lb2} and the previous display we recover the claim.
\end{proof}

The latter result shows that
\begin{align}
\begin{split}\label{eq.final_bound1}
    \Pi\Big( B_2\big(P_{f^*}, 4K^{q^*}(q^*+1)Q\eps_n(\eta^*)\big) \Big) 
	&\geq e^{- |\bd^*|_1 Q^2 n\eps_n(\eta^*)^2} \int_{I_n^*} \pi(\lambda^*,\bbeta) \, d\bbeta,
\end{split}
\end{align}
with $I_n^* = \{\bbeta=(\beta_0,\dots,\beta_{q^*}): \beta_i \in [\beta_i^*-1/\log^2 n, \beta_i^*], \forall i \}$. Recall that, by construction~\eqref{def.prior}, $\pi(\eta) \propto e^{-\Psi_n(\eta)} \gamma(\eta)$ with $\Psi_n(\eta) = n \eps_n(\eta)^2 + e^{e^{|\bd|_1}}.$ For any $\bbeta\in I_n^*,$ we have $\Psi_n(\lambda^*,\bbeta) \leq \Psi_n(\lambda^*,\bbeta^*),$ and Assumption~\ref{ass.prior} gives $\gamma(\lambda^*,\bbeta) = \gamma(\lambda^*)\gamma(\bbeta|\lambda^*)$ with $\gamma(\lambda^*)>0$ independent of $n$ and $\gamma(\cdot|\lambda^*)$ the uniform distribution over $I(\lambda^*) = [\beta_{-},\beta_{+}]^{q^*+1}.$ Thus, $\gamma(I_n^*|\lambda^*) = |I_n^*|/|I(\lambda^*)|$ and $|I_n^*| = (1/\log^2 n)^{q^*+1},$ so that by Assumption \ref{ass.rate}
\begin{align}\label{eq.prior_bound1}
    \begin{split}
	\frac{\pi(\eta)}{\int_{I_n^*} \pi(\lambda^*,\bbeta) \, d\bbeta} &\leq \frac{e^{Q^2n\eps_n(\eta^*)^2+e^{e^{|\bd^*|_1}}-\Psi_n(\eta)} \gamma(\eta)}{\gamma(\lambda^*) \gamma(I_n^*|\lambda^*)} \\ 
    &= \frac{|I(\lambda^*)|}{\gamma(\lambda^*)} e^{Q^2n\eps_n(\eta^*)^2+e^{e^{|\bd^*|_1}}-\Psi_n(\eta)} e^{2(q^*+1) \log\log n} \gamma(\eta).
    \end{split}
\end{align} 
The quantities $\gamma(\lambda^*),$ $|I(\lambda^*)| = (\beta_{+}-\beta_{-})^{q^*+1}$ and $\exp(e^{e^{|\bd^*|_1}})$ are constants independent of $n.$ We can finally verify condition~\eqref{eq.mod_select_to_show} by showing that, with $a_* = 4K^{q^*}(q^*+1)Q,$
\begin{align*}
	e^{2 a_*^2 n\eps_n(\eta^*)^2 } e^{Q^2n\eps_n(\eta^*)^2 + 2(q^*+1) \log\log n} \sum_{\lambda} \int_{\bbeta:(\lambda,\bbeta)\notin \mM_n(C)} e^{-\Psi_n(\eta)} \gamma(\eta) \, d\bbeta \rightarrow 0.
\end{align*}
By the lower bound in Assumption~\ref{ass.rate}~(i), we have $n\eps_n(\eta^*)^2 \geq n\frakr_n(\eta^*)^2 \gg 2(q^*+1)\log\log n,$ since the quantity $n\frakr_n(\eta^*)^2$ is a positive power of $n$ by definition~\eqref{eq.rate_minimax}. The complement of the set $\mM_n(C)$ is the union of $\{\eta:\eps_n(\eta) > C \eps_n(\eta^*)\}$ and $\{\eta:|\bd|_1>\log(2\log n)\}.$ Over these sets, by construction~\eqref{def.prior}, we have either $\Psi_n(\eta) > C^2n\eps_n(\eta^*)^2$ or $\Psi_n(\eta) > n^2.$ Therefore, the term $e^{-\Psi_n(\eta)}$ decays faster than either $e^{-C^2n\eps_n(\eta^*)^2}$ or $e^{-n^2}.$ In the first case, the latter display converges to zero for sufficiently large $C>0.$ In the second case, the latter display converges to zero since $\eps_n(\eta^*)<1$ and $n\eps_n(\eta^*)^2 < n \ll n^2.$ This completes the proof.
\end{proof}

We need some preliminary notation and results before proving Theorem~\ref{thm.post_contr}. 

\textbf{Entropy bounds.} Previous bounds for the metric entropy of H\"older-balls, e.g. Proposition~C.5 in~\cite{GVDV17}, are of the form $\log\mN\big(\delta, \mC_r^\beta(K), \|\cdot\|_\infty \big) \leq Q_1(\beta,r,K) \delta^{-r/\beta}$ for some constant $Q_1(\beta,r,K)$ that is hard to control. An exception is Theorem~8 in~\cite{bolley2010quantitative} that, however, only holds for $\beta\leq1.$ We derive an explicit bound on the constant $Q_1(\beta,r,K)$ for all $\beta>0.$ The proof is given in Appendix~\ref{sec.proofs_main_auxiliary}.
\begin{lemma}\label{lem.Holder_entropy_sharp}
    For any positive integer $r,$ any $\beta>0$ and $0<\delta<1,$ we have
    \begin{align*}
        \mN\left(\delta,C_r^\beta(K),\|\cdot\|_\infty\right) &\leq \left( \frac{4eK^2 r^\beta}{\delta} + 1\right)^{(\beta+1)^r} \left(2^{\beta+2} e Kr^\beta+1\right)^{4^r (\beta+1)^r r^r (2eK)^{\frac{r}{\beta}} \delta^{-\frac{r}{\beta}}} \\
        &\leq e^{Q_1(\beta,r,K)\delta^{-\frac{r}{\beta}}}
    \end{align*}
    with $Q_1(\beta,r,K) := (1+eK) 4^{r+1} (\beta+3)^{r+1} r^{r+1} (8eK^2)^{r/\beta}.$ For any $0<\alpha\leq 1$ and any sequence   $\delta_n \geq Q_1(\beta,r,K)^{\beta/(2\beta+r)}  n^{-\beta\alpha/(2\beta\alpha+r)},$ we also have
    \begin{align}\label{eq.delta_n_entropy}
        \log \mN\left(\delta_n^{1/\alpha},C_r^\beta(K),\|\cdot\|_\infty\right) &\leq n\delta_n^2.
    \end{align}
\end{lemma}

\textbf{Support of DGP prior and local complexity.} For any graph $\lambda=(q,\bd,\bt,\mS)$ and any $\bbeta \in [\beta_-,\beta_+]^{q+1},$ denote by $\Theta_n(\lambda,\bbeta, K)$ the space of functions $f:[-1,1]^d\to [-1,1]$ for which there exists a decomposition $f=h_q\circ \ldots \circ h_0$ such that $h_{ij}:[-1,1]^{d_i}\rightarrow [-1,1]$ and $\overline h_{ij} \in \mD_i(\lambda, \bbeta,K),$ for all $i=0,\ldots,q;$ $j=1,\ldots,d_{i+1}$ and $\mD_i(\lambda, \bbeta,K)$ as defined in~\eqref{eq.def.Di}. Differently speaking 
\begin{align}\label{eq.def_Theta}
    \Theta_n(\lambda,\bbeta,K) &:= \Theta_{q,n}(\lambda,\bbeta,K) \circ \cdots \circ \Theta_{0,n}(\lambda,\bbeta,K)
\end{align}
with
\begin{align*}
    \Theta_{i,n}(\lambda,\bbeta,K) &:= \Big\{h_i:[-1,1]^{d_i}\to[-1,1]^{d_{i+1}} : \overline h_{ij} \in \mD_i(\lambda,\bbeta,K), j=1,\ldots,d_{i+1} \Big\}.
\end{align*}
By construction, the support of the deep Gaussian process $G^{(\eta)}\sim\Pi(\cdot|\eta)$ is contained in $\Theta_n(\lambda,\bbeta, K).$ For a subset $B\subseteq [\beta_-,\beta_+]^{q+1}$ we also set $\Theta_n(\lambda,B, K):= \cup_{\bbeta\in B}\Theta_n(\lambda,\bbeta,K).$ The next lemma provides a bound for the covering number of $\Theta_n(\lambda,B, K).$ Recall that $\diam(B)=\sup_{\bbeta,\bbeta'\in B}|\bbeta-\bbeta'|_\infty.$ We postpone the proof to Appendix~\ref{sec.proofs_main_auxiliary}.

\begin{lemma}\label{lem.local_entropy}
    Suppose that Assumption~\ref{ass.rate} holds and let $\lambda$ be a graph such that $|\bd|_1\leq\log(2\log n).$ Let $B\subseteq[\beta_{-},\beta_{+}]^{q+1}$ with $\diam(B) \leq 1/\log^2 n.$ Then, with $R_n := 5Q(2\log n)^{1+\log K},$
    \begin{align*}
        \sup_{\bbeta\in B} \frac{\log\mN\left(R_n \eps_n(\lambda,\bbeta), \Theta_n(\lambda,B, K),\|\cdot\|_\infty\right)}{n\eps_n(\lambda,\bbeta)^2} &\leq \frac{R_n^2}{25}.
    \end{align*}
\end{lemma}

\begin{proof}[Proof of Theorem~\ref{thm.post_contr}]
For any $\rho>0,$ introduce the complement Hellinger ball $\mH^c(f^*,\rho) := \{f:d_H(P_f,P_{f^*}) > \rho\}.$ The convergence with respect to the Hellinger distance $d_H$ implies convergence in $L^2(\mu)$ thanks to \eqref{eq.Hell_lb}. As a consequence, we show that
    \begin{align*}
    	\sup_{f^* \in \mF(\eta^*,K)} \ \E_{f^*}\Big[\Pi\Big( \mH^c\big(f^*,L_n \eps_n(\eta^*)\big) \Big| (\bX,\bY)\Big)\Big] \rightarrow 0,
    \end{align*}
    with $L_n := MR_n,$ $R_n := 10QC(2\log n)^{1+\log K}$ and $M>0$ a sufficiently large universal constant to be determined. With $\mM_n(C) = \{\eta:\eps_n(\eta)\leq C\eps_n(\eta^*)\}\cap\{\eta: |\bd|_1\leq \log(2\log n)\}$ and the notation in~\eqref{def.post_model}, we denote by $\Pi(\cdot\cap\mM_n(C)|\bX,\bY)$ the contribution of the composition structures $\mM_n(C)$ to the posterior mass.

    We denote by $\mL_n(C)$ the set of graphs that are realized by some good composition structure, that is,
    \begin{align}\label{def.Ln_lambda}
        \mL_n(C) &:= \big\{\lambda\in\Lambda: \exists\bbeta\in I(\lambda),\ \eta=(\lambda,\bbeta)\in \mM_n(C) \big\}.
    \end{align}
    By construction, all graphs in $\mL_n(C)$ have bounded dimension $|\bd|_1\leq\log(2\log n)$ and so part~\ref{ass.rate_comparison} in Assumption~\ref{ass.rate} can be applied.  For any $\lambda\in \mL_n(C),$ partition $I(\lambda) = [\beta_{-},\beta_{+}]^{q+1}$ into hypercubes of diameter $1/\log^2 n$ and let $B_1(\lambda),\ldots,B_{N(\lambda)}(\lambda)$ be the $N(\lambda)$ blocks that contain at least one $\bbeta\in I(\lambda)$ that is realized by some composition structure in $\mM_n(C).$ The blocks may contain also values of $\bbeta$ for which $(\lambda,\bbeta)\notin\mM_n(C).$ The set $\mM_n(C)$ is contained in the enlargement
    \begin{align*}
        \widetilde\mM_n(C) := \bigcup_{\lambda\in\mL_n(C)} \bigcup_{k=1}^{N(\lambda)} \big(\lambda,B_k(\lambda) \big).
    \end{align*}
    Thanks to Theorem~\ref{thm.mod_selection} and the enlarged set of structures $\widetilde\mM_n(C),$ it is enough to show, for sufficiently large constants $M,C,$ 
    \begin{align}\label{eq.posterior_mass_to_bound_strong}
    	\sup_{f^* \in \mF(\eta^*,K)} \ \E_{f^*}\Big[\Pi\Big( \mH^c\big(f^*,MR_n \eps_n(\eta^*)\big) \cap \widetilde\mM_n(C)\Big| (\bX,\bY)\Big)\Big] \rightarrow 0.
    \end{align}
    
    Fix any $f^* \in \mF(\eta^*,K).$ Since there is no ambiguity, we shorten the notation to $\mH_n^c = \mH^c(f^*,MR_n \eps_n(\eta^*))$ and rewrite
    \begin{align*}
        \E_{f^*}\left[ \Pi\Big( \mH_n^c \cap \widetilde\mM_n(C)\Big| (\bX,\bY)\Big) \right] = \E_{f^*}\left[ \frac{ \int_{\mH_n^c} \Pi(df\cap\widetilde\mM_n(C)|\bX,\bY)}{ \int \Pi(df|\bX,\bY)} \right].
    \end{align*}
    We follow the steps of the proof of Theorem~8.14 in~\cite{GVDV17}. In their notation we use $\eps_n = R_n\eps_n(\eta^*)$ and $\xi=1/2$ for contraction with respect to Hellinger loss. Set
    \begin{align*}
    	A_n^* = \bigg\{ \int \Pi(df|\bX,\bY) \geq \Pi\Big(B_2\big(f^*,R_n\eps_n(\eta^*)\big) \Big) e^{-2R_n^2 n \eps_n(\eta^*)^2} \bigg\}.
    \end{align*}
    Then $\P_{f^*}(A_n^*)$ tends to $1,$ thanks to Lemma~8.10 in \cite{GVDV17} applied with $D=1.$ Since $1 = \Ind(A_n^*) + \Ind(A_n^{*,c}),$ we have
    \begin{align*}
        \E_{f^*}\left[ \Pi\Big( \mH_n^c \cap \widetilde\mM_n(C)\Big| (\bX,\bY)\Big) \right] \leq \P_{f^*}(A_n^{*,c}) 
        + \E_{f^*}\left[ \Ind(A_n^*) \frac{ \int_{\mH_n^c} \Pi(df\cap \widetilde\mM_n(C)|\bX,\bY)}{ \int \Pi(df|\bX,\bY)} \right]
    \end{align*}
    and $\P_{f^*}(A_n^{*,c}) \to 0$ when $n\to+\infty.$ It remains to show that the second terms on the right side tends to zero. 
    
    Let $\phi_{n,k}(\lambda)$ be arbitrary statistical tests to be chosen later. Test are to be understood as $\phi_{n,k}(\lambda) = \phi_{n,k}(\lambda)(\bX,\bY)$ measurable functions of the sample $(\bX,\bY),$ taking values in $[0,1].$ Then, $1 = \phi_{n,k}(\lambda) + (1-\phi_{n,k}(\lambda)).$ 
    Using the definition of $\Pi(df\cap \widetilde\mM_n(C)|\bX,\bY)$ and Fubini's theorem, we find
    \begin{align}\label{eq.T1_T2}
        \E_{f^*}\left[ \Pi\Big( \mH_n^c \cap \widetilde\mM_n(C)\Big| (\bX,\bY)\Big) \right] \leq \P_{f^*}(A_n^{*,c}) 
        + \E_{f^*}\left[ T_1 + T_2 \right],
    \end{align}
    where
    \begin{align*}
        T_1 &:= \Ind(A_n^*) \frac{ \sum_{\lambda\in\mL_n(C)} \sum_{k=1}^{N(\lambda)} \phi_{n,k}(\lambda) \int_{B_k(\lambda)}  \pi(\lambda,\bbeta) \left(\int_{ \mH_n^c} \frac{p_f}{p_{f^*}}(\bX,\bY)  \Pi(df|\lambda,\bbeta) \right) d\bbeta}{\int \Pi(df|\bX,\bY)}, \\
        T_2 &:= \Ind(A_n^*) \frac{ \sum_{\lambda\in\mL_n(C)} \sum_{k=1}^{N(\lambda)} \int_{B_k(\lambda)}  \pi(\lambda,\bbeta) \left(\int_{ \mH_n^c} (1-\phi_{n,k}(\lambda)) \frac{p_f}{p_{f^*}}(\bX,\bY)  \Pi(df|\lambda,\bbeta) \right) d\bbeta}{\int \Pi(df|\bX,\bY)}.
    \end{align*}
    We bound $T_1$ by using $\Ind(A_n^*) \leq 1,$ together with
    \begin{align*}
        \frac{ \int_{B_k(\lambda)}  \pi(\lambda,\bbeta) \left(\int_{ \mH_n^c} \frac{p_f}{p_{f^*}}(\bX,\bY)  \Pi(df|\lambda,\bbeta) \right) d\bbeta}{\int \Pi(df|\bX,\bY)} \leq 1,
    \end{align*}
    so that
    \begin{align}\label{eq.T_1_bound}
        \E_{f^*}\left[ T_1 \right] &\leq \sum_{\lambda\in\mL_n(C)} \sum_{k=1}^{N(\lambda)} \E_{f^*}\left[ \phi_{n,k}(\lambda) \right].
    \end{align}
    
    We bound $T_2$ using the definition of $A_n^*,$ and obtain
    \begin{align*}
        T_2 \leq \frac{ \sum_{\lambda\in\mL_n(C)} \sum_{k=1}^{N(\lambda)} \int_{B_k(\lambda)}  \pi(\lambda,\bbeta) \left(\int_{ \mH_n^c} (1-\phi_{n,k}(\lambda)) \frac{p_f}{p_{f^*}}(\bX,\bY) \Pi(df|\lambda,\bbeta) \right) d\bbeta}{\Pi\left(B_2\big(f^*,R_n\eps_n(\eta^*)\big) \right) e^{-2R_n^2 n \eps_n(\eta^*)^2}}.
    \end{align*}
    For large $n,$ $R_n = 10QC(2\log n)^{1+\log K} \geq 4QK^{q^*}(q^*+1).$ By Lemma~\ref{lem.small_ball_condition}, we find
    \begin{align}\label{eq.small_ball_first}
        \Pi\left(B_2\big(f^*,R_n\eps_n(\eta^*)\big) \right) \geq e^{-R_n^2 n\eps_n(\eta^*)^2} \pi(\lambda^*,I_n^*),
    \end{align}
    with $I_n^* := \{\bbeta=(\beta_0,\dots,\beta_{q^*}): \beta_i \in [\beta_i^*-b_n, \beta_i^*], \forall i \}$ and $b_n := 1/\log^2 n.$ By the construction of the prior $\pi$ in~\eqref{def.prior}, the denominator term $e^{-\Psi_n(\eta)}$ is bounded above by $1$ and $\int_\Omega \gamma(\eta)\, d\eta = 1,$ thus 
    \begin{align*}
        \pi(\lambda^*,I_n^*) &= \frac{\int_{I_n^*} e^{-\Psi_n(\lambda^*,\bbeta)} \gamma(\lambda^*,\bbeta) d\bbeta}{\int_\Omega e^{-\Psi_n(\eta)} \gamma(\eta) \, d\eta} \geq \int_{I_n^*} e^{-\Psi_n(\lambda^*,\bbeta)} \gamma(\lambda^*,\bbeta) d\bbeta.
    \end{align*}
    Furthermore, by condition~\ref{ass.rate_comparison} in Assumption~\ref{ass.rate}, we have $\eps_n(\lambda^*,\bbeta) \leq Q\eps_n(\lambda^*,\bbeta^*)$ for any $\bbeta \in I_n^*,$ which gives 
    \begin{align*}
        \pi(\lambda^*,I_n^*) \geq \exp(-e^{e^{|\bd^*|_1}})\, e^{-Q^2n\eps_n(\eta^*)^2} \gamma(\lambda^*,I_n^*),
    \end{align*}
    with proportionality constant $\exp(-e^{e^{|\bd^*|_1}}) > 0$ independent of $n.$ Again by construction, the measure $\gamma(\lambda^*,I_n^*)$ can be split into the product $\gamma(\lambda^*)\gamma(I_n^*|\lambda^*)$ where $\gamma(\cdot|\lambda^*)$ is the uniform measure on $I(\lambda^*) = [\beta_{-},\beta_{+}]^{q^*+1}$ and the quantity $\gamma(\lambda^*) > 0$ is a constant independent of $n.$ Thus, with $c(\eta^*) := \exp(-e^{e^{|\bd^*|_1}})\gamma(\lambda^*)/|I(\lambda^*)|,$ we obtain $\pi(\lambda^*,I_n^*) \geq c(\eta^*) |I_n^*| e^{-Q^2n\eps_n(\eta^*)^2}.$ In the discussion after~\eqref{eq.prior_bound1} we have shown that $n\eps_n(\eta^*)^2$ is a positive power of $n$ and so, for a large enough $n$ depending only on $\eta^*,$ one has $|I_n^*| = (1/\log^2 n)^{q^*+1} > e^{-n\eps_n(\eta^*)^2}.$ This results in
    \begin{align}\label{eq.hyperprior_lowbound}
        \pi(\lambda^*,I_n^*) \geq c(\eta^*) e^{-(Q^2+1)n\eps_n(\eta^*)^2} \geq c(\eta^*) e^{-R_n^2n\eps_n(\eta^*)^2}.
    \end{align}
    By putting together the small-ball probability bound~\eqref{eq.small_ball_first} and the hyperprior bound~\eqref{eq.hyperprior_lowbound}, we recover
    \begin{align}\label{eq.T_2_bound}
        T_2 \leq c(\eta^*)^{-1} \frac{ \sum_{\lambda\in\mL_n(C)} \sum_{k=1}^{N(\lambda)} \int_{B_k(\lambda)}  \pi(\lambda,\bbeta) \left(\int_{ \mH_n^c} (1-\phi_{n,k}(\lambda)) \frac{p_f}{p_{f^*}}(\bX,\bY) \Pi(df|\lambda,\bbeta) \right) d\bbeta}{e^{-4R_n^2 n \eps_n(\eta^*)^2}},
    \end{align}
    with proportionality constant $c(\eta^*)^{-1}$ independent of $n$ and depending only on $\eta^*,\ \beta_{-},\ \beta_{+}.$ We now bound the numerator of $T_2$ using Fubini's theorem and the inequality $\E_{f^*}[(1 - \phi_{n,k}(\lambda)) (p_f/p_{f_*})(\bX,\bY)] \leq \E_f[1 - \phi_{n,k}(\lambda)],$
    \begin{align}
        \E_{f^*}&\left[ \sum_{\lambda\in\mL_n(C)} \sum_{k=1}^{N(\lambda)} \int_{B_k(\lambda)}  \pi(\lambda,\bbeta) \left(\int_{ \mH_n^c} (1-\phi_{n,k}(\lambda)) \frac{p_f}{p_{f^*}}(\bX,\bY) \Pi(df|\lambda,\bbeta) \right) d\bbeta \right] \notag \\
        &= \sum_{\lambda\in\mL_n(C)} \sum_{k=1}^{N(\lambda)} \int_{B_k(\lambda)}  \pi(\lambda,\bbeta) \left(\int_{ \mH_n^c} \E_{f^*}\left[ (1 - \phi_{n,k}(\lambda)) \frac{p_f}{p_{f^*}}(\bX,\bY) \right] \Pi(df|\lambda,\bbeta)\right) d\bbeta \notag \\
        &\leq \sum_{\lambda\in\mL_n(C)} \sum_{k=1}^{N(\lambda)} \int_{B_k(\lambda)}  \pi(\lambda,\bbeta) \int_{ \mH_n^c} \E_{f}\left[ (1 - \phi_{n,k}(\lambda)) \right] \Pi(df|\lambda,\bbeta) d\bbeta. \label{eq.T_2_bound_fubini}
    \end{align}
    
    With the supports $\Theta_n(\lambda,\bbeta,K)$ in \eqref{eq.def_Theta} and any $k=1,\ldots,N(\lambda),$ consider $\Theta_n(\lambda,B_k(\lambda),K) = \cup_{\bbeta\in B_k(\lambda)} \Theta_n(\lambda,\bbeta,K).$ Now, for any fixed $\lambda,k$ choose tests $\phi_{n,k}(\lambda)$ according to Theorem D.5 in \cite{GVDV17}, so that for $f\in \Theta_n(\lambda,B_k(\lambda),K)\cap\mH^c(f^*,MR_n \eps_n(\eta^*))$ we have, for some universal constant $\widetilde K>0,$
    \begin{align*}
    	\E_{f^*} [\phi_{n,k}(\lambda)] &\leq c_k(\lambda) \mN\left(\frac{R_n\eps_n(\eta^*)}{2}, \Theta_n(\lambda,B_k(\lambda),K),\|\cdot\|_{\infty}\right)  \frac{e^{- \widetilde K M^2 R_n^2 n \eps_n(\eta^*)^2}}{1-e^{-\widetilde K M^2 R_n^2 n \eps_n(\eta^*)^2}},  \\
    	\E_f[1-\phi_{n,k}(\lambda)] &\leq c_k(\lambda)^{-1} e^{-\widetilde K M^2 R_n^2 n \eps_n(\eta^*)^2},
    \end{align*}
     with choice of coefficients
    \begin{align*}
        c_k(\lambda)^2 := \frac{\pi(\lambda,B_k(\lambda))}{\mN\left(\frac{R_n\eps_n(\eta^*)}{2}, \Theta_n(\lambda,B_k(\lambda),K),\|\cdot\|_{\infty}\right) }.
    \end{align*}
    Let us denote by $\rho_k(\lambda)$ the local complexities
    \begin{align*}
        \rho_k(\lambda) := \sqrt{\pi(\lambda,B_k(\lambda)) } \cdot \sqrt{\mN\left(\frac{R_n\eps_n(\eta^*)}{2}, \Theta_n(\lambda,B_k(\lambda),K),\|\cdot\|_{\infty}\right)}.
    \end{align*}
    Combining this with the bound on $T_1$ in~\eqref{eq.T_1_bound} and the bounds on $T_2$ in~\eqref{eq.T_2_bound}--\eqref{eq.T_2_bound_fubini}, gives
    \begin{align*}
        \E_{f^*}\left[ T_1 \right] &\leq \frac{e^{-\widetilde K M^2 R_n^2 n \eps_n(\eta^*)^2}}{1-e^{-\widetilde K M^2 R_n^2 n \eps_n(\eta^*)^2}} \sum_{\lambda\in\mL_n(C)} \sum_{k=1}^{N(\lambda)} \rho_k(\lambda), \\
        \E_{f^*}\left[ T_2 \right] &\leq c(\eta^*)^{-1} e^{(4 - \widetilde K M^2)R_n^2 n \eps_n(\eta^*)^2} \sum_{\lambda\in\mL_n(C)} \sum_{k=1}^{N(\lambda)} \rho_k(\lambda).
    \end{align*}
    It remains to show that both expectations in the latter display tend to zero when $n\to+\infty.$ 
    
    Since $M > 0$ can be chosen arbitrarily large, we choose it in such a way that $\widetilde K M^2 > 5$ and the proof is complete if we can show that
    \begin{align}\label{eq.final_to_show}
        \sum_{\lambda\in\mL_n(C)} \sum_{k=1}^{N(\lambda)} \rho_k(\lambda) \lesssim e^{R_n^2 n\eps_n(\eta^*)^2},
    \end{align}
    for some proportionality constant independent of $n.$ Fix $\lambda\in\mL_n(C)$ and $k=1,\ldots,N(\lambda).$ By construction, there exists $\bbeta\in B_k(\lambda)$ such that $(\lambda,\bbeta)\in\mM_n(C).$ Since $R_n = 10QC (2\log n)^{1+\log K},$ using $C\eps_n(\eta^*) \geq \eps_n(\eta)$ together with Lemma~\ref{lem.local_entropy}, we find
    \begin{align*}
        \log\mN\left( \frac{R_n\eps_n(\eta^*)}{2}, \Theta_n(\lambda,B, K),\|\cdot\|_\infty\right) &\leq \frac{R_n^2}{25} n\eps_n(\eta^*)^2.
    \end{align*}
    This results in 
    \begin{align*}
        \sum_{\lambda\in\mL_n(C)} \sum_{k=1}^{N(\lambda)} \rho_k(\lambda) &\leq \exp\left(\frac{1}{50} R_n^2 n\eps_n(\eta^*)^2\right) \sum_{\lambda\in\mL_n(C)} \sum_{k=1}^{N(\lambda)} \sqrt{\pi(\lambda,B_k(\lambda)) } \\
        &= \exp\left(\frac{1}{50} R_n^2 n\eps_n(\eta^*)^2\right) z_n^{-\frac{1}{2}} \sum_{\lambda\in \mL_n(C)} \sum_{k=1}^{N(\lambda)} \sqrt{{\gamma}(\lambda,B_k(\lambda))},
    \end{align*}
    where $z_n = \sum_\lambda\int_{I(\lambda)} e^{-\Psi_n(\lambda,\bbeta)} \gamma(\lambda,\bbeta) d\bbeta$ is the normalization term in~\eqref{def.prior}. By the localization argument in~\eqref{eq.hyperprior_lowbound}, we know that $z_n \geq \pi(\lambda^*,I_n^*) \gtrsim e^{-R_n^2 n\eps_n(\eta^*)^2},$ with proportionality constant independent of $n.$ Thus,
    \begin{align*}
        \sum_{\lambda\in\mL_n(C)} \sum_{k=1}^{N(\lambda)} \rho_k(\lambda) &\lesssim \exp\left(\frac{26}{50} R_n^2 n\eps_n(\eta^*)^2\right) \sum_{\lambda\in \mL_n(C)} \sum_{k=1}^{N(\lambda)} \sqrt{{\gamma}(\lambda,B_k(\lambda))}.
    \end{align*}
    Since $R_n = 10QC (2\log n)^{1+\log K} \gg 1,$ it is sufficient for~\eqref{eq.final_to_show} that
    \begin{align}\label{eq.rieman_sum_stronger}
        \sum_{\lambda\in \mL_n(C)} \sum_{k=1}^{N(\lambda)} \sqrt{{\gamma}(\lambda,B_k(\lambda))} \lesssim e^{\frac{1}{2} n\eps_n(\eta^*)^2},
    \end{align}
    for some proportionality constant that can be chosen independent of $n.$ To see this, observe that by Assumption~\ref{ass.prior}, we have $\gamma(\lambda,\bbeta) = \gamma(\lambda)\gamma(\bbeta|\lambda)$ with $\gamma(\cdot|\lambda)$ the uniform distribution over $I(\lambda).$ Thus
    \begin{align*}
        \sum_{\lambda\in \mL_n(C)} \sum_{k=1}^{N(\lambda)} \sqrt{{\gamma}(\lambda,B_k(\lambda))} &= \sum_{\lambda\in \mL_n(C)} \sum_{k=1}^{N(\lambda)} \sqrt{|B_k(\lambda)|} \sqrt{{\gamma}(\lambda,\bbeta_k(\lambda))},
    \end{align*}
    where $\bbeta_k(\lambda)$ is the center point of the hypercube $B_k(\lambda).$ Since $|B_k(\lambda)| = (1/\log^2 n)^{q+1}$ and $(q+1) \leq |\bd|_1 \leq \log(2\log n)$, we have $\log (|B_k(\lambda)|^{-1}) = 2(q+1)(\log n) \leq 4\log^2 n \ll n\eps_n(\eta^*)^2.$ Therefore, for a sufficiently large $n$ depending only on $\eta^*,$ we find $|B_k(\lambda)|^{-1} \leq e^{n\eps_n(\eta^*)^2}.$ The above discussion yields the following bound on the latter display,
    \begin{align*}
        \sum_{\lambda\in \mL_n(C)} \sum_{k=1}^{N(\lambda)} \sqrt{{\gamma}(\lambda,B_k(\lambda))} &= \sum_{\lambda\in \mL_n(C)} \sum_{k=1}^{N(\lambda)} \frac{1}{\sqrt{|B_k(\lambda)|}} |B_k(\lambda)| \sqrt{{\gamma}(\lambda,\bbeta_k(\lambda))} \\
        &\leq e^{\frac{1}{2} n\eps_n(\eta^*)^2} \sum_{\lambda\in \mL_n(C)} \sum_{k=1}^{N(\lambda)} |B_k(\lambda)| \sqrt{{\gamma}(\lambda,\bbeta_k(\lambda))}.
    \end{align*}
    This is enough to obtain~\eqref{eq.rieman_sum_stronger} since, by Assumption~\ref{ass.prior},
    \begin{align*}
        \sum_{\lambda\in \mL_n(C)} \sum_{k=1}^{N(\lambda)} |B_k(\lambda)| \sqrt{{\gamma}(\lambda,\bbeta_k(\lambda))} \leq \sum_{\lambda\in \Lambda} \int_{I(\lambda)} \sqrt{{\gamma}(\lambda,\bbeta)}d\bbeta = \int_\Omega \sqrt{\gamma(\eta)}\, d\eta
    \end{align*}
    is a finite constant independent of $n.$ 
    
    Since all bounds are independent of the particular choice of $f^*$ and only depend on the function class $\mF(\eta^*,K),$ this concludes the proof of the uniform statement~\eqref{eq.posterior_mass_to_bound_strong}.
\end{proof}

\subsubsection{Proofs of auxiliary results}\label{sec.proofs_main_auxiliary}

\begin{proof}[Proof of Lemma~\ref{lem.Holder_entropy_sharp}]
    We follow the proof of Theorem~2.7.1 in~\cite{van1996weak} and provide explicit expressions for all constants. We start by covering the interval $[-1,1]^r$ with a grid of width $\tau = (\delta/c(\beta))^{1/\beta},$ where $c(\beta) := e r^\beta + 2K.$ 
    The grid consists of $M$ points $\bx_1,\ldots,\bx_M$ with 
    \begin{align}\label{eq.bound_on_M}
        M \leq \frac{vol([-2,2]^r)}{\tau^r} = 4^r c(\beta)^{\frac{r}{\beta}} \delta^{-\frac{r}{\beta}}.
    \end{align}
    For any $h\in\mC_r^\beta(K)$ and any $\balpha = (\alpha_1,\ldots,\alpha_r)\in\N^r$ with $|\balpha|_1=\alpha_1+\ldots+\alpha_r \leq \lfloor\beta\rfloor,$ set
    \begin{align}\label{eq.def_A_alpha}
        A^{\balpha}h := \left(\left\lfloor \frac{\partial^{\balpha}h(\bx_1)}{\tau^{\beta-|\balpha|_1}} \right\rfloor,\ldots, \left\lfloor \frac{\partial^{\balpha}h(\bx_M)}{\tau^{\beta-|\balpha|_1}} \right\rfloor \right).
    \end{align}
    The vector $\tau^{\beta-|\balpha|_1} A^{\balpha}h$ consists of the values $\partial^{\balpha}h(\bx_i)$ discretized on a grid of mesh-width $\tau^{\beta-|\balpha|_1}.$
    Since $\partial^{\balpha}h \in \mC_r^{\beta-|\balpha|_1}(K)$ and $\tau < 1$ by construction, the entries of the vector in the latter display are integers bounded in absolute value by
    \begin{align}\label{eq.bound_on_coeff}
        \left\lfloor \frac{|\partial^{\balpha}h(\bx_i)|}{\tau^{\beta-|\balpha|_1}} \right\rfloor \leq \left\lfloor \frac{K}{\tau^{\beta-|\balpha|_1}} \right\rfloor \leq \left\lfloor \frac{K}{\tau^{\beta}} \right\rfloor.
    \end{align}
    Let $h,\widetilde h \in \mC_r^\beta(K)$ be two functions such that $A^{\balpha}h = A^{\balpha}\widetilde h$ for all $\balpha$ with $|\balpha|_1\leq\lfloor\beta\rfloor.$ We now show that $\|h-\widetilde h\|_\infty \leq \delta.$ For any $\bx\in[-1,1]^r,$ let $\bx_i$ be the closest grid vertex, so that $|\bx-\bx_i|_\infty \leq \tau.$ Taylor expansion around $\bx_i$ gives
    \begin{align}
    \begin{split}\label{eq.main_part_remainder}
        (h - \widetilde h)(\bx) &= \sum_{\balpha:|\balpha|_1\leq \lfloor\beta\rfloor} \partial^{\balpha}(h-\widetilde h)(\bx_i) \frac{(\bx-\bx_i)^{\balpha}}{\balpha!} + R, \\
        R &= \sum_{\balpha:|\balpha|_1=\lfloor\beta\rfloor} \Big[\partial^{\balpha}(h-\widetilde h)(\bx_\xi)-\partial^{\balpha}(h-\widetilde h)(\bx_i)\Big] \frac{(\bx-\bx_i)^{\balpha}}{\balpha!},
    \end{split}
    \end{align}
    with $\bx_{\xi_i} = \bx_i + \xi_i(\bx-\bx_i)$ and a suitable $\xi_i\in[0,1].$ With $|\partial^{\balpha}(h-\widetilde h)|_{\beta-|\balpha|_1}$ the H\"older seminorm of the function $\partial^{\balpha}(h-\widetilde h) \in \mC_r^{\beta-|\balpha|_1}(2K),$ we bound the remainder $R$ by
    \begin{align*}
        |R| &\leq \sum_{\balpha:|\balpha|_1=\lfloor\beta\rfloor} \left| \partial^{\balpha}(h-\widetilde h)(\bx_\xi) - \partial^{\balpha}(h-\widetilde h)(\bx_i) \right| \frac{\tau^{|\balpha|_1}}{\balpha!} \\
        &\leq \sum_{\balpha:|\balpha|_1=\lfloor\beta\rfloor} \left| \partial^{\balpha}(h-\widetilde h)\right|_{\beta-|\balpha|_1} \tau^{\beta-|\balpha|_1} \frac{\tau^{|\balpha|_1}}{\balpha!} \\
        &\leq 2K\tau^\beta.
    \end{align*}
    Plugging this into the bound~\eqref{eq.main_part_remainder} gives
    \begin{align*}
        |(h - \widetilde h)(\bx)| &\leq \sum_{\balpha:|\balpha|_1\leq\lfloor\beta\rfloor} \left| \partial^{\balpha}(h-\widetilde h)(\bx_i) \right| \frac{\tau^{|\balpha|_1}}{\balpha!}  + 2K\tau^\beta \\
        &= \sum_{\balpha:|\balpha|_1\leq\lfloor\beta\rfloor} \tau^{\beta-|\balpha|_1} \left| \frac{\partial^{\balpha}(h-\widetilde h)(\bx_i)}{\tau^{\beta-|\balpha|_1}} \right| \frac{\tau^{|\balpha|_1}}{\balpha!}  + 2K\tau^\beta.
    \end{align*}
    In view of definition~\eqref{eq.def_A_alpha}, we denote $A^{\balpha}(h-\widetilde h)(\bx_i) = \lfloor \partial^{\balpha}(h-\widetilde h)(\bx_i)/\tau^{\beta-|\balpha|_1} \rfloor$ and $B^{\balpha}(h-\widetilde h)(\bx_i) = \partial^{\balpha}(h-\widetilde h)(\bx_i)/\tau^{\beta-|\balpha|_1} - A^{\balpha}(h-\widetilde h)(\bx_i).$ Thus $A^{\balpha}(h-\widetilde h)(\bx_i) = 0$ by assumption on $h,\widetilde h$ and $|B^{\balpha}(h-\widetilde h)(\bx_i)| < 1.$ We now prove 
    \begin{align}
    \sum_{\balpha:|\balpha|_1\leq\lfloor\beta\rfloor} \frac{1}{\balpha!} \leq e r^\beta.
    \label{eq.identity_567}
    \end{align}
    In fact, for any positive integer $k,$ consider the multinomial distribution induced by a fair $r$-sided die over $k$ independent rolls. The corresponding p.m.f. is $\balpha \mapsto r^{-k} k! /\balpha!$ and is supported on $\{\balpha:|\balpha|_1=k\}.$ Since the p.m.f. sums to one, we have $\sum_{\balpha:|\balpha|_1=k} 1/\balpha! = r^k/k!.$ By summing over $k=0,\ldots,\lfloor\beta\rfloor,$ one finds $\sum_{\balpha:|\balpha|_1\leq\lfloor\beta\rfloor} 1/\balpha! = \sum_{k=0}^{\lfloor\beta\rfloor} \sum_{\balpha:|\balpha|_1=k} 1/\balpha! \leq e r^\beta,$ proving \eqref{eq.identity_567}. Combining the discussion above together with the previous bounds, yields 
    \begin{align}
        |(h - \widetilde h)(\bx)| &\leq\tau^\beta \sum_{\balpha:|\balpha|_1\leq\lfloor\beta\rfloor} \frac{1}{\balpha!}  + 2K\tau^\beta
        \leq \tau^\beta \left(e r^\beta + 2K \right) 
        = \delta,
        \label{eq.890}
    \end{align}
    proving that, if two functions $h,\widetilde h \in \mC_r^\beta(K)$ have $A^{\balpha}h = A^{\balpha}\widetilde h$ for all $\balpha$ with $|\balpha|_1\leq\lfloor\beta\rfloor,$ then $\|h-\widetilde h\|_\infty \leq \delta.$
    
    The quantity $\mN(\delta,\mC_r^\beta(K),\|\cdot\|)$ is bounded above by cardinality $\#\mA$ of the set of matrices
    \begin{align*}
        \mA = \left\{ Ah = \left( A^{\balpha}h \right)_{\balpha:|\balpha|_1\leq \lfloor\beta\rfloor}^\top : h\in\mC_r^\beta(K) \right\}.
    \end{align*}
    The rows of the matrix $Ah$ consist of the row vectors $A^{\balpha}h.$ Since we consider $\balpha:|\balpha|_1\leq \lfloor\beta\rfloor,$ the matrix $Ah$ can have at most $(\lfloor\beta\rfloor+1)^r$ rows. Any matrix $Ah$ has moreover $M$ columns. 
    
    To complete the counting argument, we first explain the underlying idea. If two neighboring grid points $\bx_i,\bx_j$ are selected such that $|\bx_i-\bx_j|_\infty < 2\tau,$ say, then, $\partial^{\balpha} h(\bx_i) \approx \partial^{\balpha} h(\bx_j),$ whenever $|\balpha|_1<\beta.$ Since the $\ell$-th column of $Ah$ contains the discretized entries $(\partial^{\balpha} h(\bx_\ell))_{\balpha:|\balpha|_1\leq \lfloor\beta\rfloor},$ the number of possible realizations of the $i$-th and $j$-th column vector can be bounded by the possible realizations of the $i$-th column vector times a factor that describes the number of possible deviations of the values in the $j$-th column vector. 
    
    We now show that this factor is bounded by $2^{\beta+1}c(\beta).$ To see this, observe that Taylor expansion gives
    \begin{align*}
        \partial^{\balpha}h(\bx_j) &= \sum_{\bk:|\balpha|_1+|\bk|_1<\lfloor\beta\rfloor} \partial^{\balpha+\bk}h(\bx_i) \frac{(\bx_j-\bx_i)^{\bk}}{\bk!} + \sum_{\bk:|\balpha|_1+|\bk|_1=\lfloor\beta\rfloor} \partial^{\balpha+\bk}h(\bx_\xi) \frac{(\bx_j-\bx_i)^{\bk}}{\bk!},
    \end{align*}
    for some $\bx_\xi$ on the line with endpoints $\bx_i,\bx_j.$ By replacing $\wt h$ by $0$, $h$ by $\partial^{\balpha} h$, $\beta$ by $\beta-|\balpha|_1$ and $\tau$ by $2\tau,$ we can argue as for \eqref{eq.890} to find
    \begin{align*}
        \left| \partial^{\balpha}h(\bx_j) - \sum_{\bk:|\balpha|_1+|\bk|_1\leq\lfloor\beta\rfloor} \tau^{\beta-|\balpha|_1-|\bk|_1} A^{\balpha+\bk}h(\bx_i) \frac{(\bx_j-\bx_i)^{\bk}}{\bk!} \right|
        &\leq 2^{\beta-|\balpha|_1} c(\beta-|\balpha|_1) \tau^{\beta-|\balpha|_1} \\
        &\leq 2^{\beta} c(\beta) \tau^{\beta-|\balpha|_1}.
    \end{align*}
    This shows that, if the $i$-th column of $Ah$ is fixed, the values $\partial^{\balpha}h(\bx_j)$ range over an interval of length at most $2\cdot2^{\beta} c(\beta) \tau^{\beta-|\balpha|_1}.$ The entry $\lfloor \partial^{\balpha}h(\bx_j)/\tau^{\beta-|\balpha|_1}\rfloor$ can attain therefore at most  $2^{\beta+1} c(\beta) \tau^{\beta-|\balpha|_1}/\tau^{\beta-|\balpha|_1} +1 = 2^{\beta+1} c(\beta) + 1$ different values. As there are at most $(\beta+1)^r$ many rows, for fixed $i$-th column of $Ah,$ the $j$-th column of $Ah$ can attain at most $(2^{\beta+1} c(\beta)+1)^{(\beta+1)^r}$ different values.  
    
    Without loss of generality, assume that the points $\bx_1,\ldots,\bx_M$ are ordered in such a way that for each $j>1,$ there exists $i<j,$ such that $|\bx_i-\bx_j|_\infty < 2\tau.$ This determines then also the ordering of the columns of the matrix $Ah.$ In view of Equation~\eqref{eq.bound_on_coeff}, the first column of $Ah$ can attain at most  $(2K\tau^{-\beta} + 1)^{(\beta+1)^r}$ different values. For each of the $M-1$ remaining columns, we can use the argument above and find
    \begin{align*}
        \mN(\delta,\mC_r^\beta(K),\|\cdot\|_\infty) \leq \# \mA \leq (2K\tau^{-\beta} + 1)^{(\beta+1)^r} \cdot (2^{\beta+1} c(\beta)+1)^{(M-1)(\beta+1)^r}.
    \end{align*}
    Since $x+y\leq xy$ for all $x,y\geq 2,$ using $c(\beta) = er^\beta + 2K \leq 2eK r^\beta,$ the bound on $M$ in~\eqref{eq.bound_on_M} and the definition of $\tau = (\delta/c(\beta))^{1/\beta},$ the first assertion of the lemma follows. 
    
    For the bound on the constant $Q_1(\beta,r,K),$ we take the logarithm. With $\log(x+1) \leq \log (2x)$ for all $x>1,$ we get
    \begin{align*}
        \log&\left( \left( \frac{4eK^2 r^\beta}{\delta} + 1\right)^{(\beta+1)^r} \left(2^{\beta+2} e Kr^\beta+1\right)^{4^r (\beta+1)^r r^r (2eK)^{\frac{r}{\beta}} \delta^{-\frac{r}{\beta}}} \right) \\
        &\leq (\beta+1)^r \log\left( \frac{8eK^2r^\beta}{\delta}\right) + 4^r(\beta+1)^r r^r (2eK)^{\frac{r}{\beta}} \delta^{-\frac{r}{\beta}} \log\left(2^{\beta+3}eKr^\beta\right) \\
        & =: A_1 + A_2.
    \end{align*}
    Observe that $\log(x) < x^a/a$ for all $a,x>0,$ then
    \begin{align*}
        \log\left( \frac{8eK^2r^\beta}{\delta}\right) \leq \frac{\beta}{r} (8eK^2r^\beta)^{\frac{r}{\beta}} \delta^{-\frac{r}{\beta}} = \beta r^{r-1} (8eK^2)^{\frac{r}{\beta}} \delta^{-\frac{r}{\beta}},
    \end{align*}
    which yields $A_1 \leq (\beta+1)^{r+1} r^{r-1} (8eK^2)^{r/\beta} \delta^{-r/\beta}.$ Furthermore, using that $\log x < x$ for all $x>0,$
    \begin{align*}
        \log\left(2^{\beta+3}eKr^\beta\right) \leq 2(\beta+3) + \beta r + eK \leq (r+2)(\beta+3) + eK.
    \end{align*}
    Since $r>1,$ the latter display is smaller than $4(\beta+3)r+eK \leq 4eK(\beta+3)r$ and so $A_2 \leq 4eK(\beta+3)r\cdot 4^r(\beta+1)^r r^r (2eK)^{r/\beta} \delta^{-r/\beta}.$ Combining the bounds for $A_1$ and $A_2,$ we find
    \begin{align*}
        \frac{A_1 + A_2}{\delta^{-\frac{r}{\beta}}} &\leq (\beta+1)^{r+1} r^{r-1} (8eK^2)^{\frac{r}{\beta}} + eK(\beta+3) 4^{r+1} r^{r+1}(\beta+1)^r (2eK)^{\frac{r}{\beta}} \\
        &\leq (1+eK) 4^{r+1} (\beta+3)^{r+1} r^{r+1} (8eK^2)^{\frac{r}{\beta}},
    \end{align*}
    and the right hand side is $Q_1(\beta,r,K)$ by definition.
    
    We now prove the entropy bound in \eqref{eq.delta_n_entropy}. Let $Q_1 = Q_1(\beta,r,K),$ $C_1 = Q_1^{\beta\alpha/(2\beta\alpha+r)}$ and $\frakr_n = n^{-\beta\alpha/(2\beta\alpha+r)}.$ By construction, $\frakr_n^{-r/\beta\alpha} = n \frakr_n^2$ and $Q_1 C_1^{-r/\beta\alpha} = C_1^2.$ For any sequence $\delta_n \geq C_1 \frakr_n,$ the first part of the proof gives
    \begin{align*}
        \log\mN\left(\delta_n^{\frac{1}{\alpha}}, \mC_r^{\beta}(K),\|\cdot\|_\infty\right) &\leq \log\mN\left((C_1\frakr_n)^{\frac{1}{\alpha}}, \mC_{t_i}^{\beta_i}(K),\|\cdot\|_\infty\right) \\
        &\leq Q_1  C_1^{-\frac{r}{\beta\alpha}} \frakr_n^{-\frac{r}{\beta\alpha}} \\
        &= C_1^2 n \frakr_n^2 \\
        &\leq n \delta_n^2.
    \end{align*}
    The proof is complete.
\end{proof}

\begin{proof}[Proof of Lemma \ref{lem.local_entropy}]
   Fix any $\bbeta\in [\beta_{-},\beta_{+}]^{q+1}.$ In a first step we show that, with $R := 5K^q(q+1)$ and $\delta_{in}(\lambda,\bbeta) = \eps_n(\alpha_i,\beta_i,t_i)^{1/\alpha_i},$ we have
    \begin{align}\label{eq.Theta_entropy_first}
        \mN\left(R\eps_n(\lambda,\bbeta), \Theta_n(\lambda,\bbeta, K),\|\cdot\|_\infty\right) \leq \prod_{i=0}^q \mN\left(3\delta_{in}(\lambda,\bbeta), \Theta_{i,n}(\lambda,\bbeta, K),\|\cdot\|_\infty\right).
    \end{align}
     For any $i=0,\ldots,q,$ let $g_{i,1},\ldots,g_{i,N_i}$ be the centers of a $3\delta_{in}(\lambda,\bbeta)$-covering of $\Theta_{i,n}(\lambda,\bbeta, K).$ Then, any function $g_q\circ \dots \circ g_0 \in \Theta_n(\lambda,\bbeta, K)$ belongs to a ball around a composition of centers $g_{q,k_q}\circ\cdots\circ g_{0,k_0}$ for some $\bk=(k_0,\ldots,k_q)$ and such that $\|g_i-g_{i,k_i}\|_\infty \leq 3\delta_{in}(\lambda,\bbeta).$ By definition of $\Theta_{i,n}(\lambda,\bbeta, K),$ the components $(g_{ij,k_i})_j$ of $g_{i,k_i}$ satisfy $|g_{ij,k_i}(\bx)-g_{ij,k_i}(\by)| \leq 2\delta_{in}(\lambda,\bbeta) + K|\bx-\by|_\infty^{\beta_i\wedge 1},$ for all $\bx,\by\in[-1,1]^{d_i}.$ Using Lemma~\ref{lem.comp_approx}, the definition of $\delta_{in}(\lambda,\bbeta)$ and the fact that $\alpha_i \leq 1,$ gives
    \begin{align*}
        \left\|g_q\circ\cdots\circ g_0- g_{q,k_q}\circ\cdots\circ g_{0,k_0}\right\|_\infty 
        &\leq K^q \sum_{i=0}^q (2\delta_{in}(\lambda,\bbeta))^{\alpha_i} + \left(3\delta_{in}(\lambda,\bbeta)\right)^{\alpha_i} \\
        &\leq 5 K^q(q+1) \eps_n(\lambda,\bbeta) \\
        &= R\eps_n(\lambda,\bbeta).
    \end{align*}
    Since there are $N_0\times\cdots\times N_q$ centers, this concludes the first part of the proof.

    We now cover the set $\Theta_n(\lambda,B, K) = \cup_{\bbeta\in B} \Theta_n(\lambda,\bbeta, K).$ Denote by $\underline\bbeta$ the smallest element in $B$, that is, for any $\bbeta\in B$ we have $\underline\beta_i \leq \beta_i,$ for all $i=0,\ldots,q.$ Since $B$ is contained in the closed hypercube $[\beta_{-},\beta_{+}]^{q+1},$ we have $\underline\bbeta\in[\beta_{-},\beta_{+}]^{q+1}.$ We now show that, for any $i=0,\ldots,q$ and $\bbeta\in B,$ the set $\Theta_{i,n}(\lambda,\bbeta, K)$ is contained in the set $\Theta_{i,n}(\lambda,\underline\bbeta, K).$ In fact, by definition, any function $h_i\in \Theta_{i,n}(\lambda,\bbeta, K)$ satisfies $\overline h_{ij} \in \mD_i(\lambda,\bbeta,K)$ and, as a consequence of the embedding in Lemma~\ref{lem.Di_embedding} and the fact that $\delta_{in}(\lambda,\bbeta) \leq \delta_{in}(\lambda,\underline\bbeta)$ by the rate comparison condition~\ref{ass.rate_comparison} in Assumption~\ref{ass.rate}, one has $\mC_{t_i}^{\beta_i}(K) + \B_\infty(2\delta_{in}(\lambda,\bbeta)) \subseteq \mC_{t_i}^{\underline\beta_i}(K) + \B_\infty(2\delta_{in}(\lambda,\underline\bbeta)).$ Thus, $\mD_i(\lambda,\bbeta,K)\subseteq\mD_i(\lambda,\underline\bbeta,K)$ and $\Theta_{i,n}(\lambda,\bbeta, K)\subseteq\Theta_{i,n}(\lambda,\underline\bbeta, K).$ Together with~\eqref{eq.Theta_entropy_first}, we obtain
    \begin{align*}
        \mN\left(R\eps_n(\lambda,\underline\bbeta), \Theta_n(\lambda,B, K),\|\cdot\|_\infty\right) &\leq \mN\left(R\eps_n(\lambda,\underline\bbeta), \Theta_n(\lambda,\underline\bbeta, K),\|\cdot\|_\infty\right) \\
        &\leq \prod_{i=0}^q \mN\left(3\delta_{in}(\lambda,\underline\bbeta), \Theta_{i,n}(\lambda,\underline\bbeta, K),\|\cdot\|_\infty\right).
    \end{align*}
    We now use the definition of $\Theta_{i,n}(\lambda,\underline\bbeta, K)$ and upper bound the metric entropy by removing the constraint $\B_\infty(1)$ in the definition of $\mD_i(\lambda,\underline\bbeta,K).$ This gives
    \begin{align*}
        \mN\left(3\delta_{in}(\lambda,\underline\bbeta), \Theta_{i,n}(\lambda,\underline\bbeta, K),\|\cdot\|_\infty\right) 
        \leq \prod_{j=i}^{d_{i+1}} \mN\left(3\delta_{in}(\lambda,\underline\bbeta), \mC_{t_i}^{\underline\beta_i}(K) + \B_\infty\left(2\delta_{in}(\lambda,\underline\bbeta)\right),\|\cdot\|_\infty\right).
    \end{align*}
    Any function in $\mC_{t_i}^{\underline\beta_i}(K) + \B_\infty(2\delta_{in}\left(\lambda,\underline\bbeta)\right)$ is at most, in sup-norm distance, $2\delta_{in}(\lambda,\underline\bbeta)$-away from some function in $\mC_{t_i}^{\underline\beta_i}(K).$ Therefore, by applying Lemma~\ref{lem.Holder_entropy_sharp} with $r=t_i,$ $\beta=\underline\beta_i,$ $\alpha=\underline\alpha_i,$ and $\delta_n = \delta_{in}(\eta),$
    \begin{align*}
        \mN\left(3\delta_{in}(\lambda,\underline\bbeta), \mC_{t_i}^{\underline\beta_i}(K) + \B_\infty\left(2\delta_{in}(\lambda,\underline\bbeta)\right),\|\cdot\|_\infty\right)
        &\leq \mN\left(\delta_{in}(\lambda,\underline\bbeta), \mC_{t_i}^{\underline\beta_i}(K),\|\cdot\|_\infty\right) \\
        &\leq  e^{n\eps_n(\lambda,\underline\bbeta)^2}.
    \end{align*}
    Assumption~\ref{ass.rate} ensures that $\eps_n(\lambda,\underline\bbeta) \leq Q \eps_n(\lambda,\bbeta),$ thus combining the last inequalities gives
    \begin{align*}
        \log\mN\left(R Q \eps_n(\lambda,\bbeta), \Theta_n(\lambda,B, K),\|\cdot\|_\infty\right)
        &\leq 
        \log\mN\left(R\eps_n(\lambda,\underline\bbeta), \Theta_n(\lambda,B, K),\|\cdot\|_\infty\right) \\ &\leq \sum_{i=0}^q \sum_{j=1}^{d_{i+1}} n\eps_n(\lambda,\underline\bbeta)^2 \\
        &= |\bd|_1 n\eps_n(\lambda,\underline\bbeta)^2 \\
        &\leq Q^2 |\bd|_1  n\eps_n(\lambda,\bbeta)^2.
    \end{align*}
    Since $R = 5K^q(q+1),$ we use that $(q+1) \leq |\bd|_1\leq \log(2\log n),$ together with $\log(2\log n) \leq 2\log n.$ Thus, $K^q(q+1) \leq (2\log n)^{1+\log K}$ and, with $R_n = 5Q(2\log n)^{1+\log K},$
    \begin{align*}
        \log\mN\left(R_n\eps_n(\lambda,\bbeta), \Theta_n(\lambda,B, K),\|\cdot\|_\infty\right) &\leq \frac{R_n^2}{25} n\eps_n(\lambda,\bbeta)^2,
    \end{align*}
    which concludes the proof.
\end{proof}

\subsection{Proofs for Section~\ref{sec.dgp_minimax}} \label{sec.proofs_dgp_optimal}

\begin{proof}[Proof of Lemma~\ref{lem.smallest_solution}]
    Fix $\beta,r$ and let $\eps_T$ be such that $\varphi^{(\beta,r,K)}(\eps_T) \leq T \eps_T^2$ for all $T\geq 1.$ For this choice of $(\beta,r)$, we show that the concentration function inequality~\eqref{eq.def_eps_alpha} holds for any $0<\alpha\leq 1$ with $\eps_n(\alpha,\beta,r) := \eps_{m_n}^\alpha,$ where the sequence $m_n$ is chosen such that $m_n \eps_{m_n}^{2-2\alpha} \leq n.$ To see this, observe that 
    \begin{align*}
        \varphi^{(\beta,r,K)}\big(\eps_n(\alpha,\beta,r)^{1/\alpha}\big)
        &= \varphi^{(\beta,r,K)}\big(\eps_{m_n}\big)
        \leq m_n \eps_{m_n}^2 \leq n \eps_{m_n}^{2\alpha}
        =
        n \eps_n(\alpha,\beta,r)^2. 
    \end{align*}
    
    By Lemma~3 in~\cite{castillo2008lower}, the function $u \mapsto\varphi^{(\beta,r,K)}(u)$ is strictly decreasing on $u\in(0,+\infty),$ thus any $\overline\eps_n(\alpha,\beta,r) \geq \eps_n(\alpha,\beta,r)$ satisfies
    \begin{align*}
        \varphi^{(\beta,r,K)}\big(\overline\eps_n(\alpha,\beta,r)^{1/\alpha}\big) \leq \varphi^{(\beta,r,K)}\big(\eps_n(\alpha,\beta,r)^{1/\alpha}\big) \leq n \eps_n(\alpha,\beta,r)^2 \leq n \overline\eps_n(\alpha,\beta,r)^2,
    \end{align*}
    which concludes the proof.
\end{proof}

\begin{proof}[Proof of Lemma~\ref{lem.rates}]
{\it (i):} 
    By Lemma~\ref{lem.smallest_solution}, the sequence $\eps_n(\alpha,\beta,r)$ can be obtained from $\eps_n(1,\beta,r)$ via $\eps_n(\alpha,\beta,r) = \eps_{m_n}(1,\beta,r)^\alpha,$ for any sequence $m_n$ satisfying $m_n \eps_{m_n}(1,\beta,r)^{2-2\alpha} \leq n.$ We verify this for the sequence $m_n = C_3(\log n)^{-C_4} n^{(2\beta+r)/(2\beta\alpha+r)}$ with
    \begin{align*}
        C_3 &:= \left[ C_1(2\beta+1)^{C_2} \right]^{-\frac{(2-2\alpha)(2\beta+r)}{2\beta\alpha + r}}\quad \text{and} \ \ 
        C_4 := \frac{(2-2\alpha)(2\beta+r)}{2\beta\alpha+r} C_2.
    \end{align*}
    Since $C_1\geq 1$ and $n\geq 3,$ we must have $C_3\leq 1,$ $(\log n)^{-C_4}\leq 1$ and thus also $\log(m_n) \leq (2\beta+1)\log(n).$ Consequently,
    \begin{align*}
        m_n \eps_{m_n}(1,\beta,r)^{2-2\alpha} &= m_n \left[ C_1 (\log m_n)^{C_2} m_n^{-\frac{\beta}{2\beta+r}} \right]^{2-2\alpha} \\
        &= C_1^{2-2\alpha} (\log m_n)^{C_2(2-2\alpha)} m_n^{\frac{2\beta + r -2\beta + 2\beta\alpha}{2\beta+r}} \\
        &\leq C_1^{2-2\alpha} \big((2\beta+1)\log(n) \big)^{C_2(2-2\alpha)} m_n^{\frac{2\beta\alpha+r}{2\beta+r}} \\
        &\leq C_1^{2-2\alpha} (2\beta+1)^{C_2(2-2\alpha)} C_3^{\frac{2\beta\alpha+r}{2\beta+r}} (\log n)^{C_2(2-2\alpha) - C_4\frac{2\beta\alpha+r}{2\beta+r}} n \\
        &=n.
    \end{align*}
    As shown in the proof of Lemma~\ref{lem.smallest_solution}, any sequence $\overline \eps_n(\alpha,\beta,r) \geq \eps_n(\alpha,\beta,r)$ is still a solution to the concentration function inequality~\eqref{eq.def_eps_alpha}. We now derive a simple upper bound for $\eps_n(\alpha,\beta,r).$ Using that $\alpha \leq 1$ and $\log(m_n) \leq (2\beta+1)\log(n),$ we find
    \begin{align*}
        \eps_n(\alpha,\beta,r) &= \eps_{m_n}(1,\beta,r)^\alpha \\
        &\leq C_1^\alpha \log(m_n)^{\alpha C_2} m_n^{-\frac{\beta\alpha}{2\beta+r}} \\
        &\leq C_1 (2\beta+1)^{C_2} C_3^{-\frac{\beta\alpha}{2\beta+r}} (\log n)^{C_2+ C_4\frac{\beta\alpha}{2\beta\alpha+r}}
        n^{-\frac{\beta\alpha}{2\beta\alpha+r}}.
    \end{align*}
    Using the definition of $C_3$ together with $0<\alpha\leq1,$ $(2-2\alpha)\leq2$ and  $2\beta\alpha/(2\beta\alpha+r) \leq 1$ yields
    \begin{align*}
        C_3^{-\frac{\beta\alpha}{2\beta+r}} &= \left[ C_1(2\beta+1)^{C_2} \right]^{ \frac{(2-2\alpha)\beta\alpha}{2\beta\alpha + r}} \leq C_1(2\beta+1)^{C_2}.
    \end{align*}
    Similarly, we get
    \begin{align*}
        C_4\frac{\beta\alpha}{2\beta\alpha+r} &\leq C_2 \frac{(2-2\alpha)(2\beta+r)}{2\beta\alpha+r} \cdot \frac{1}{2} \leq  \frac{2\beta+r}{2\beta\alpha+r} \leq C_2(2\beta+1).
    \end{align*}
    The two previous displays recover the first assertion.

{\it (ii):}    
    By assumption, for any $\delta\in(0,1),$ we have $\varphi^{(\beta,r,K)}(\delta) \leq C_1' (\log \delta^{-1})^{C_2'} \delta^{-\frac{r}{\beta}}.$ We now choose $\delta = \eps_n(\alpha,\beta,r)^{1/\alpha}$ and $\eps_n(\alpha,\beta,r) = C_1' (\beta\log n)^{C_2'} n^{-\beta\alpha/(2\beta\alpha+r)}.$ Since $C_1' \geq 1,$ $C_2' \geq 0,$ and $\log n \geq 1,$
    \begin{align*}
        \log \eps_n(\alpha,\beta,r)^{-\frac{1}{\alpha}} &\leq - \frac{1}{\alpha}\log n^{-\frac{\beta\alpha}{2\beta\alpha+r}}
        \leq \frac{\beta}{2\beta\alpha+r} \log n
        \leq \beta \log n.
    \end{align*}
    Similarly, $$\eps_n(\alpha,\beta,r)^{-\frac{r}{\beta\alpha}}\leq \left(n^{-\frac{\beta\alpha}{2\beta\alpha+r}}\right)^{-\frac{r}{\beta\alpha}} = n \cdot n^{-\frac{2\beta\alpha}{2\beta\alpha+r}}\leq \frac{n \eps_n(\alpha,\beta,r)^2}{ (C_1')^2(\beta\log n)^{2C_2'}},$$
    and therefore,
    \begin{align*}
        \varphi^{(\beta,r,K)}\left( \eps_n(\alpha,\beta,r)^{\frac{1}{\alpha}} \right) &\leq C_1' \left(\log \eps_n(\alpha,\beta,r)^{-\frac{1}{\alpha}}\right)^{C_2'} \eps_n(\alpha,\beta,r)^{-\frac{r}{\beta\alpha}} \leq n\eps_n(\alpha,\beta,r)^2.
    \end{align*}
\end{proof}

\begin{proof}[Proof of Lemma~\ref{lem.eps_ratio}.] 
It is sufficient to show that for $n>1,$ any composition graph $\lambda = (q,\bd,\bt,\mS),$ and any $\bbeta'=(\beta_0',\dots,\beta_q'),\bbeta=(\beta_0,\dots,\beta_q)\in I(\lambda)$ satisfying $\beta_i'\leq \beta_i \leq \beta_i'+1/\log^2 n$ for all $i=0,\dots,q,$ the rates relative to the composition structures $\eta=(\lambda,\bbeta)$ and $\eta'=(\lambda,\bbeta')$ satisfy $\eps_n(\eta) \leq \eps_n(\eta') \leq e^{\beta_+} \eps_n(\eta).$ 

Since $\eps_n(\eta) = \widetilde C_1(\eta)(\log n)^{\widetilde C_2(\eta)} \frakr_n(\eta)$ with $\widetilde C_j(\eta):=\max_{i=0,\dots,q} \sup_{\beta \in [\beta_-,\beta_+]} C_j(\beta,t_i),$ $j\in \{1,2\},$ we have that $\widetilde C_j(\eta) = \widetilde C_j(\eta')$ and it is thus sufficient to prove $\frakr_n(\eta) \leq \frakr_n(\eta') \leq e^{\beta_+} \frakr_n(\eta).$
    
Using that $\frakr_n(\eta)=\max_{i=0,\ldots,q} n^{-\beta_i\alpha_i/(2\beta_i\alpha_i + t_i)}$ and the fact that the function $x\mapsto x/(2x+t_i)$ is strictly increasing for $x>0$ (its derivative is $x\mapsto t_i/(2x+t_i)^2$), the first inequality $\frakr_n(\eta) \leq \frakr_n(\eta')$ follows. For the second inequality, rewriting the expressions and simplifying the exponents gives
\begin{align*}
    \frac{\frakr_n(\eta')}{\frakr_n(\eta)}
    \leq \max_{i =0,\dots,q} \ \min_{j=0,\dots,q } \, n^{-\frac{\beta_i' \alpha_i'}{2\beta_i' \alpha_i'+t_i}+\frac{ \beta_j \alpha_j}{2 \beta_j \alpha_j+t_j}}
    \leq \max_{i =0,\dots,q} \, n^{-\frac{\beta_i' \alpha_i'}{2\beta_i' \alpha_i'+t_i}+\frac{ \beta_i \alpha_i}{2 \beta_i \alpha_i+t_i}}
    \leq \max_{i =0,\dots,q} \, n^{\beta_i\alpha_i-\beta_i'\alpha_i'}.
\end{align*}
We conclude the proof by showing that $|\beta_i\alpha_i-\beta_i'\alpha_i'| \leq \beta_{+}/\log n.$ For $u,u',v,v'\geq 0,$ we have that $|uv-u'v'|\leq u|v-v'|+v'|u-u'|.$ In particular, if $u,v'\leq 1,$ then also $|uv-u'v'|\leq |v-v'|+|u-u'|.$ By iterating this argument, we find that
\begin{align*}
    |\alpha_i-\alpha_i'|
    \leq \sum_{\ell=i+1}^q \big| (1\wedge \beta_\ell) - (1\wedge \beta_\ell')\big|
    \leq \sum_{\ell=i+1}^q \big|\beta_\ell
    -  \beta_\ell'\big|\leq \frac{q-i}{\log^2 n}
\end{align*}
and thus, $\beta_i\alpha_i-\beta_i'\alpha_i'\leq \beta_+|\alpha_i-\alpha_i'|+\alpha_i'|\beta_i-\beta_i'| \leq \beta_+ (q+1)/\log^2 n.$ Since we are restricting ourselves to graphs $\lambda$ such that $|\bd|_1 = 1 + \sum_{i=0}^{q} d_i \leq \log(2\log n),$ we must have $q + 1 \leq \log(2\log n).$ Since $\log x \leq x/2$ for all $x>0,$ we find $(q+1)\leq \log n.$ This gives $(q+1)/\log^2 n \leq 1/\log n$ and thus Assumption~\ref{ass.rate} (ii) holds with $Q = e^{\beta_{+}}.$ 
\end{proof}

\subsection{Proofs for Section~\ref{sec.examples}} \label{sec.proofs_ex}

\begin{proof}[Proof of Lemma~\ref{lem.fbm_rkhs}.]
    For the first part of the proof, take a kernel $\phi$ with $R_\beta := \int_{\R^r} |\bv|_\infty^{\beta} \phi(\bv) d\bv < +\infty.$ Using that $h\in \mC_r^{\beta}(K),$ $\beta\leq 1$ and the change of variable $\bv'=\bv/\sigma,$ we immediately get
    \begin{align*}
        \big|(h*\phi_\sigma)(\bu) - h(\bu)\big| &\leq  \int_{\R^r} \phi_\sigma(\bv)|h(\bu-\bv) - h(\bu)| \, d\bv \\
        &\leq K \int_{\R^r} |\bv|_\infty^{\beta} \phi_\sigma(\bv) \, d\bv \\
        &= K \int_{\R^r} |\bv|_\infty^{\beta} \sigma^{-r}\phi(\bv/\sigma) \, d\bv \\
        &= K \int_{\R^r} |\sigma\bv'|_\infty^{\beta} \phi(\bv') \, d\bv' \\
        &\leq K R_\beta \sigma^{\beta}.
    \end{align*}
    This shows $\|h*\phi_\sigma-h\|_\infty \leq K R_\beta \sigma^{\beta}$ and concludes the first part of the proof.
    
    We now deal with the RKHS norm. Notice that
	\begin{align*}
		(h*\phi_\sigma)(\bu) &= \int_{\R^r} \widehat{h}(\bxi)\widehat{\phi_\sigma}(\bxi) e^{-i\bu^\top\bxi} \frac{d\bxi}{(2\pi)^{r/2}} \\
		&= C_\beta^{1/2}\int_{\R^r} \widehat{h}(\bxi)\widehat{\phi_\sigma}(\bxi)  |\bxi|_2^{\beta+r/2} \frac{e^{-i\bu^\top\bxi}-1}{C_\beta^{1/2}|\bxi|_2^{\beta+r/2}} \frac{d\bxi}{(2\pi)^{r/2}} + (h*\phi_\sigma)(0).
	\end{align*}
	The RKHS of $Z+X^\beta$ is the direct sum of the space of constant functions $\bbH^Z$ and $\bbH^\beta.$ If the term $(h*\phi_\sigma)(0)$ is finite, it is a constant and thus belongs to $\bbH^Z.$ Then, the function $h*\phi_\sigma$ is a candidate element of $\bbH^Z\oplus\bbH^\beta$ since $h*\phi_\sigma-(h*\phi_\sigma)(0)$ has been represented as a potential element of the RKHS of $X^\beta.$ We now bound their norm using the isometry property of the norm $\|\cdot\|_{\bbH^\beta},$ so that
	\begin{align*}
		\|h*\phi_\sigma\|_{\bbH^Z\oplus\bbH^\beta}^2 \leq 2|(h*\phi_\sigma)(0)|^2 + 2C_\beta\int_{\R^r} |\widehat{h}(\bxi)|^2 |\widehat{\phi_\sigma}(\bxi)|^2  |\bxi|_2^{2\beta+r} \frac{d\bxi}{(2\pi)^{r/2}}.
	\end{align*}
	By the change of variable $\bxi'=\sigma\bxi,$ the fact that $\bxi \mapsto (1+|\bxi|_2)^{\beta}\widehat{h}(\bxi)$ has $L^2$-norm bounded by $K,$ the property $\widehat\phi_\sigma(\bxi) = \widehat\phi(\sigma\bxi),$ and choosing $\phi$ such that $M^2 := \sup_{\bxi\in\R^r} |\widehat{\phi}(\bxi)|^2|\bxi|_2^r < +\infty,$ we can bound
	\begin{align*}
		\int &|\widehat{h}(\bxi)|^2 |\widehat{\phi_\sigma}(\bxi)|^2  |\bxi|_2^{2\beta+r} \frac{d\bxi}{(2\pi)^{r/2}} \\
		&= \sigma^{-2\beta-2r} \int_{\R^r} |\widehat{h}(\bxi/\sigma)|^2 |\widehat{\phi}(\bxi)|^2  |\bxi|_2^{2\beta+r} \frac{d\bxi}{(2\pi)^{r/2}} \\
		&= \sigma^{-2\beta-2r} \int_{\R^r} \frac{(1+|\bxi/\sigma|_2)^{2\beta}}{(1+|\bxi/\sigma|_2)^{2\beta}} |\widehat{h}(\bxi/\sigma)|^2 |\widehat{\phi}(\bxi)|^2 |\bxi|_2^{2\beta+r} \frac{d\bxi}{(2\pi)^{r/2}} \\
		&\leq \sigma^{-2\beta-2r} \sup_{\bxi\in\R^r} \frac{|\widehat{\phi}(\bxi)|^2|\bxi|_2^{2\beta+r}}{(1+|\bxi/\sigma|_2)^{2\beta}} \int_{\R^r} (1+|\bxi/\sigma|_2)^{2\beta} |\widehat{h}(\bxi/\sigma)|^2 \frac{d\bxi}{(2\pi)^{r/2}} \\
		&\leq \sigma^{-2\beta-{2r}} \sup_{\bxi\in\R^r} \frac{|\widehat{\phi}(\bxi)|^2|\bxi|_2^{2\beta+r}}{|\bxi/\sigma|_2^{2\beta}} \int_{\R^r} \sigma^r (1+|\bxi|_2)^{2\beta} |\widehat{h}(\bxi)|^2 \frac{d\bxi}{(2\pi)^{r/2}} \\
		&\leq K^2 M^2 \sigma^{-r}.
	\end{align*}

	Similarly, by choosing $\phi$ such that $N^2 := (2\pi)^{-r/2} \int_{\R^r} |\widehat{\phi}(\bxi)|^2 \, d\bxi < +\infty,$ we obtain
	\begin{align*}
		|(h*\phi_\sigma)(0)|^2 &\leq \bigg( \int_{\R^r} |\widehat{h}(\bxi)| |\widehat{\phi_\sigma}(\bxi)| \frac{d\bxi}{(2\pi)^{r/2}} \bigg)^2 \\
		&= \bigg( \int_{\R^r} |\widehat{h}(\bxi)| (1+|\bxi|_2)^{\beta} \frac{|\widehat{\phi}(\sigma\bxi)|}{(1+|\bxi|_2)^{\beta}}  \frac{d\bxi}{(2\pi)^{r/2}} \bigg)^2 \\
		&\leq \bigg(\int_{\R^r} |\widehat{h}(\bxi)|^2 (1+|\bxi|_2)^{2\beta} \frac{d\bxi}{(2\pi)^{r/2}}\bigg) \bigg(\int_{\R^r} \frac{|\widehat{\phi}(\sigma\bxi)|^2}{(1+|\bxi|_2)^{2\beta}} \frac{d\bxi}{(2\pi)^{r/2}}\bigg) \\
		&\leq K^2 \sigma^{-r} \int_{\R^r} \frac{|\widehat{\phi}(\bxi)|^2}{(1+|\bxi/\sigma|_2)^{2\beta}} \frac{d\bxi}{(2\pi)^{r/2}} \\
		&\leq K^2 N^2 \sigma^{-r}.
	\end{align*}
	The proof is complete by taking $L_\beta^2 := 2(C_\beta + 1)(M^2 \vee N^2).$ Since this will be useful for the proof of Lemma~\ref{lem.fbm_ass2}, the explicit form of the constant $C_\beta$ is given in~(3.67) in~\cite{cohen2013fractional} as 
	\begin{align*}
	    C_\beta = \frac{\pi^{1/2} \Gamma(\beta + 1/2)}{2^{r/2} \beta \Gamma(2\beta) \sin(\pi\beta) \Gamma(\beta+r/2)},
	\end{align*}
	and depends only on $\beta,\ r.$
\end{proof}

\begin{proof}[Proof of Lemma~\ref{lem.fbm_ass2}]
We show that:
    \begin{compactenum}[(a)]
        \item Lemma~\ref{lem.rates}~(ii) holds for some $C_1'(\beta,r)\geq1$ and $C_2'(\beta,r)=0.$
        \item Any sequence $\eps_n(\alpha,\beta,r) \geq C_1'(\beta,r) n^{-\beta\alpha/(2\beta\alpha+r)}$ solves~\eqref{eq.def_eps_alpha}.
        \item  Assumption~\ref{ass.rate}~(i) holds for $\eps_n(\alpha,\beta,r) = C_1(\beta,r) n^{-\beta\alpha/(2\beta\alpha+r)}$ with $C_1(\beta,r) := C_1'(\beta,r) \vee Q_1(\beta,r,K)^{\beta/(2\beta+r)}.$
        \item $\sup_{\beta\in[\beta_{-},\beta_{+}]} C_1(\beta,r) <+\infty.$
        \item Assumption~\ref{ass.rate}~(ii) holds for $\eps_n(\eta)$ of the form~\eqref{eq.def_rate_enlarge}.
    \end{compactenum}
    By Lemma~\ref{lem.rates}, (a) $\implies$ (b) $\implies$ (c) and by Lemma~\ref{lem.eps_ratio}, (d) $\implies$ (e). Thus it remains to prove (a) and (d).
    
    \textit{Proof of (a):} Fix $\beta\in[\beta_{-},\beta_{+}].$ We have to show that, for all $\delta\in(0,1),$ $\varphi^{(\beta,r,K)}(\delta) \leq C_1'(\beta,r) \delta^{-r/\beta},$ for some constant $C_1'(\beta,r)\geq1$ depending only on $\beta,\ r,\ K.$ We denote the small-ball probability term by $\varphi_0^{(\beta,r)}(\delta) := -\log\P(\|Z+X^\beta\|_\infty \leq \delta)$ and show that
    \begin{align*}
        (A): 
        \ \ \sup_{\delta\in(0,1)} \delta^{r/\beta} \varphi_0^{(\beta,r)}(\delta) < +\infty,\quad (B): \ \  \sup_{\delta\in(0,1)} \delta^{r/\beta} \big( \varphi^{(\beta,r,K)}(\delta) - \varphi_0^{(\beta,r)}(\delta)\big) < +\infty.
    \end{align*}
    To prove (A), observe that the process $Z + X^\beta$ is the sum of two independent processes, thus its small-ball probability can be bounded by $\log\P(\|Z+X^\beta\|_\infty < \delta) \geq \log\P(\|Z\|_\infty < \delta/2) + \log\P(\|X^\beta\|_\infty < \delta/2).$ It is then sufficient to study the small-ball probabilities of $Z$ and $X^\beta$, separately. We now show the following condition, which implies (A),
    \begin{align}\label{eq.levy-fbm_c_Z}
        \sup_{\delta\in(0,1)} -\delta^{r/\beta} \log\P(\|X^\beta\|_\infty < \delta) < +\infty,\quad \sup_{\delta\in(0,1)} - \delta^{r/\beta} \log\P(|Z| < \delta) < +\infty.
    \end{align}
    Sharp bounds are known, see Theorem~5.1 in \cite{MR1861734}, for the small-ball probability of the fractional Brownian motion $X^{\beta}.$ In particular, for $0<\delta<1,$ we have $-\log \P(\|X^{\beta}\|_\infty\leq \delta) \leq c_X(\beta,r) \delta^{-r/\beta}$ for a finite constant $c_X(\beta,r)$ depending only on $\beta,\ r.$ Since $Z$ is a standard normal, we have $\P(|Z|\leq\delta) =(2\pi)^{-1/2} \int_{-\delta}^\delta e^{-x^2/2}dx \geq 2\delta e^{-\delta^2/2}/\sqrt{2\pi}.$ With the universal constant $c := 2/\sqrt{2\pi},$ this gives $-\log\P(|Z|\leq\delta) \leq \log(c^{-1}\delta^{-1}) + \delta^2/2.$ Therefore, using $\log(x) \leq x^a/a$ for all $x>1,a>0,$ we get
    \begin{align*}
        \sup_{\delta\in(0,1)} - \delta^{r/\beta} \log\P(|Z| < \delta) \leq \sup_{\delta\in(0,1)} \delta^{r/\beta} \Big(  \frac{\beta}{r} \left(\frac{1}{c\delta}\right)^{r/\beta} + \frac{\delta^2}{2}\Big) = \frac{\beta}{c^{r/\beta} r}  + \frac{1}{2} =: c_Z(\beta,r).
    \end{align*}
    This concludes the proof of (A).
    
    To prove (B), we apply Lemma~\ref{lem.fbm_rkhs}. In particular, with finite constants $R(\beta,r),L(\beta,r)$ depending only on $\beta,\ r,$ take $\sigma = (KR(\beta,r))^{-1/\beta} \delta^{1/\beta}.$ Then, any function $h\in\mC_r^\beta(K)\cap\mW_r^\beta(K)$ can be well approximated by the convolution $h*\phi_\sigma$ in such a way that $\|h-h*\phi_\sigma\|_\infty \leq \delta$ and $\|h*\phi_\sigma\|_{\bbH^Z\oplus\bbH^\beta}^2 \leq K^2L(\beta,r)^2 \delta^{-r/\beta}.$ This proves (B) because it gives 
    \begin{align*}
        \sup_{\delta\in(0,1)} \frac{\varphi^{(\beta,r,K)}(\delta) - \varphi_0^{(\beta,r)}(\delta)}{\delta^{-r/\beta}} \leq \sup_{\delta\in(0,1)} \frac{K^2L(\beta,r)^2 \delta^{-r/\beta}}{\delta^{-r/\beta}} = K^2L(\beta,r)^2.
    \end{align*}
    
    We have thus concluded the proof of (a), that is, for any $\beta\in[\beta_{-},\beta_{+}],$ condition~(ii) in Lemma~\ref{lem.rates} holds with finite constants
    \begin{align*}
        C_1'(\beta,r) := c_X(\beta,r) + c_Z(\beta,r) + K^2 L(\beta,r)^2 ,\quad C_2'(\beta,r)=0.
    \end{align*}
    
    \textit{Proof of (d):} With the definition of $C_1(\beta,r),$ we want to show that
    \begin{align*}
        \sup_{\beta\in[\beta_{-},\beta_{+}]} C_1'(\beta,r) \vee Q_1(\beta,r,K)^{\beta/(2\beta+r)} < +\infty.
    \end{align*}
    The constant $Q_1(\beta,r,K)$ is given explicitly in Lemma~\ref{lem.Holder_entropy_sharp} and depends continuously on $\beta>0.$ Thus, $\sup_{\beta\in[\beta_{-},\beta_{+}]} Q_1(\beta,r,K) =: \widetilde Q_1 <+\infty.$ Since the function $\beta\mapsto \beta/(2\beta+r)$ is increasing for $\beta>0,$ we also have $Q_1(\beta,r,K)^{\beta/(2\beta+r)} \leq \widetilde Q_1^{\beta_{+}/(2\beta_{+}+r)}.$ 
    
    In the previous part of the proof we have found $C_1'(\beta,r) = c_Z(\beta,r) + c_X(\beta,r) + K^2 L(\beta,r)^2,$ thus it remains to prove
    \begin{align}\label{eq.levy-fbm_to_show_4}
        \sup_{\beta\in[\beta_{-},\beta_{+}]} c_Z(\beta,r) + c_X(\beta,r) + K^2 L(\beta,r)^2 <+\infty.
    \end{align}
    By examining the proof of Lemma~\ref{lem.fbm_rkhs}, we know that $\sup_{\beta\in[\beta_{-},\beta_{+}]} K^2 L(\beta,r)^2 < +\infty.$ The explicit form of $c_Z(\beta)$ is given in~\eqref{eq.levy-fbm_c_Z} and so, with $c=2/\sqrt{2\pi},$
    \begin{align*}
        \sup_{\beta\in[\beta_{-},\beta_{+}]} c_Z(\beta,r) \leq \frac{\beta_{+}}{r c^{\beta_{+}} } + \frac{1}{2} <+\infty.
    \end{align*}
     We now show that the properties of $c_X(\beta,r)$ can be deduced from Theorem~5.2 in~\cite{MR1861734}. We observe that $\E[|X^\beta(\bu)-X^\beta(\bu')|^2] = |\bu-\bu'|_2^{2\beta}.$ Furthermore, the function $\beta \mapsto \E[X^\beta(\bu)X^\beta(\bu')]$ is continuous for all fixed $\bu,\bu'\in[-1,1]^r.$ In the notation of Theorem~5.2 in~\cite{MR1861734}, we can take $\sigma_\beta(\delta) := \delta^\beta$ and check that, with $c_1:=1/2^{\beta_{+}}$ and $c_2 := 1,$ $c_1 \sigma_\beta(2\delta\wedge1) \leq \sigma_\beta(\delta) \leq c_2 \sigma_\beta(2\delta\wedge1)$ for all $\delta\in(0,1).$ The constants $c_1,\ c_2$ are chosen to be independent of $\beta$ in the compact interval $[\beta_{-},\beta_{+}].$ From this, one obtains a constant $c_X(r)$ that only depends on $c_1,c_2$ and such that $-\log\P(\|X^\beta\|_\infty < \delta) \leq c_X(r)\delta^{-r/\beta}$ for all $\delta>0.$ This shows that the quantity $c_X(\beta,r)$ in~\eqref{eq.levy-fbm_to_show_4} can be replaced by $c_X(r)$ and thus is bounded, concluding the proof of (d).
\end{proof}

\begin{proof}[Proof of Lemma \ref{lem.gp_conditioned}]
    Using the definition of $X^\beta$ yields
    \begin{align*}
        \|X^\beta\|_{\infty,\infty,\beta} = \max_{j=1,\ldots,J_\beta} \frac{1}{\sqrt{jr}} \max_{k = 1,\ldots,2^{jr}} |Z_{j,k}|.
    \end{align*}
    It is known that $\E[\max_{k = 1,\ldots,2^{jr}} Z_{j,k}] \leq \sqrt{2\log(2^{jr})},$ a reference is Lemma~2.3 in~\cite{MR2319879}. For $K'>1,$ using symmetry of $Z_{j,k}$ and the Borell-TIS inequality, e.g. Theorem~2.1.1 in~\cite{MR2319516},
    \begin{align*}
        \P\left( \max_{k = 1,\ldots,2^{jr}} |Z_{j,k}| \geq (1+K') \sqrt{2\log(2^{jr})}\right) &\leq 2 \P\left( \max_{k = 1,\ldots,2^{jr}} Z_{j,k} \geq (1+K') \sqrt{2\log(2^{jr})}\right) \\
        &\leq 4 \exp\left(-{K'}^2 \log(2^{jr}) \right).
    \end{align*}
    Combining this with the union bound and the formula for the geometric sum, we obtain for any $K'>2/\sqrt{r},$ 
    \begin{align*}
        \P\left(\exists j=1,\ldots,J_\beta, \max_{k = 1,\ldots,2^{jr}} |Z_{j,k}| \geq (1+K') \sqrt{2\log(2^{jr})}\right) \leq \sum_{j=1}^{J_\beta} 2^{j(2-r{K'}^2)} \leq \frac{1}{1-2^{2-r{K'}^2}} - 1.
    \end{align*}
    Therefore, with $K'>\sqrt{3},$ on an event with probability at least $1-4/(2^{r{K'}^2}-4)$, we find
    \begin{align*}
        \|X^\beta\|_{\infty,\infty,\beta} \leq \max_{j=1,\ldots,J_\beta} \frac{(1+K') \sqrt{2\log(2^{jr})}}{\sqrt{jr}} = (1+K') \sqrt{2\log2}.
    \end{align*}
\end{proof}

\begin{proof}[Proof of Lemma \ref{lem.trunc_conditioned_ass2}] 
    In view of Section~4.3.6 in~\cite{gine2016mathematical}, the Besov space $\mB_{\infty,\infty,\beta}$ contains the H\"older space $\mC_r^\beta$ for any $\beta>0,$ and they coincide whenever $\beta\notin\N.$ Thus there exists $K'$ such that $\mC_r^\beta(K)\subseteq \mB_{\infty,\infty,\beta}(K').$ We show that:
    \begin{compactenum}[(a)]
        \item The assumption of Lemma~\ref{lem.rates}~(i) holds for some $C_1'(\beta,r)\geq1$ and $C_2'(\beta,r)=3/2$. 
        \item Any sequence $\eps_n(\alpha,\beta,r) \geq C_1'(\beta,r)^2(2\beta+1)^3 (\log n)^{3(\beta+1)} n^{-\beta\alpha/(2\beta\alpha+r)}$ solves the concentration function inequality~\eqref{eq.def_eps_alpha}.
        \item Assumption~\ref{ass.rate}~(i) holds by taking $\eps_n(\alpha,\beta,r) = C_1(\beta,r) (\log n)^{C_2(\beta,r)} n^{-\beta\alpha/(2\beta\alpha+r)}$ with $C_1(\beta,r) := C_1'(\beta,r)^2(2\beta+1)^3 \vee Q_1(\beta,r,K)^{\beta/(2\beta+r)}$ and $C_2(\beta,r) := 3(\beta+1).$
        \item $\sup_{\beta\in[\beta_{-},\beta_{+}]} C_1(\beta,r) <+\infty$ and $\sup_{\beta\in[\beta_{-},\beta_{+}]} C_2(\beta,r) <+\infty.$
        \item Assumption~\ref{ass.rate}~(ii) holds for an $\eps_n(\eta)$ of the form~\eqref{eq.def_rate_enlarge}.
    \end{compactenum}
    By Lemma~\ref{lem.rates}, (a) $\implies$ (b) $\implies$ (c) and by Lemma~\ref{lem.eps_ratio}, (d) $\implies$ (e). Thus it remains to prove (a) and (d).
    
    \textit{Proof of (a):}  We denote the small-ball probability term by $\varphi_0^{(\beta,r)}(\delta) := -\log\P(\|X^\beta\|_\infty \leq \delta)$ and the RKHS term by $\varphi^{(\beta,r,K)}(\delta) - \varphi_0^{(\beta,r)}(\delta).$ We start with the RKHS term. The proof of Theorem~4.5 in~\cite{vvvz} shows that any function $h\in\mB_{\infty,\infty,\beta}(K')$ can be well approximated by its projection $h^{J_\beta}$ at truncation level $J_\beta.$ In fact, one has $\|h-h^{J_\beta}\|_\infty \leq K' 2^{-J_\beta \beta} /(2^\beta-1) $ and, with coefficients $\omega_j = 2^{-j(\beta+r/2)}/\sqrt{jr},$
    \begin{align*}
        \|h^{J_\beta}\|_{\bbH^\beta}^2 = \sum_{j=1}^{J_\beta} \sum_{k=1}^{2^{jr}} \lambda_{j,k}(h)^2 \omega_j^{-2} \leq {K'}^2 r J_\beta \sum_{j=1}^{J_\beta} 2^{jr} \leq {K'}^2 r J_\beta^2 2^{J_\beta r}.
    \end{align*}
    Recall that $J_\beta$ is defined as the closest integer to the solution $J$ of $2^{J} = n^{1/(2\beta+r)}.$ By definition, we always have $J_\beta \leq 1 + \log_2 n/(2\beta+r)$ and so $2^{J_\beta } \leq 2 n^{1/(2\beta+r)}$ and $2^{-J_\beta \beta} \geq 2^{-\beta} n^{-\beta/(2\beta+r)}.$ With all the above, the choice
    \begin{align*}
        \delta_n := K' \frac{(2^\beta+1)^2}{(2^\beta-1)} \sqrt{r 2^r} J_\beta^{3/2} 2^{-J_\beta \beta},
    \end{align*}
    implies $\|h-h^{J_\beta}\|_\infty < \delta_n$ and
    \begin{align*}
        \varphi^{(\beta,r,K)}(\delta_n) - \varphi_0^{(\beta,r)}(\delta_n) &\leq {K'}^2 r J_\beta^2 2^{J_\beta r} \\
        &\leq {K'}^2 r 2^r J_\beta^2  n^{\frac{r}{2\beta+r}} \\
        &\leq {K'}^2 \frac{(2^\beta+1)^2}{(2^\beta-1)^{2}} r 2^r J_\beta^2 n^{\frac{r}{2\beta+r}} \\
        &\leq n {K'}^2 \frac{(2^\beta+1)^4}{(2^\beta-1)^{2}} r 2^r J_\beta^2 2^{-2J_\beta \beta} \\
        &\leq n\delta_n^2.
    \end{align*}
    
    We now study the small-ball probability. The proof of Theorem~4.5 in~\cite{vvvz} shows that, for any sequence $\delta_n\in(0,1),$
    \begin{align*}
        \varphi_0^{(\beta,r)}(\delta_n) \leq - \sum_{j=1}^{J_\beta} 2^{jr} \log\left( 2\Phi\left(\frac{\delta_n 2^{j\beta}}{\widetilde K(\beta) + j^2 r^2} \right) - 1\right),
    \end{align*}
    where $\widetilde K(\beta)$ is chosen in such a way that the function $x\mapsto x^{\beta/r}/(\widetilde K(\beta) + \log_2^2(x))$ is increasing for $x\geq1.$ Taking the derivative and imposing it to be positive for $x>1$ yields $\beta \widetilde K(\beta)/r + \log_2^2(x) - 2\log(x)/\log_2^2(2) > 0,$ which is solved for any $\widetilde K(\beta) \geq 4r/\beta.$ \cite{vvvz} also shows that the function $f(y) = -\log(2\Phi(y)-1)$ is decreasing and can be bounded above by $f(y) \leq 1+ |\log y|$ on any interval $y\in[0,c].$ Thus we find
    \begin{align*}
        \varphi_0^{(\beta,r)}(\delta_n) &\leq  \sum_{j=1}^{J_\beta} 2^{jr} \left( 1 + \left|\log\left( \frac{\delta_n 2^{j\beta}}{\widetilde K(\beta) + j^2 r^2} \right) \right| \right).
    \end{align*}
    We now show that for sufficiently large $n,$ we have $\delta_n \leq (\widetilde K(\beta) + J_\beta^2 r^2) 2^{-J_\beta \beta}.$ For this, one has to verify the inequality
    \begin{align*}
        {K'} \frac{(2^\beta+1)^2}{(2^\beta-1)} \sqrt{r 2^r} J_\beta^{3/2} \leq \widetilde K(\beta) + J_\beta^2 r^2,
    \end{align*}
    which holds for sufficiently large $n$ since $J_\beta^{3/2} \ll J_\beta^2.$ Therefore, using that the function $j\mapsto 1 + \log( (\widetilde K(\beta) + j^2 r^2)/(\delta_n 2^{j\beta}))$ is decreasing, together with $1+\log(x)\leq2\log(x)$ for all $x\geq e$, and $\log(a/b)=\log(a)+\log(1/b)$,
    \begin{align*}
        \varphi_0^{(\beta,r)}(\delta_n) &\leq  \sum_{j=1}^{J_\beta} 2^{jr} \left( 1 + \log\left( \frac{\widetilde K(\beta) + j^2 r^2}{\delta_n 2^{j\beta}} \right) \right) \\
        &\leq J_\beta 2^{J_\beta r} \left( 1 + \log\left( \frac{\widetilde K(\beta) + r^2}{\delta_n 2^{\beta}} \right) \right) \\
        &\leq 2 J_\beta 2^{J_\beta r} \left[ \log\left( \frac{\widetilde K(\beta) + r^2}{2^{\beta}} \right) + \log\left( \frac{1}{\delta_n} \right) \right].
    \end{align*}
    Then, with the definition of $\delta_n,$ 
    \begin{align*}
        \varphi_0^{(\beta,r)}(\delta_n) &\leq 2 \left[ \log\left( \frac{\widetilde K(\beta) + r^2}{2^{\beta}} \right) + \log\left( \frac{1}{\delta_n} \right) \right] J_\beta 2^{J_\beta r} \\ 
        &\leq 2 \left[ \log\left( \frac{\widetilde K(\beta) + r^2}{2^{\beta}} \right) + \log\left( \frac{2^{J_\beta \beta}}{{K'} \frac{(2^\beta+1)^2}{(2^\beta-1)} \sqrt{r 2^r J_\beta^3} } \right) \right] J_\beta 2^{J_\beta r} \\
        &\leq 2 \left[ \log\left( \frac{\widetilde K(\beta) + r^2}{2^{\beta}} \right) + \log\left( 2^{J_\beta \beta} \right)\right] J_\beta 2^{J_\beta r} \\
        &\leq 2 \left[ \log\left( \frac{\widetilde K(\beta) + r^2}{2^{\beta}} \right) + \beta\log\left( 2 \right)\right] J_\beta^2 2^{J_\beta r} \\
        &\leq 2 \left[ \log\left( \frac{\widetilde K(\beta) + r^2}{2^{\beta}} \right) + \beta\log\left( 2 \right)\right] 2^r J_\beta^2 n^{\frac{r}{2\beta+r}}.
    \end{align*}
    The last term in the latter display only depends on $n$ through the quantity $J_\beta^2 n^{\frac{r}{2\beta+r}}$. Since $J_\beta^2 \ll J_\beta^3,$ for large enough $n$ the latter display is smaller than ${K'}^2 \frac{(2^\beta+1)^2}{(2^\beta-1)^{2}} 2^r J_\beta^3 n^{\frac{r}{2\beta+r}} \leq n\delta_n^2.$
    
    This concludes the proof of (a), since we have shown that the sequence
    \begin{align*}
        \eps_n(1,\beta,r) := {K'} \frac{(2^\beta+1)^2}{(2^\beta-1)} \sqrt{r 2^r} J_\beta^{3/2} 2^{-J_\beta \beta},
    \end{align*}
    solves the concentration function inequality~\eqref{eq.def_eps_alpha} for $\alpha=1.$
    
    \textit{Proof of (d):} 
    The constant $Q_1(\beta,r,K)$ is given explicitly in Lemma~\ref{lem.Holder_entropy_sharp} and depends continuously on $\beta>0,$ so $\sup_{\beta\in[\beta_{-},\beta_{+}]} Q_1(\beta,r,K) =: \widetilde Q_1 <+\infty.$ Since the function $\beta\mapsto \beta/(2\beta+r)$ is increasing for $\beta>0,$ we also have $Q_1(\beta,r,K)^{\beta/(2\beta+r)} \leq \widetilde Q_1^{\beta_{+}/(2\beta_{+}+r)}.$ 
    
    All the quantities involved in the construction of the rate $\eps_n(1,\beta,r)$ are explicit. It is immediate to see that they are all bounded on the compact interval $\beta\in[\beta_{-},\beta_{+}]$ since $\beta_{-}>0.$
    
    This concludes the proof of (d).
\end{proof}

\begin{proof}[Proof of Lemma \ref{lem.stat_ass2}]
    We show that:
    \begin{compactenum}[(a)]
        \item Lemma~\ref{lem.rates}~(i) holds for some $C_1'(\beta,r)\geq1$ and $C_2'(\beta,r)=(1+r)\beta/(2\beta+r).$
        \item Any sequence $\eps_n(\alpha,\beta,r) \geq C_1'(\beta,r)^2(2\beta+1)^{2C_2'(\beta,r)} (\log n)^{(2\beta+2)C_2'(\beta,r)} n^{-\beta\alpha/(2\beta\alpha+r)}$ solves the concentration function inequality~\eqref{eq.def_eps_alpha}.
        \item Assumption~\ref{ass.rate}~(i) holds for any $\eps_n(\alpha,\beta,r) = C_1(\beta,r) (\log n)^{C_2(\beta,r)} n^{-\beta\alpha/(2\beta\alpha+r)}$ with $C_1(\beta,r) := C_1'(\beta,r)^2(2\beta+1)^{2C_2'(\beta,r)} \vee Q_1(\beta,r,K)^{\beta/(2\beta+r)}$ and $C_2(\beta,r) := (2\beta+2)(1+r)\beta/(2\beta+r).$
        \item $\sup_{\beta\in[\beta_{-},\beta_{+}]} C_1(\beta,r) <+\infty$ and $\sup_{\beta\in[\beta_{-},\beta_{+}]} C_2(\beta,r) <+\infty.$
        \item Assumption~\ref{ass.rate}~(ii) holds for an $\eps_n(\eta)$ of the form~\eqref{eq.def_rate_enlarge}.
    \end{compactenum}
    By Lemma~\ref{lem.rates}, (a) $\implies$ (b) $\implies$ (c) and by Lemma~\ref{lem.eps_ratio}, (d) $\implies$ (e). Thus it remains to prove (a) and (d).

    \textit{Proof of (a):} Let $\varphi_a^{(\beta,r,K)}$ the concentration function of the rescaled process $X^{\nu}(a\cdot),$ then Lemma~11.55 and Lemma~11.56 in \cite{GVDV17} show that, for all $0<\delta<1,$
    \begin{align*}
        \varphi_a^{(\beta,r,K)}(\delta) \leq \Big(C(r)  \big(\log(a\delta^{-1})\big)^{1+r} + D(r)\Big) a^r,
    \end{align*}
    where $C(r)$ and $D(r)$ are constants that only depend on $r$ and the spectral measure $\nu$ of $X^{\nu}.$ It is sufficient to solve the concentration function inequality~\eqref{eq.def_eps_alpha} for $\alpha=1.$ The solution is given in Section~11.5.2 in~\cite{GVDV17} as $\eps_n(1,\beta,r) = C_1'(\beta,r) (\log n)^{(1+r)\beta/(2\beta+r)} n^{-\beta/(2\beta+r)},$ for some constant $C_1'(\beta,r)\geq1$ depending on $\beta,\ r,\ K.$ 
    
    \textit{Proof of (d):} 
    The constant $Q_1(\beta,r,K)$ is given explicitly in Lemma~\ref{lem.Holder_entropy_sharp} and depends continuously on $\beta>0,$ so $\sup_{\beta\in[\beta_{-},\beta_{+}]} Q_1(\beta,r,K) =: \widetilde Q_1 <+\infty.$ Since the function $\beta\mapsto \beta/(2\beta+r)$ is increasing for $\beta>0,$ we also have $Q_1(\beta,r,K)^{\beta/(2\beta+r)} \leq \widetilde Q_1^{\beta_{+}/(2\beta_{+}+r)}.$ 
    
    The dependence on $\beta$ in the concentration function bound only appears in the scaling $a = a(\beta,r),$ since the constants $C(r)$ and $D(r)$ are independent of $\beta.$ The right side of the latter display depends continuously on the scaling $a = a(\beta,r),$ which in turn is continuous in $\beta\in[\beta_{-},\beta_{+}]$ by construction~\eqref{eq.scaling}. This gives
    \begin{align*}
        \sup_{\beta\in\beta\in[\beta_{-},\beta_{+}]} C_1(\beta,r) < +\infty,\quad \sup_{\beta\in\beta\in[\beta_{-},\beta_{+}]} C_2(\beta,r) \leq (2\beta_{+}+1)(1+r)\frac{\beta_{+}}{2\beta_{+}+r}.
    \end{align*}
\end{proof}

\subsection{Proofs for Section~\ref{sec.disc_dgp_prior}} \label{sec.proofs_disc}
\begin{proof}[Proof of Lemma \ref{lem.redundant}]
For two functions $g_k \in \mC_1^{\beta_k}(1),$ $k=1,2,$ and $\beta_1, \beta_2 \leq 1,$ we have that $|g_2(g_1(x))-g_2(g_1(y))|\leq |g_1(x)-g_1(y)|^{\beta_2}\leq |x-y|^{\beta_1\beta_2}.$ Hence, $g_2\circ g_1 \in \mC_1^{\beta_1\beta_2}(1).$ We now write 
\begin{align*}
    f= h_q\circ \dots \circ h_0
    = h_q \circ \dots \circ h_{j+1} \circ \wt h_j \circ h_{j-2}\circ \dots \circ h_0,
\end{align*}
with $\wt h_j:=h_j \circ h_{j-1}.$ The right hand side can be written as composition structure $\eta':= (q-1 ,\bd_{-j} ,\bt_{-j} ,\mS_{-j},\bbeta'),$ with $\bd_{-j} ,\bt_{-j} ,\mS_{-j},\bbeta'$ as defined in the statement of the lemma. Due to $\beta_+\leq 1,$ we have $\mathfrak{r}_n(\eta)=\max_{i=0,\dots,q} n^{-\frac{\gamma_i}{2\gamma_i+t_i}},$ with $\gamma_i = \prod_{\ell=i}^q \beta_\ell$ and it follows that $\mathfrak{r}_n(\eta)=\mathfrak{r}_n(\eta').$
\end{proof}

\newpage 

\bibliographystyle{acm}       % (uses file "plain.bst")
\bibliography{biblio}           % expects file "refsPart1.bib"

\end{document}